\documentclass[a4paper,12pt,reqno]{amsart}
\usepackage{amsfonts,amsmath,amssymb,amstext,amsxtra}
\usepackage{mathabx}
\usepackage{euscript}
\usepackage{inputenc}
\usepackage{color,graphicx,psfrag,xcolor}
\usepackage{enumitem}
\usepackage{polski}
\usepackage[polish,english]{babel}
\usepackage{pgf}
\usepackage[T1]{fontenc}

\setlist[itemize]{leftmargin=*}

\newtheorem{theorem}{Theorem}[section]
\newtheorem{proposition}[theorem]{Proposition}
\newtheorem{lemma}[theorem]{Lemma}
\newtheorem{corollary}[theorem]{Corollary}
\newtheorem{remark}[theorem]{Remark}
\newtheorem{definition}[theorem]{Definition}
\newtheorem{example}[theorem]{Example}
\newtheorem{problem}[theorem]{Problem}

\newcommand\1{1\hspace{-1mm}\text{\rm I}}

\newcommand\Ad{\operatorname{Ad}}
\newcommand\apartment{\mathfrak{a}}
\newcommand\artanh{\operatorname{artanh}}

\newcommand\C{\mathbb{C}}
\newcommand\chamber{\apartment^+}

\newcommand\co{\operatorname{co}}
\newcommand\const{\operatorname{const.}}
\newcommand\Cosh{\operatorname{Cosh}}
\renewcommand\cot{\operatorname{cot}}
\renewcommand\coth{\operatorname{coth}}
\newcommand\crat{c_{\ssf\text{\rm rat}}}
\newcommand\ctrig{c_{\ssf\text{\rm trig}}}
\renewcommand\epsilon{\varepsilon}
\newcommand\Hn{\mathbb{H}^{\ssf n}}
\newcommand\Id{\operatorname{I}}
\renewcommand\Im{\operatorname{Im}}
\newcommand\N{\mathbb{N}}
\renewcommand\phi{\varphi}
\newcommand\Q{\mathbb{Q}}

\renewcommand\Re{\operatorname{Re}}
\newcommand\R{\mathbb{R}}
\newcommand\Rn{\mathbb{R}^n}
\newcommand\Rnminusone{\mathbb{R}^{\ssf n-1}}
\newcommand\Rnplusone{\mathbb{R}^{1+\ssf n}}
\newcommand\sign{\operatorname{sign}}
\newcommand\Sn{\mathbb{S}^{\ssf n}}
\newcommand\Snminusone{\mathbb{S}^{\ssf n-1}}
\newcommand\ssb{\hspace{-.25mm}}
\newcommand\ssf{\hspace{.25mm}}

\renewcommand\tanh{\operatorname{tanh}}
\newcommand\Tq{\mathbb{T}_q}
\newcommand\vol{\operatorname{vol}}
\newcommand\vsb{\hspace{-.1mm}}
\newcommand\vsf{\hspace{.1mm}}
\newcommand\Z{\mathbb{Z}}

\title[Dunkl theory]
{An introduction to Dunkl theory\\
and its analytic aspects}
\author[Anker]
{Jean--Philippe Anker}
\address{Universit\'e d'Orl\'eans \& CNRS,
F\'ed\'eration Denis Poisson (FR 2964),
Laboratoire MAPMO (UMR 7349),
B\^atiment de Math\'ematiques,
B.P.~6759, 45067 Orl\'eans cedex 2, France}
\email{anker@univ-orleans.fr}
\thanks{Work partially supported by the regional project MADACA
(Marches Al\'eatoires et processus de Dunkl\;--\;Approches
Combinatoires et Alg\'ebriques, www.fdpoisson.fr/madaca).
We thank the organizers of the school and conference AAGADE 2015
(\textit{Analytic, Algebraic and Geometric Aspects of Differential Equations\/},
Mathematical Research and Conference Center, B\k{e}dlewo, Poland, September 2015)
for their invitation and the warm welcome at Pa{\l}ac B\k{e}dlewo.
We also thank
B\'echir Amri,
Nizar Demni,
L\'eonard Gallardo,
Chaabane Rejeb,
Margit R\"osler,
Simon Ruijsenaars,
Patrice Sawyer
and Michael Voit
for helpful comments and discussions.}
\begin{document}

${}$
\vspace{-5mm}

\maketitle
\vspace{-3mm}

\centerline{
To appear in
}\centerline{\it
Analytic, Algebraic and Geometric Aspects of Differential Equations\/}
\centerline{
(\textit{\foreignlanguage{polish}{B"edlewo}, Poland, September 2015\/}),
}\centerline{
G.~Filipuk, Y.~Haraoka \& S.~Michalik (eds.),
}\centerline{
Trends Math., Birkh\"auser
}

\tableofcontents

\section{Introduction}

Dunkl theory is a far reaching generalization
of Fourier analysis and special function theory
related to root systems.
During the sixties and seventies,
it became gradually clear that
radial Fourier analysis on rank one symmetric spaces
was closely connected
with certain classes of special functions in one variable\,:
\begin{itemize}

\item
Bessel functions in connection with radial Fourier analysis on Euclidean spaces,

\item
Jacobi polynomials
in connection with radial Fourier analysis on spheres,

\item
Jacobi functions
(i.e.~the Gauss hypergeometric function
\ssf${}_2\ssf\text{\rm F}_{\ssb1}$)
in connection with radial Fourier analysis on hyperbolic spaces.

\end{itemize}
See \cite{Koornwinder1984} for a survey.
During the eighties, several attempts were made,
mainly by the Dutch school (Koornwinder, Heckman, Opdam),
to extend these results in higher rank (i.e.~in several variables),
until the discovery of Dunkl operators in the rational case
and Cherednik operators in the trigonometric case.
Together with $q$--special functions introduced by Macdonald,
this has led to a beautiful theory, developed by several authors
which encompasses in a unified way
harmonic analysis on all Riemannian symmetric spaces
and spherical functions thereon\,:
\begin{itemize}

\item
generalized Bessel functions on flat symmetric spaces,
and their asymmetric version, known as the Dunkl kernel,

\item
Heckman--Opdam hypergeometric functions
on positively or negatively curved symmetric spaces,
and their asymmetric version, due to Opdam,

\item
Macdonald polynomials on affine buildings.

\end{itemize}

Beside Fourier analysis and special functions,
this theory has also deep and fruitful interactions with
\begin{itemize}

\item
algebra (double affine Hecke algebras),

\item
mathematical physics
(Calogero-Moser-Sutherland models,
quantum many body problems),

\item
probability theory (Feller processes with jumps).

\end{itemize}

There are already several surveys about Dunkl theory available in the literature\,:
\begin{itemize}

\item
\cite{Roesler2002} (see also \cite{DunklXu2001})
about rational Dunkl theory (state of the art in 2002),

\item
\cite{Opdam2000} about trigonometric Dunkl theory (state of the art in 1998),

\item
\cite{Etingof2007} about integrable systems related to Dunkl theory,

\item
\cite{Macdonald2003} and \cite{Cherednik2005}
about \ssf$q$\hspace{.4mm}--\ssf Dunkl theory and affine Hecke algebras,

\item
\cite{Graczyk2008} about probabilistic aspects of Dunkl theory
(state of the art in 2006).

\end{itemize}

\begin{figure}[ht]
\psfrag{affine buildings}[c]{\tiny$\text{affine buildings}$}
\psfrag{Bessel functions (p)}[c]{\tiny$\text{Bessel functions}$}
\psfrag{Bessel functions (Rn)}[c]{\color{red}\tiny$\text{Bessel functions}$}
\psfrag{circular}[c]{\footnotesize$\text{circular}\vphantom{|}$}
\psfrag{compact}[c]{\color{blue}\footnotesize$\text{compact}$}
\psfrag{DAHA}[c]{$\text{DAHA}$}
\psfrag{double affine}[c]{\tiny$\text{double affine}$}
\psfrag{Euclidean}[c]{\color{blue}\footnotesize$\text{Euclidean}$}
\psfrag{Dunkl theory (rational)}[c]{\color{red}\footnotesize$\text{Dunkl theory}$}
\psfrag{Dunkl theory (trigonometric)}[c]{\color{red}\footnotesize$\text{Dunkl theory}$}
\psfrag{GK}[c]{\tiny$G/K$}
\psfrag{generalized}[c]{\tiny$\text{generalized}$}
\psfrag{Hecke algebras}[c]{\tiny$\text{Hecke algebras}$}
\psfrag{higher rank}[c]{\color{blue}\footnotesize$\text{higher rank}$}
\psfrag{hyperbolic}[c]{\color{red}\footnotesize$\text{hyperbolic}\vphantom{|}$}
\psfrag{hyperbolic spaces}[c]{\color{red}\tiny$\text{hyperbolic spaces}$}
\psfrag{Jacobi functions}[c]{\color{red}\tiny$\text{Jacobi functions}$}
\psfrag{Jacobi polynomials}[c]{\color{red}\tiny$\text{Jacobi polynomials}$}
\psfrag{Macdonald (buildings)}[c]{\tiny$\text{Macdonald}$}
\psfrag{Macdonald (q)}[c]{\footnotesize$\text{Macdonald}$}
\psfrag{non compact}[c]{\color{blue}\footnotesize$\text{non compact}$}
\psfrag{p}[c]{\tiny$\mathfrak{p}$}
\psfrag{p-adic}[c]{\color{blue}\footnotesize$\text{$p$\ssf--\ssf adic}$}
\psfrag{q-polynomials}[c]{\footnotesize$\text{$q$\ssf--\ssf polynomials}$}
\psfrag{Rn}[c]{\color{red}\tiny$\R^n$}
\psfrag{rank 1}[c]{\color{blue}\footnotesize$\text{rank 1}$}
\psfrag{rational}[c]{\color{red}\footnotesize$\text{rational}$}
\psfrag{spheres}[c]{\color{red}\tiny$\text{spheres}$}
\psfrag{spherical functions (GK)}[c]{\tiny$\text{spherical functions}$}
\psfrag{spherical functions (UK)}[c]{\tiny$\text{spherical functions}$}
\psfrag{spherical functions (buildings)}[c]{\tiny$\text{spherical functions}$}
\psfrag{spherical functions (tree)}[c]{\color{red}\tiny$\text{spherical functions}$}
\psfrag{theory}[c]{\footnotesize$\text{theory}$}
\psfrag{trees}[c]{\color{red}\tiny$\text{homogeneous trees}$}
\psfrag{trigonometric}[c]{\color{red}\footnotesize$\text{trigonometric}$}
\psfrag{UK}[c]{\tiny$U/K$}
\centerline{\hspace{26mm}\includegraphics[width=165mm]{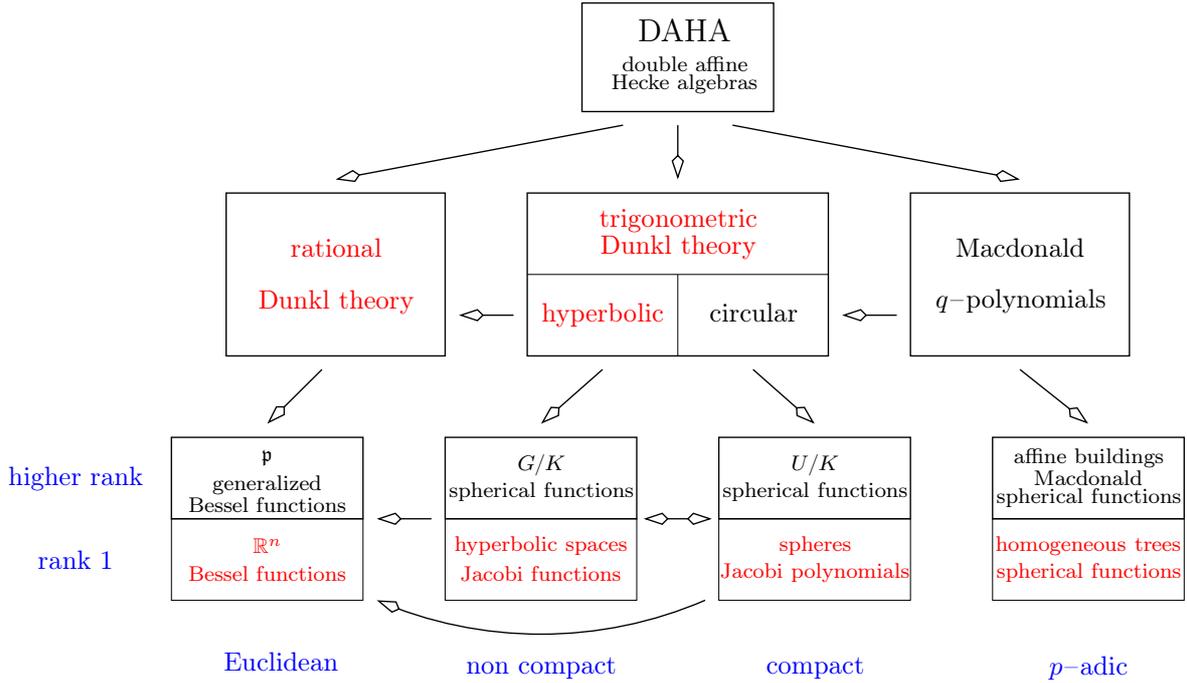}}
\label{BigPicture}
\caption{Relation between various special function theories}
\end{figure}

These lectures are intended to give
an overview of some analytic aspects of Dunkl theory.
The topics are indicated in red in Figure \ref{BigPicture},
where we have tried to summarize relations
between several theories of special functions,
which were alluded to above,
and where arrows mean limits.

\smallskip

Let us describe the content of our notes.
In Section 2, we consider several geometric settings
(Euclidean spaces, spheres, hyperbolic spaces, homogeneous trees, \dots)
where radial Fourier analysis is available and can be applied successfully,
for instance to study evolutions equations
(heat equation, wave equation, Schr\"odinger equation, \dots ).
Section 3 is devoted to the rational Dunkl theory
and Section 4 to the trigonometric Dunkl theory.
In both cases, we first review the basics
and next address some important analytic issues.
We conclude with an appendix about root systems
and with a comprehensive bibliography.
For lack of time and competence,
we haven't touched upon other aforementioned aspects of Dunkl theory,
for which we refer to the bibliography.

\section{Spherical Fourier analysis in rank $1$}

\subsection{Cosine transform}
Let us start with an elementary example.
Within the framework of even functions on the real line \ssf$\R$\ssf,
the Fourier transform is given by
\begin{equation*}\label{CosineTransform}
\widehat{f}\ssf(\lambda)
=\!\int_{\ssf\R}\ssb dx\,f(x)\,\cos\lambda\ssf x
\end{equation*}
and the inverse Fourier transform by
\begin{equation*}\label{InverseCosineTransform}
f(x)=\tfrac1{2\ssf\pi}\int_{\ssf\R}\ssb d\lambda\,
\widehat{f}\ssf(\lambda)\,\cos\lambda\ssf x\,.
\end{equation*}
The cosine functions
\ssf$\phi_\lambda(x)\ssb=\ssb\cos\lambda\ssf x$ ($\lambda\!\in\!\C$)
occurring in these expressions can be characterized in various ways.
Let us mention
\begin{itemize}

\item
\textit{Power series expansion}\,:
\vspace{-2mm}
\begin{equation*}
\phi_\lambda(x)=\sum\nolimits_{\ssf\ell=0}^{\ssf+\infty}
\tfrac{(-1)^\ell}{(2\ssf\ell)\ssf!}\,(\lambda\ssf x)^{2\ssf\ell}
\qquad\forall\;\lambda,x\!\in\!\C\,.
\end{equation*}
\vspace{-5mm}

\item
\textit{Differential equation}\,:
the functions \ssf$\phi\ssb=\ssb\phi_\lambda$
are the smooth eigenfunctions of \ssf$(\frac\partial{\partial x})^2$,
which are even and normalized by  \ssf$\phi(0)\!=\!1$\ssf.

\item
\textit{Functional equation}\,:
the functions \ssf$\phi\ssb=\ssb\phi_\lambda$
\ssf are the nonzero continuous functions on~\ssf$\R$
\ssf which satisfy
\vspace{-1mm}
\begin{equation*}
\tfrac{\phi(x\ssf+\ssf y)\,+\,\phi(x\ssf-\ssf y)\vphantom{\frac00}}
{2\vphantom{\frac00}}
=\phi(x)\,\phi(y)
\qquad\forall\;x,y\!\in\!\R\ssf.
\end{equation*}

\end{itemize}

\subsection{Hankel transform on Euclidean spaces}
\label{Hankel}

The Fourier transform on \ssf$\Rn$ and its inverse are given by
\vspace{-1mm}
\begin{equation}\label{FourierTransformRn}
\widehat{f}\ssf(\xi)=\ssb{\displaystyle\int_{\ssf\Rn}}\!dx\,
f(x)\,e^{-\ssf i\ssf\langle\ssf\xi,\ssf x\ssf\rangle}
\end{equation}
\vspace{-5mm}

and
\vspace{-1mm}
\begin{equation}\label{InverseFourierTransformRn}
f(x)=(2\pi)^{-n}{\displaystyle\int_{\ssf\Rn}}\!d\lambda\,
\widehat{f}\ssf(\xi)\,e^{\,i\ssf\langle\ssf\xi,\ssf x\ssf\rangle}
\end{equation}
Notice that the Fourier transform of
a radial function \ssf$f\!=\!f(r)$ \ssf on \ssf$\Rn$
is again a radial function \ssf$\widehat{f}\!=\!\widehat{f}\ssf(\lambda)$\ssf.
In this case,
\eqref{FourierTransformRn} and \eqref{InverseFourierTransformRn} become
\begin{equation}\label{HankelTransform}
\widehat{f}\ssf(\lambda)
=\tfrac{2\,\pi^{\frac n2}}{\Gamma(\frac n2)}
\int_{\,0}^{+\infty}\hspace{-1.5mm}dr\;
r^{\ssf n-1}\,f(r)\;j_{\frac {n-2}2}(i\ssf\lambda\hspace{.4mm}r)
\end{equation}
\vspace{-5mm}

and
\vspace{-1mm}
\begin{equation}\label{InverseHankelTransform}
f(r)=\tfrac{1\vphantom{\frac00}}
{2^{\ssf n-1\vphantom{\frac00}}\,\pi^{\frac n2}\,\Gamma(\frac n2)}
\int_{\,0}^{+\infty}\hspace{-1.5mm}d\lambda\;
\lambda^{\ssf n-1}\,\widehat{f}\ssf(\lambda)\;
j_{\frac {n-2}2}(i\ssf\lambda\hspace{.4mm}r)\,.
\end{equation}

Instead of the exponential function or the cosine function,
\eqref{HankelTransform} and \eqref{InverseHankelTransform} involve now
the modified Bessel function \ssf$j_{\frac {n-2}2}$\ssf,
which can be characterized again in various ways\,:
\begin{itemize}

\item
\textit{Relation with classical special functions and power series expansion}\ssf.
For every \ssf$z\!\in\!\C$\ssf,
\begin{align*}
j_{\frac {n-2}2}(z)&=\Gamma(\tfrac n2)\,
\bigl(\tfrac{i\ssf z}2\bigr)^{\!\frac{2-n}2}\ssf
J_{\frac{n-2}2}(i\ssf z)\\
&=\sum\nolimits_{\ssf\ell=0}^{\ssf+\infty}\,
\tfrac{\Gamma(\frac n2)\vphantom{\frac0|}}
{\ell\,!\;\Gamma(\frac n2+\ssf\ell\ssf)\vphantom{\frac|0}}\,
\bigl(\tfrac z2\bigr)^{\ssb2\ssf\ell}\\
&={}_0\ssf\text{\rm F}_{\ssb1}\vsb(\tfrac n2\ssf;\ssb\tfrac{z^2}4)
=e^{-z}\ssf{}_1\vsf\text{\rm F}_{\ssb1}\vsb
(\tfrac{n-1}2\ssf;\ssb n\hspace{-.5mm}-\!1\ssf;\ssb2\ssf z)\hspace{.5mm},
\end{align*}
where $J_\nu$ denotes the classical Bessel function of the first kind and
\begin{equation*}
{}_p\vsf\text{\rm F}_{\hspace{-.6mm}q}
(a_1,\dots,\ssb a_{\vsf p}\ssf;
\ssb b_{\vsf1},\dots,\ssb b_{\vsf q}\ssf;\ssb z)
={\displaystyle\sum\nolimits_{\ssf\ell=0}^{+\infty}}\,
\tfrac{(a_1\ssb)_\ell\ssf\dots\ssf(a_{\vsf p}\ssb)_\ell}
{(b_1\ssb)_\ell\ssf\dots\ssf(b_{\vsf q}\ssb)_\ell}\,
\tfrac{z^{\ssf\ell}}{\ell\,!}
\end{equation*}
the generalized hypergeometric function.
\item
\textit{Differential equations}\ssf.
The function
\,$\phi_\lambda(r)\ssb=\ssb\smash{j_{\frac {n-2}2}}(i\ssf\lambda\hspace{.4mm}r)$
\ssf is the unique smooth solution to the differential equation
\begin{equation*}\label{BesselEquation}
\bigl(\tfrac\partial{\partial\ssf r}\bigr)^{\ssb2}\ssf\phi_\lambda
+\tfrac{n-1}r\,\bigl(\tfrac\partial{\partial\ssf r}\bigr)\,\phi_\lambda
+\lambda^2\,\phi_\lambda=0\,,
\end{equation*}
which is normalized by \ssf$\phi_\lambda(0)\!=\!1$\ssf.
Equivalently, the function
\begin{equation}\label{BesselRn}
x\longmapsto\phi_\lambda(|x|)=j_{\frac {n-2}2}(i\ssf\lambda\ssf|x|)
\end{equation}
is the unique smooth radial normalized eigenfunction of the Euclidean Laplacian
\begin{equation*}
\Delta_{\ssf\Rn}
=\ssf\sum\nolimits_{\ssf j=1}^{\,n}
\bigl(\tfrac\partial{\partial\ssf x_j}\bigr)^{\ssb2}\ssb
=\ssf\bigl(\tfrac\partial{\partial\ssf r}\bigr)^{\ssb2}\!
+\tfrac{n-1}r\ssf\bigl(\tfrac\partial{\partial\ssf r}\bigr)
+\tfrac1{r^2}\,\Delta_{\ssf\Snminusone}
\end{equation*}
corresponding to the eigenvalue \ssf$-\lambda^2$.

\end{itemize}

\begin{remark}
The function \eqref{BesselRn}
\ssf is a matrix coefficient of a continuous unitary representation
of the Euclidean motion group \,$\Rn\!\rtimes\ssb\text{O}(n)$\ssf.
\end{remark}

The function  \eqref{BesselRn} is a spherical average of plane waves.
Specifically,
\begin{equation*}
\phi_\lambda(|x|)\ssf
=\ssb\int_{\text{O}(n)}\hspace{-1mm}dk\;
e^{\,i\ssf\lambda\ssf\langle\ssf u,\ssf k.x\ssf\rangle}
=\ssf\tfrac{\Gamma(\frac n2)\vphantom{\frac0|}}
{2\,\pi^{\frac n2\vphantom{\frac00}}}
\int_{\ssf\Snminusone}\hspace{-1mm}dv\;
e^{\,i\ssf\lambda\ssf\langle\ssf v,\ssf x\ssf\rangle}\,,
\end{equation*}
where $u$ is any unit vector in \ssf$\Rn$.
Hence the integral representation
\begin{equation}\begin{aligned}\label{SphericalFunctionsRn}
\phi_\lambda(r)\ssf
&=\ssf\tfrac{\Gamma(\frac n2)\vphantom{\frac0|}}
{\sqrt{\pi\ssf}\;\Gamma(\frac{n-1}2)\vphantom{\frac|0}}
\int_{\,0}^{\ssf\pi}\hspace{-1mm}d\ssf\theta\,
(\sin\theta\ssf)^{n-2}\,e^{\,i\ssf\lambda\,r\cos\theta}\\
&=\ssf\tfrac{2\;\Gamma(\frac n2)\vphantom{\frac0|}}
{\sqrt{\pi\ssf}\;\Gamma(\frac{n-1}2)\vphantom{\frac|0}}
\;r^{\ssf2-n}\ssb\int_{\,0}^{\ssf r}\!ds\;
(\ssf r^2\!-\ssb s^2\ssf)^{\ssb\frac{n-3}2}\ssf\cos\lambda\ssf s\,.
\end{aligned}\end{equation}

\subsection{Spherical Fourier analysis on real spheres}

\begin{figure}[ht]
\psfrag{0}[l]{$0$}
\psfrag{e0}[l]{$e_{\ssf0}$}
\psfrag{x0}[l]{$x_{\ssf0}$}
\psfrag{Rn}[c]{$\Rn$}
\psfrag{Korbit}[l]{\color{red}$K$\!--\ssf orbit}
\includegraphics[height=70mm]{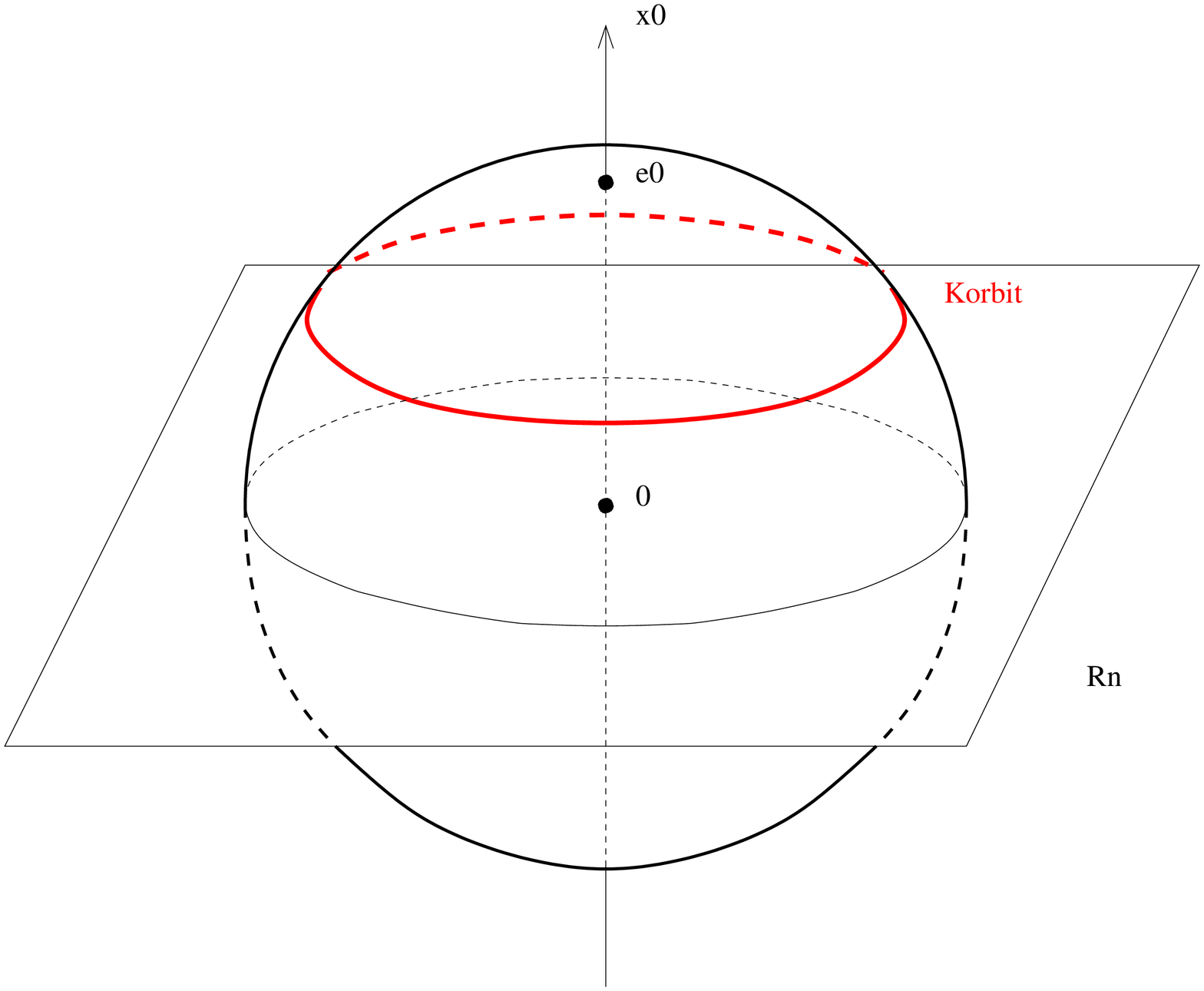}
\caption{Real sphere \ssf$\Sn$}
\end{figure}

Real spheres
\begin{equation*}
\Sn\ssb=\{\,x\ssb=\ssb(x_0,x_1,\dots,x_n)\ssb\in\ssb\Rnplusone\,|\,
|x|^2\ssb=\ssb x_0^2\ssb+\ssf\dots\ssf+\ssb x_n^2\ssb=\ssb1\,\}
\end{equation*}
of dimension \ssf$n\ssb\ge\ssb2$
\ssf are the simplest examples
of Riemannian symmetric spaces of compact type.
They are simply connected Riemannian manifolds,
with constant positive sectional curvature.
The Riemannian structure on \ssf$\Sn$
is induced by the Euclidean metric in \ssf$\Rnplusone$,
restricted to the tangent bundle of \ssf$\Sn$,
and the Laplacian on \ssf$\Sn$ is given by \,$\Delta f\ssb
=\ssb\widetilde{\Delta}\ssf\widetilde{f}\,\big|_{\ssf\Sn}$\ssf,
where \ssf$\vphantom{\Big|}\widetilde{\Delta}\ssb
=\!\smash{\displaystyle\sum\nolimits_{\ssf j=0}^{\,n}}
\bigl(\tfrac{\partial}{\partial\ssf x_j}\bigr)^{\ssb2}$
\ssf denotes the Euclidean Laplacian in \ssf$\Rnplusone$
and \ssf$\smash{\widetilde{f}(x)\ssb=\ssb f\bigl(\tfrac x{|x|}\bigr)}$
\ssf the homogeneous extension of \ssf$f$
\ssf to \ssf$\Rnplusone\!\smallsetminus\!\{0\}$\ssf.
In spherical coordinates
\begin{equation*}\begin{cases}
\,x_0&\hspace{-3.5mm}=\cos\theta_1\ssf,\\
\,x_1&\hspace{-3.5mm}=\sin\theta_1\cos\theta_2\ssf,\\
&\hspace{-2.5mm}\vdots\\
\,x_{\ssf n-1}&\hspace{-3.5mm}=
\sin\theta_1\sin\theta_2\,\dots\,\sin\theta_{n-1}\cos\theta_n\ssf,\\
\,x_{\ssf n}&\hspace{-3.5mm}=
\sin\theta_1\sin\theta_2\,\dots\,\sin\theta_{n-1}\sin\theta_n\ssf,\\
\end{cases}\end{equation*}
the Riemannian metric,
the Riemannian volume
and the Laplacian
read respectively
\begin{align*}
ds^2\ssb
&=\sum\nolimits_{\ssf j=1}^{\,n}
(\ssf\sin\theta_1)^2\dots\ssf(\ssf\sin\theta_{j-1})^2\,(d\ssf\theta_j)^2
\ssf,\\
d\ssb\vol\ssb
&=(\ssf\sin\theta_1)^{n-1}\dots\ssf(\ssf\sin\theta_{n-1})
d\ssf\theta_1\dots d\ssf\theta_n
\end{align*}
\vspace{-5.5mm}

\noindent
and
\vspace{-.5mm}
\begin{equation*}
\Delta=\sum\nolimits_{\ssf j=1}^{\,n}
\tfrac{1\vphantom{\frac00}}
{(\ssf\sin\theta_1)^2\ssf\dots\,(\ssf\sin\theta_{j-1})^2\vphantom{\frac00}}\,
\bigl\{\bigl(\tfrac\partial{\partial\ssf\theta_j}\bigr)^{\ssb2}\!
+(\vsf n\hspace{-.4mm}-\hspace{-.4mm}j\vsf)\ssf
(\ssf\cot\theta_j)\,\tfrac\partial{\partial\ssf\theta_j}\bigr\}\,.
\end{equation*}
Let \ssf$G\ssb=\ssb\text{O}(n\hspace{-.4mm}+\!1)$
be the isometry group of \ssf$\Sn$
and let \ssf$K\ssb\approx\ssb\text{O}(n)$
be the stabilizer of \ssf$e_{\ssf0}\ssb=\ssb(1,0,\dots,0)$\ssf.
Then \ssf$\Sn$ can be realized as the homogeneous space \ssf$G/K$.
As usual, we identify right\ssf--\ssf$K$\ssb--\ssf invariant functions on \ssf$G$
\ssf with functions on \ssf$\Sn$,
and bi\ssf--\ssf$K$\ssb--\ssf invariant functions on \ssf$G$
\ssf with radial functions on \ssf$\Sn$ 
i.e.~functions on \ssf$\Sn$
which depend only on \ssf$x_{\ssf0}\ssb=\ssb\cos\theta_1$\ssf.
For such functions,
\begin{equation*}\label{IntegrationSn}
\int_{\ssf\Sn}\!d\vol\,f\ssf
=\,2\,\tfrac{\pi^{\frac n2}\vphantom{\frac00}}
{\Gamma(\frac n2)\vphantom{\frac|0}}
\int_{\,0}^{\ssf\pi}\ssb d\theta_1\,
(\ssf\sin\theta_1)^{n-1}\ssf f(\cos\theta_1)
\end{equation*}
\vspace{-4mm}

\noindent
and
\vspace{-.5mm}
\begin{equation*}\label{RadialLaplacianSn}
\Delta\ssf f=\tfrac{\partial^{\ssf2}\ssb f}{\partial\ssf\theta_1^{\ssf2}}
+(n\hspace{-.4mm}-\!1)\ssf(\ssf\cot\theta_1)\,
\tfrac{\partial f}{\partial\ssf\theta_1}\,.
\end{equation*}
The spherical functions on \ssf$\Sn$ are
the smooth normalized radial eigenfunctions of the Laplacian on \ssf$\Sn$.
Specifically,
\begin{equation*}\begin{cases}
\;\Delta\,\phi_{\ssf\ell}
=-\ssf\ell\,(\ell\hspace{-.4mm}+\hspace{-.4mm}n\hspace{-.4mm}-\!1)\,
\phi_{\ssf\ell}\,,\\
\;\phi_{\ssf\ell\ssf}(e_{\ssf0})=1\,,
\end{cases}\end{equation*}
where \ssf$\ell\!\in\!\N$\ssf.
They can be expressed in terms of classical special functions, namely
\begin{equation*}
\phi_{\ssf\ell\ssf}(x_{\ssf0})
=\tfrac{\ell\,!\,(n-2)\ssf!}{(\ell\ssf+\ssf n-2)\ssf!}\,
C_{\ssf\ell}^{\ssf(\frac{n-1}2)}(x_{\ssf0})
=\tfrac{\ell\,!}{(\frac n2)_\ell}\,
P_{\ssf\ell}^{\ssf(\frac n2-1,\frac n2-1)}(x_{\ssf0})
\end{equation*}
\vspace{-6mm}

or
\vspace{-1mm}
\begin{equation*}
\phi_{\ssf\ell\ssf}(\cos\theta_1)
={}_2\ssf\text{\rm F}_{\ssb1}\bigl(-\ssf\ell,
\ell\hspace{-.4mm}+\hspace{-.4mm}n\hspace{-.4mm}-\!1\ssf;
\ssb\tfrac n2\ssf;\sin^2\ssb\theta_1\bigr)\ssf,
\end{equation*}
where \ssf$C_{\ssf\ell}^{\ssf(\lambda)}$
are the Gegenbauer or ultraspherical polynomials,
\ssf$P_{\ssf\ell}^{(\alpha,\beta)}$ the Jacobi polynomials
and ${}_2\ssf\text{\rm F}_{\ssb1}$ the Gauss hypergeometric function.

\begin{remark}
We have emphasized the characterization
of spherical functions on \,$\Sn$
by a differential equation.
Here are other characterizations\,$:$
\begin{itemize}

\item
The spherical functions are
the continuous bi\ssf--\ssf$K$\!--\ssf invariant functions \,$\phi$ on \,$G$
which satisfy the functional equation
\begin{equation}\label{FunctionalEquation1Sn}
\int_K\ssb dk\;\phi(x\ssf k\ssf y)=\phi(x)\,\phi(y)
\qquad\forall\;x,y\!\in\!G\ssf.
\end{equation}

\item
The spherical functions are
the continuous bi\ssf--\ssf$K$\!--\ssf invariant functions \,$\phi$ on \,$G$
such that
\begin{equation}\label{FunctionalEquation2Sn}
f\,\longmapsto\int_G\ssb dx\,f(x)\,\phi(x)
\end{equation}
defines a character of the (commutative) convolution algebra
\,$C_c(K\backslash G/K)$\ssf.

\item
The spherical functions are the matrix coefficients
\begin{equation*}
\phi(x)=\langle\ssf\pi(x)\ssf v,v\ssf\rangle\,,
\end{equation*}
where \,$\pi$ is a continuous unitary representation of \,$G$,
which has nonzero \,$K$\!--\ssf fixed vectors and which is irreducible,
and \,$v$ is a \,$K$\!--\ssf fixed vector,
which is normalized by \,$|v|\!=\!1\ssf$.

\item
{\rm Integral representation\,:}
\begin{equation}\begin{aligned}\label{SphericalFunctionsSn}
\phi_{\ssf\ell\ssf}(\cos\theta_1)
&=\tfrac{\Gamma(\frac n2)\vphantom{\frac0|}}
{\sqrt{\pi\ssf}\,\Gamma(\frac{n-1}2)\vphantom{\frac|0}}
\int_{\,0}^{\ssf\pi}\hspace{-1mm}d\ssf\theta_2\,(\ssf\sin\theta_2)^{n-2}\ssf
\bigl[\ssf\cos\theta_1\!
+\ssb i\,(\ssf\sin\theta_1)\ssf(\cos\theta_2)\ssf\bigr]^{\ssf\ell}\\
&=\tfrac{\Gamma(\frac n2)\vphantom{\frac0|}}
{\sqrt{\pi\ssf}\,\Gamma(\frac{n-1}2)\vphantom{\frac|0}}\,
(\ssf\sin\theta_1)^{2-n}\ssb
\int_{-\sin\theta_1}^{\,\sin\theta_1}\hspace{-1mm}ds\,
(\ssf\sin^2\ssb\theta_1\!-\ssb s^2\ssf)^{\ssb\frac{n-3}2}\ssf
(\cos\theta_1\!+\ssb i\,s\ssf)^\ell\,.
\end{aligned}\end{equation}

\end{itemize}
\end{remark}
The spherical Fourier expansion of radial functions on \ssf$\Sn$ reads
\begin{equation*}
f(x)=\sum\nolimits_{\ssf\ell\in\N}
d_{\ssf\ell}\,\langle\ssf f,\phi_{\ssf\ell\ssf}\rangle\,\phi_{\ssf\ell\ssf}(x)\,,
\end{equation*}
\vspace{-5mm}

where
\vspace{-1mm}
\begin{equation*}
d_{\ssf\ell}=\tfrac{n\,(n\ssf+\ssf2\ssf\ell\ssf-1)\,(n\ssf+\ssf\ell\ssf-\ssf2)\ssf!}{n\,!\;\ell\,!}
\end{equation*}
\vspace{-5.5mm}

and
\begin{equation*}
\langle\ssf f,\phi_{\ssf\ell\ssf}\rangle
=\tfrac{\Gamma(\frac{n+1}2)\vphantom{\frac0|}}
{2\,\pi^{\frac{n+1}2}\vphantom{\frac|0}}
\int_{\ssf\Sn}\!dx\,f(x)\,\phi_{\ssf\ell\ssf}(x)
=\tfrac{\Gamma(\frac{n+1}2)\vphantom{\frac0|}}
{\sqrt{\pi\ssf}\,\Gamma(\frac n2)\vphantom{\frac|0}}
\int_{\,0}^{\ssf\pi}\!d\ssf\theta_1\,(\ssf\sin\theta_1)^{n-1}\ssf f(\cos\theta_1)\,\phi_{\ssf\ell\ssf}(\cos\theta_1)\,.
\end{equation*}

\subsection{Spherical Fourier analysis on real hyperbolic spaces}
\label{SphericalFourierAnalysisHn}

Real hyperbolic spaces \ssf$\Hn$ are the simplest examples
of Riemannian symmetric spaces of noncompact type.
They are simply connected Riemannian manifolds,
with constant negative sectional curvature.
Let us recall the following three models of \ssf$\Hn$.
\smallskip

$\bullet$
\textbf{\,Model 1\,: Hyperboloid}
\vspace{.5mm}

\begin{figure}[ht]
\psfrag{0}[l]{$0$}
\psfrag{e0}[l]{$e_{\ssf0}$}
\psfrag{en}[l]{$e_{\ssf n}$}
\psfrag{ej}[r]{$e_{\ssf j}$}
\psfrag{Korbit}[c]{\color{red}$K$\!--\ssf orbit}
\psfrag{Norbit}[c]{\color{blue}$N$\ssb--\ssf orbit}
\psfrag{Aorbit}[c]{\color{green}$A$--\ssf orbit}
\includegraphics[height=60mm]{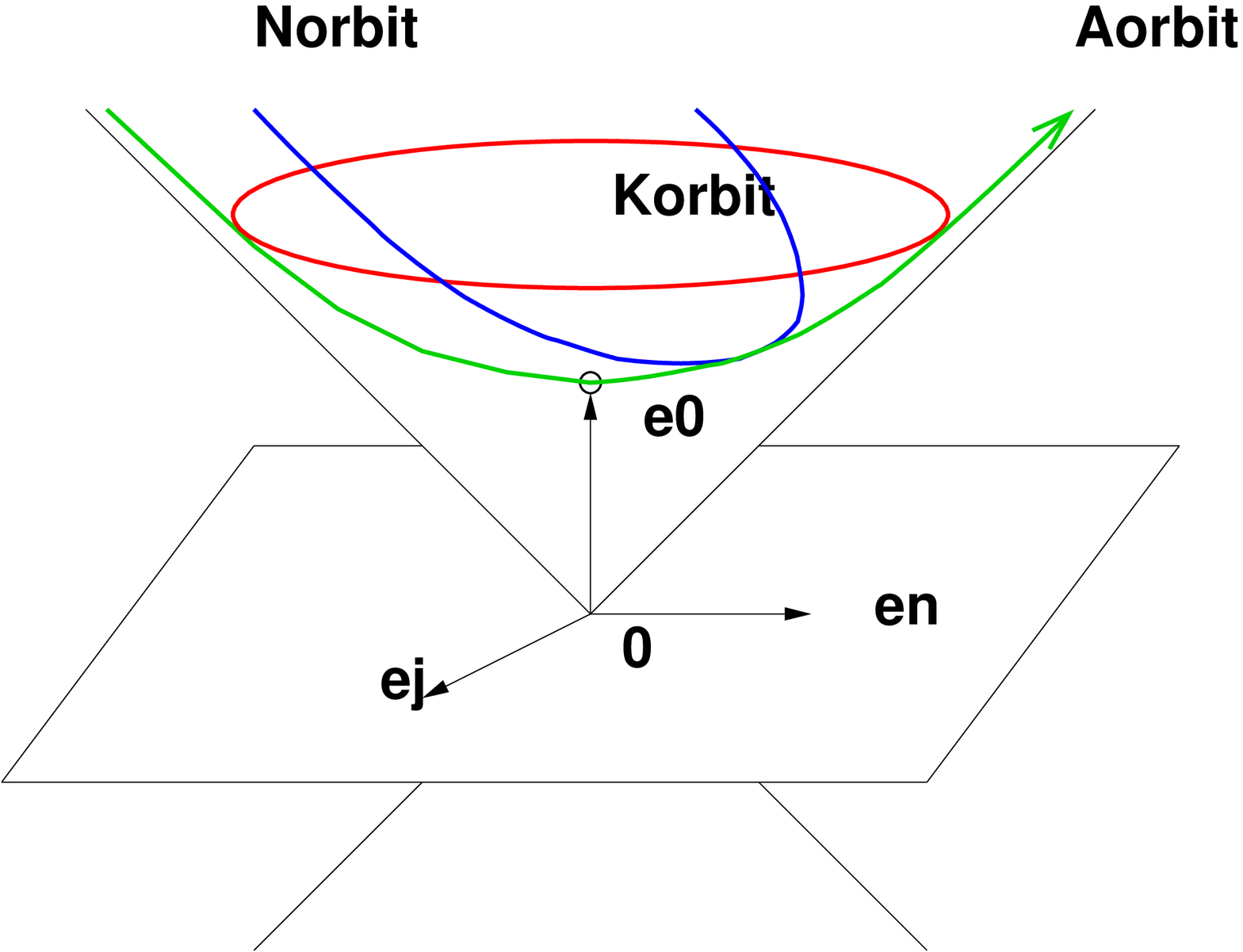}
\caption{Hyperboloid model of \ssf$\Hn$}
\end{figure}

In this model, \ssf$\Hn$ consists of the hyperboloid sheed
\begin{equation*}
\{\,x\ssb=\ssb(x_0,x_1,\dots,x_n)\ssb\in\ssb\Rnplusone\!=\ssb\R\!\times\!\Rn\,|\,L(x,x)\ssb=\ssb-\,1\ssf,\,x_0\ssb\ge\ssb1\,\}
\end{equation*}
defined by the Lorentz quadratic form
\,$L(x,x)\ssb=\ssb-\,x_0^2\ssb+\ssb x_1^2\ssb+\ssf\dots\ssf+\ssb x_n^2$\ssf.
The Riemannian structure is given by the metric
\ssf$ds^{\ssf2}\!=\ssb L(dx,dx)$\ssf,
restricted to the tangent bundle of \ssf$\Hn$,
and the Laplacian by \ssf$\Delta f\ssb=\ssb
L\bigl(\tfrac\partial{\partial x},\tfrac\partial{\partial x}\bigr)
\ssf\widetilde{f}\,\big|_{\ssf\Hn}$\ssf,
where \ssf$\widetilde{f}(x)\ssb=\ssb f\bigl(\tfrac x{\sqrt{-L(x,x)}\ssf}\bigr)$
\ssf denotes the homogeneous extension of \ssf$f$ \ssf to the light cone
\ssf$\{\,x\!\in\!\Rnplusone\,|\,L(x,x)\!<\ssb0\ssf,\,x_0\ssb>\ssb0\,\}$\ssf.
\smallskip

$\bullet$
\textbf{\,Model 2\,: Upper half--space}
\vspace{.5mm}

\begin{figure}[b]
\psfrag{en}[l]{$e_{\ssf n}$}
\psfrag{yn}[l]{$y_{\ssf n}$}
\psfrag{Rn-1}[l]{$\Rnminusone$}
\psfrag{Korbit}[l]{\color{red}$K$\!--\ssf orbit}
\psfrag{Norbit}[l]{\color{blue}$N$\ssb--\ssf orbit}
\psfrag{Aorbit}[l]{\color{green}$A$--\ssf orbit}
\includegraphics[height=60mm]{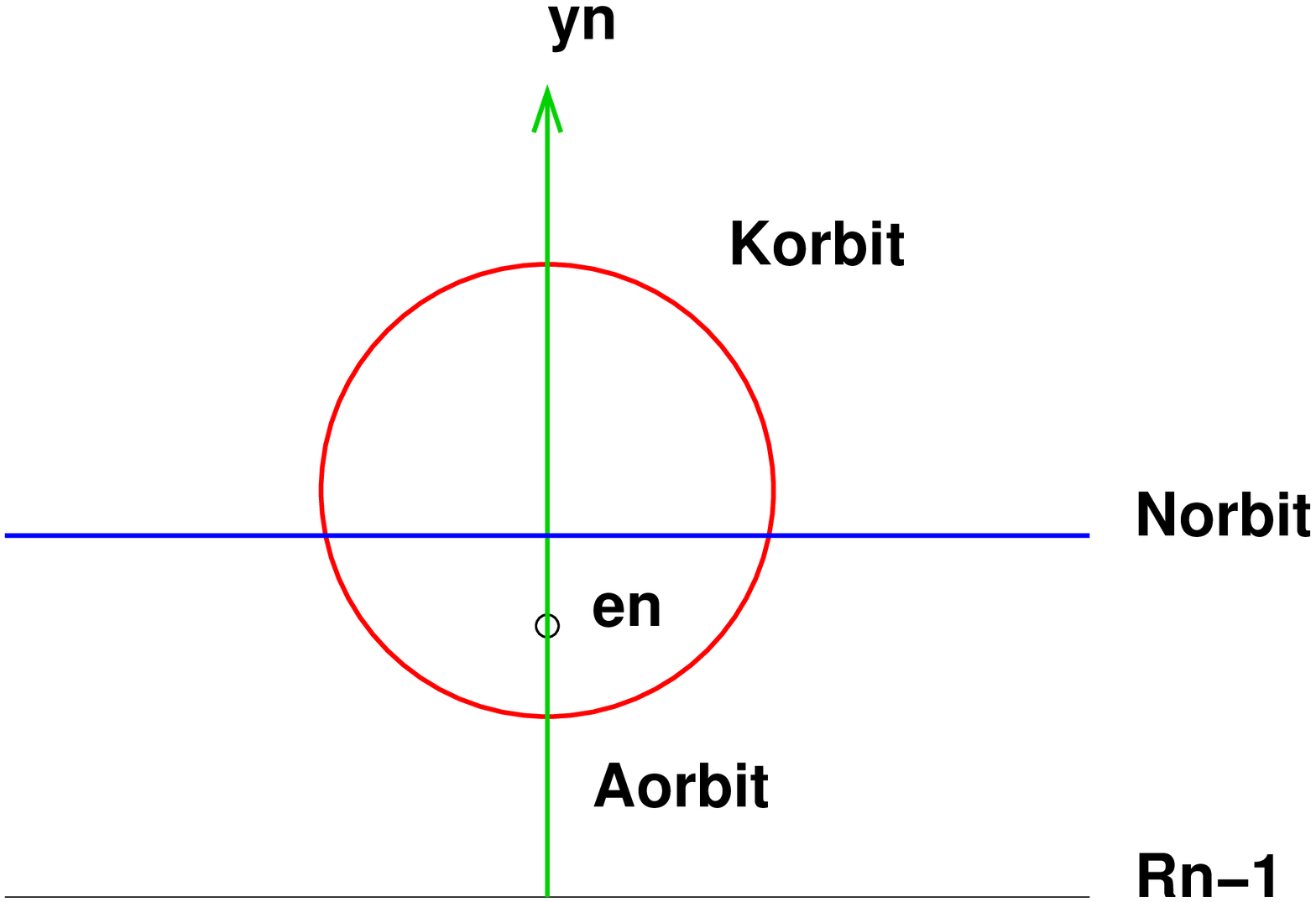}
\caption{Upper half--space model of \ssf$\Hn$}
\end{figure}

In this model, \ssf$\Hn$ consists of the upper half--space
\ssf$\R_+^{\ssf n}\!=\ssb\{\,y\!\in\!\R^n\,|\,y_n\ssb>\ssb0\ssf\}$
\ssf equipped with the Riemannian metric
\,$ds^{\ssf2}\ssb=y_n^{-2}\,|dy|^2$\ssf.
The volume is given by \ssf$d\ssb\vol\ssb=\ssb y_n^{-n}\,dy_1\dots\ssf dy_n$
and the Laplacian by
\begin{equation*}
\Delta
=y_n^2\,\sum\nolimits_{\ssf j=1}^{\,n}\tfrac{\partial^{\ssf2}}{\partial y_j^2}
-(n\!-\!2)\,y_n\,\tfrac{\partial}{\partial y_n}\,.
\end{equation*}

\pagebreak

$\bullet$
\textbf{\,Model 3\,: Ball}
\vspace{.5mm}

\begin{figure}[t]
\psfrag{0}[l]{$0$}
\psfrag{en}[l]{$e_{\ssf n}$}
\psfrag{Korbit}[c]{\color{red}$K$\!--\ssf orbit}
\psfrag{Norbit}[l]{\color{blue}$N$\ssb--\ssf orbit}
\psfrag{Aorbit}[c]{\color{green}$A$--\ssf orbit}
\includegraphics[height=60mm]{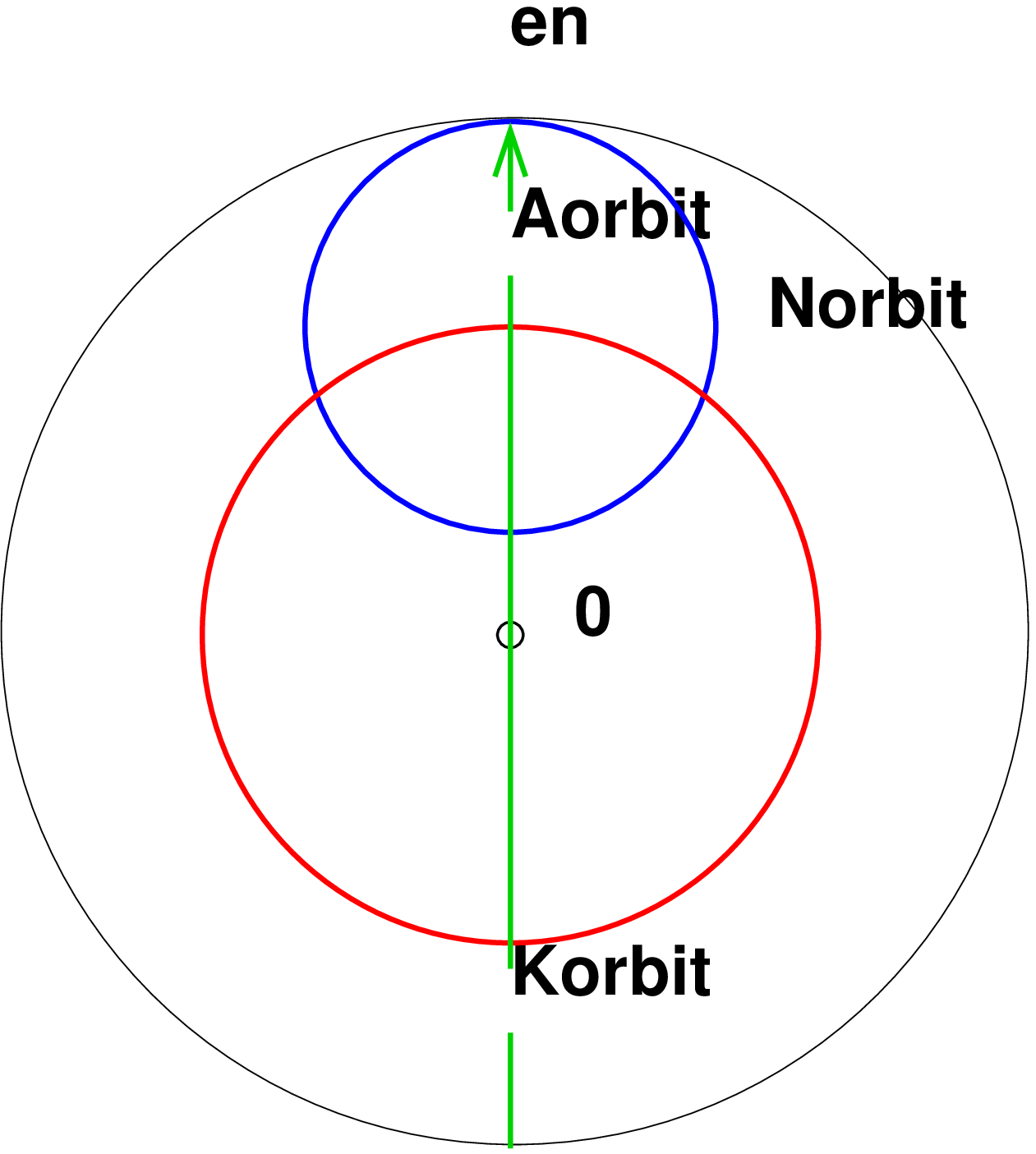}
\caption{Ball model of \ssf$\Hn$}
\end{figure}

In this model, \ssf$\Hn$ consists of the unit ball
\ssf$\mathbb{B}^n\!=\ssb\{\,z\!\in\!\R^n\,|\,|z|\!<\!1\ssf\}$\ssf.
The Riemannian metric is given by
\ssf$ds^2\ssb=\ssb\bigl(\tfrac{1\ssf-\ssf|z|^2}2\bigr)^{\!-2}\,|dz|^2$\ssf,
the volume by \ssf$d\ssb\vol\ssb
=\ssb\bigl(\tfrac{1\ssf-\ssf|z|^2}2\bigr)^{\!-n}\,dz_1\dots\ssf dz_n$\,,
the distance to the origin by
\,$r\ssb=\ssb2\ssf\artanh|z|\ssb
=\ssb\log\tfrac{1\ssf+\ssf|z|}{1\ssf-\ssf|z|}$\,,
and the Laplacian by
\begin{equation*}
\Delta=\bigl(\tfrac{1\ssf-\ssf|z|^2}2\bigr)^{\ssb2}\,
\sum\nolimits_{\ssf j=1}^{\,n}\tfrac{\partial^{\ssf2}}{\partial z_j^2}
+(n\!-\!2)\,\tfrac{1\ssf-\ssf|z|^2}2\,
\sum\nolimits_{\ssf j=1}^{\,n}z_j\,\tfrac{\partial}{\partial z_j}\,.
\end{equation*}

\begin{remark}
Model 1 and Model 3 are mapped onto each other by
the stereographic projection with respect to \ssf$-\ssf e_0$\ssf,
while Model 2 and Model 3 are mapped onto each other
by the inversion with respect to the sphere
\,$S(\ssb-\ssf e_n,\ssb\sqrt{2\ssf})$\ssf.
This leads to the following formulae
\begin{align*}
&\begin{cases}
\,x_0=\frac{1\,+\,|y|^2}{2\ssf y_n}=\frac{1\,+\,|z|^2}{1\,-\,|z|^2}\\
\,x_j=\frac{y_j}{y_n}=\frac{2\ssf z_j}{1\,-\,|z|^2}
&(\ssf j\!=\!1,\dots,\ssb n\!-\!1\ssf)\\
\,x_n=\frac{1\,-\,|y|^2}{2\ssf y_n}=\frac{2\ssf z_n}{1\,-\,|z|^2}
\end{cases}\\
&\begin{cases}
\,y_j=\frac{x_j}{x_0\ssf+\,x_n}=\frac{2\ssf z_j}{1\,+\,|z|^2\ssf+\,2\ssf z_n}
&(\ssf j\!=\!1,\dots,\ssb n\!-\!1\ssf)\\
\,y_n=\frac1{x_0\ssf+\,x_n}=\frac{1\,-\,|z|^2}{1\,+\,|z|^2\ssf+\,2\ssf z_n}
\end{cases}\\
&\begin{cases}
\,z_j
=\frac{x_j}{1\,+\,x_0}
=\frac{2\ssf y_j}{1\,+\,|y|^2\ssf+\,2\ssf y_n}
&(\ssf j\!=\!1,\dots,\ssb n\!-\!1\ssf)\\
\,z_n
=\frac{x_n}{1\,+\,x_0}
=\frac{1\,-\,|y|^2}{1\,+\,|y|^2\ssf+\,2\ssf y_n}\,.
\end{cases}
\end{align*}
\end{remark}

Let \ssf$G$ \ssf be the isometry group of \ssf$\Hn$
and let \ssf$K$ be the stabilizer of a base point in \ssf$\Hn$.
Then \ssf$\Hn$ can be realized as the homogeneous space \ssf$G/K$.
In Model 1, \ssf$G$ \ssf is made up of two among the four connected components
of the Lorentz group \ssf$\text{O}(1,n)$\ssf,
and the stabilizer of \ssf$e_0$ is \ssf$K\!=\ssb\text{O}(n)$\ssf.
Consider the subgroup \ssf$A\ssb\approx\ssb\R$ \ssf in \ssf$G$
\ssf consisting of
\begin{itemize}

\item
the matrices
\begin{equation*}
a_{\ssf r}=\left(\;\begin{matrix}
\cosh r & 0 & \sinh r \\ 0 & I & 0 \\ \sinh r & 0 & \cosh r
\end{matrix}\;\right)
\qquad(\ssf r\!\in\!\R\ssf)
\end{equation*}
in Model 1,

\item
the dilations \ssf$a_{\ssf r}\!:\ssb y\longmapsto\ssb e^{-r}y$
\ssf in Model 2,

\end{itemize}

and the subgroup \,$N\!\approx\ssb\R^{n-1}$ \ssf
consisting of horizontal translations
\,$n_\upsilon\!:\ssb y\ssb\longmapsto\ssb y\!+\!\upsilon$
\linebreak
\ssf(\ssf$\upsilon\!\in\!\R^{n-1}$)
\ssf in Model 2.
Then we have
\begin{itemize}

\item
the Cartan decomposition \ssf$G\ssb=\ssb K\overline{A^+}\ssb K$,
which corresponds to polar coordinates in Model 3,

\item
the Iwasawa decomposition \ssf$G\ssb=\ssb N\ssb A\ssf K$,
which corresponds to Cartesian coordinates in Model 2.

\end{itemize}
We shall denote by \ssf$a_{\ssf r(g)}$ and \ssf$a_{\vsf h(g)}$
the $\overline{A^+}$ and $A$ components of \ssf$g\!\in\!G$
in the Cartan and Iwasawa decompositions.

\begin{remark}
In small dimensions,
\,$\mathbb{H}^2\ssb\approx\text{\rm SL}(2,\R)\ssf/\,\text{\rm SO}(2)$
\,and
\,$\mathbb{H}^3\ssb\approx\text{\rm SL}(2,\C)\ssf/\,\text{\rm SU}(2)$\ssf.
\end{remark}

As usual, we identify right\ssf--\ssf$K$\ssb--\ssf invariant functions on \ssf$G$
\ssf with functions on \ssf$\Hn$,
and bi\ssf--\ssf$K$\ssb--\ssf invariant functions on \ssf$G$
\ssf with radial functions on \ssf$\Hn$ 
i.e.~functions on \ssf$\Hn$
which depend only on the distance \ssf$r$ \ssf to the origin.
For radial functions \ssf$f\ssb=\ssb f(r)$\ssf,
\vspace{-1mm}
\begin{equation*}
\int_{\ssf\Hn}\!d\vol\,f\ssf
=\,2\,\tfrac{\pi^{\frac n2}\vphantom{\frac00}}
{\Gamma(\frac n2)\vphantom{\frac|0}}
\int_{\,0}^{+\infty}\hspace{-1mm}dr\,(\sinh r)^{n-1}\ssf f(r)
\end{equation*}
\vspace{-4.5mm}

and
\vspace{-.5mm}
\begin{equation*}\label{RadialLaplacianHn}
\Delta\ssf f=\tfrac{\partial^{\ssf2}\ssb f}{\partial\ssf r^{\ssf2}}
+(n\ssb-\!1)\ssf(\coth r)\,\tfrac{\partial f}{\partial\ssf r}\,.
\end{equation*}
The spherical functions \ssf$\phi_\lambda$ are
the smooth normalized radial eigenfunctions of the Laplacian on \ssf$\Hn$.
Specifically,
\begin{equation*}\begin{cases}
\;\Delta\,\phi_\lambda
=-\,\bigl\{\ssf\lambda^2\!+\ssb\rho^{\ssf2\ssf}\bigr\}\,\phi_\lambda\,,\\
\;\phi_\lambda(0)=1\,,
\end{cases}\end{equation*}
where \,$\rho\ssb=\ssb\frac{n-1}2$\ssf.
\smallskip

\begin{remark}
The spherical functions on \,$\Hn$ can be characterized again in several other ways.
Notably,
\begin{itemize}

\item
{\rm Differential equation\,:}
the function \,$\phi_\lambda(r)$ \ssf is
the unique smooth solution to the differential equation
\begin{equation*}
\bigl(\tfrac\partial{\partial\ssf r}\bigr)^{\ssb2}\ssf\phi_\lambda\ssb
+(n\ssb-\!1)\ssf(\coth r)\,\bigl(\tfrac\partial{\partial\ssf r}\bigr)\,\phi_\lambda\ssb
+(\ssf\lambda^2\!+\ssb\rho^{\ssf2\ssf})\,\phi_\lambda\ssb=0\,,
\end{equation*}
which is normalized by \ssf$\phi_\lambda(0)\!=\!1$\ssf.

\item
{\rm Relation with classical special functions\,:}
\begin{equation*}
\phi_\lambda(r)
=\phi_\lambda^{\frac{n-2}2,-\frac12}(r)
={}_2\ssf\text{\rm F}_{\ssb1}\bigl(\tfrac{\rho\ssf+\ssf i\ssf\lambda}2,
\ssb\tfrac{\rho\ssf-\ssf i\ssf\lambda}2\ssf;
\ssb\tfrac n2\ssf;-\sinh^2\ssb r\bigr)\,,
\end{equation*}
where \,$\phi_\lambda^{\alpha,\beta}$ are the Jacobi functions
and \,${}_2\ssf\text{\rm F}_{\ssb1}$ the Gauss hypergeometric function.

\item
Same functional equations
as \eqref{FunctionalEquation1Sn}
and \eqref{FunctionalEquation2Sn}.

\item
The spherical functions are the matrix coefficients
\begin{equation*}
\phi_\lambda(x)
=\tfrac{\Gamma(\frac n2)\vphantom{\frac0|}}
{2\,\pi^{\frac n2\vphantom{\frac00}}}\,
\langle\ssf\pi_\lambda(x)\ssf1\ssf,1\ssf\rangle\,,
\end{equation*}
of the spherical principal series representations of \,$G$
\ssb on \ssf$L^2(\ssf\Snminusone)$\ssf.

\item
According to the Harish--Chandra formula
\begin{equation*}\label{HarishChandraFormulalHn}
\phi_\lambda(x)=\int_K\!dk\;e^{\,(\rho\ssf-\ssf i\ssf\lambda)\,h(k\ssf x)}\,,
\end{equation*}
the function \ssf$\phi_\lambda$ is a spherical average of horocyclic waves.
Let us make this integral representation more explicit\,$:$
\begin{equation}\begin{aligned}\label{SphericalFunctionsHn}
\phi_\lambda(r)
&=\ssf\tfrac{\Gamma(\frac n2)\vphantom{\frac0|}}
{2\,\pi^{\frac n2\vphantom{\frac00}}}
\int_{\ssf\Snminusone}\hspace{-1mm}dv\;
\langle\,\cosh r\ssb-\ssb(\sinh r)\,v\ssf,e_n\,
\rangle^{\ssf i\ssf\lambda\ssf-\ssf\rho}\\
&=\ssf\tfrac{\Gamma(\frac n2)\vphantom{\frac0|}}
{\sqrt{\pi\ssf}\;\Gamma(\frac{n-1}2)\vphantom{\frac|0}}
\int_{\,0}^{\ssf\pi}\hspace{-1mm}d\theta\,
(\sin\theta)^{n-2}\,
\bigl[\ssf\cosh r\ssb-\ssb(\sinh r)\ssf(\cos\theta)\ssf
\bigr]^{\ssf i\ssf\lambda\ssf-\ssf\rho}\\
&=\ssf\tfrac{2^{\frac{n-1}2}\,\Gamma(\frac n2)\vphantom{\frac0|}}
{\sqrt{\pi\ssf}\;\Gamma(\frac{n-1}2)\vphantom{\frac|0}}\;(\sinh r)^{2-n}
\int_{\,0}^{\ssf r}\!ds\,
(\ssf\cosh r\ssb-\ssb\cosh s\ssf)^{\ssb\frac{n-3}2}\ssf
\cos\ssb\lambda\ssf s\,.
\end{aligned}\end{equation}

\end{itemize}
\end{remark}

\begin{remark}
The asymptotic behavior of the spherical functions is given
by the Harish--Chandra expansion
\begin{equation*}\label{HarishChandraExpansionHn}
\phi_\lambda(r)
=\mathbf{c}\vsf(\lambda)\,\Phi_\lambda(r)
+\mathbf{c}\vsf(\ssb-\lambda)\,\Phi_{-\lambda}(r)\,,
\end{equation*}
\vspace{-5mm}

where
\vspace{-.5mm}
\begin{equation*}\label{cfunctionHn}
\mathbf{c}\vsf(\lambda)=
\tfrac{\Gamma(2\ssf\rho)}{\Gamma(\rho)}\,
\tfrac{\Gamma(i\ssf\lambda)}{\Gamma(i\ssf\lambda\ssf+\ssf\rho)}
\end{equation*}
\vspace{-5mm}

and
\vspace{-1mm}
\begin{align*}
\Phi_\lambda(r)
&=(2\cosh r)^{i\lambda-\rho}\,
{}_2\ssf\text{\rm F}_{\ssb1}\bigl(\tfrac{\rho\ssf-\ssf i\ssf\lambda}2\ssf,
\ssb\tfrac{\rho\ssf+\ssf1\ssf-\ssf i\ssf\lambda}2\,;\ssb
1\!-\ssb i\ssf\lambda\,;\cosh^{-2}\ssb r\ssf\bigr)\\
&=\ssf e^{\,(i\ssf\lambda\ssf-\ssf\rho)\ssf r}\,
\sum\nolimits_{\ssf\ell=0}^{+\infty}\ssf
\Gamma_{\ssb\ell}(\lambda)\,e^{-2\ssf\ell\ssf r}\,,
\end{align*}
with \,$\Gamma_{\ssb0}\ssb\equiv\ssb1$\ssf.

\end{remark}

The spherical Fourier transform (or Harish--Chandra transform)
of radial functions on \ssf$\Hn$ is defined by
\begin{equation}\label{SphericalFourierTransformHn}
\mathcal{H}f(\lambda)\ssf
=\ssb\int_{\ssf\Hn}\!dx\,f(x)\,\phi_\lambda(x)
=\ssf\tfrac{\Gamma(\frac n2)\vphantom{\frac0|}}
{2\,\pi^{\frac n2\vphantom{\frac00}}}
\int_{\,0}^{+\infty}\hspace{-1mm}dr\,
(\sinh r)^{n-1}\ssf f(r)\,\phi_\lambda(r)
\end{equation}
and the inversion formula reads
\vspace{-.5mm}
\begin{equation}\label{InverseSphericalFourierTransformHn}
f(x)=\ssf2^{\ssf n-3}\,\pi^{-\frac n2-1}\,\Gamma(\tfrac n2)
\int_{\,0}^{+\infty}\hspace{-1mm}d\lambda\hspace{1mm}
|\vsf\mathbf{c}\vsf(\lambda)|^{-2}\hspace{1mm}
\mathcal{H}f(\lambda)\,\phi_\lambda(x)\,.
\end{equation}

\begin{remark}
${}$
\begin{itemize}

\item
The Plancherel density reads
\vspace{-1.5mm}
\begin{equation*}\label{PlancherelMeasureHnodd}
|\ssf\mathbf{c}\vsf(\lambda)\vsf|^{-2}=
\tfrac{\pi\vphantom{\frac00}}
{2^{\vsf2\vsf n-4}\,\Gamma(\frac n2)^2\vphantom{\frac00}}\,
\prod\nolimits_{\ssf j=0}^{\frac{n-3}2}\ssf
(\ssf\lambda^2\!+\ssb j^{\ssf2}\ssf)
\end{equation*}
\vspace{-5mm}

\noindent
in odd dimension, and
\vspace{-.5mm}
\begin{equation*}\label{PlancherelMeasureHneven}
|\ssf\mathbf{c}\vsf(\lambda)\vsf|^{-2}=
\tfrac{\pi\vphantom{\frac00}}
{2^{\vsf2\vsf n-4}\,\Gamma(\frac n2)^2\vphantom{\frac00}}\,
\lambda\ssf\tanh\pi\lambda\,
\prod\nolimits_{\,j=0}^{\frac n2-1}\ssf
\bigl[\ssf\lambda^2\!
+\hspace{-.4mm}(\vsf j\!+\!\tfrac12)^2\hspace{.4mm}\bigr]
\end{equation*}
in even dimension.
Notice the different behaviors
\vspace{-.5mm}
\begin{equation*}\label{BehaviorInfinityPlancherelMeasureHn}
|\ssf\mathbf{c}\vsf(\lambda)\vsf|^{-2}
\sim\tfrac{\pi\vphantom{\frac00}}
{2^{\vsf2\vsf n-4}\,\Gamma(\frac n2)^2\vphantom{\frac|0}}\,
|\lambda|^{\ssf n-1}
\end{equation*}
\vspace{-5mm}

\noindent
at infinity, and
\vspace{-1mm}
\begin{equation*}\label{BehaviorOriginPlancherelMeasureHn}
|\ssf\mathbf{c}\vsf(\lambda)\vsf|^{-2}
\sim\tfrac{\pi\,\Gamma(\frac{n-1}2)^2\vphantom{\frac0|}}
{2^{\vsf2\vsf n-4}\,\Gamma(\frac n2)^2\vphantom{\frac|0}}\,
\lambda^2
\end{equation*}
\vspace{-5mm}

\noindent
at the origin.

\item
Observe that \eqref{SphericalFourierTransformHn}
and \eqref{InverseSphericalFourierTransformHn}
are not symmetric,
unlike \eqref{FourierTransformRn} and \eqref{InverseFourierTransformRn},
or \eqref{HankelTransform} and \eqref{InverseHankelTransform}.

\end{itemize}
\end{remark}

The spherical Fourier transform \eqref{SphericalFourierTransformHn},
which is somewhat abstract,
can be bypassed by considering the Abel transform,
which is essentially the horocyclic Radon transform restricted to radial functions.
Specifically,
\vspace{-.5mm}
\begin{align*}\label{AbelTransformHn}
\mathcal{A}f(r)\ssf
&=\,e^{-\rho\ssf r}\int_Ndn\,f(n\ssf a_{\ssf r})\\
&=\ssf\tfrac{(2\ssf\pi)^{\ssb\frac{n-1}2}\vphantom{\frac00}}
{\Gamma(\frac{n-1}2)\vphantom{\frac|0}}\ssf
\int_{\,|r|}^{+\infty}\hspace{-1mm}ds\,\sinh\ssb s\,
(\ssf\cosh s\ssb-\ssb\cosh r\ssf)^{\frac{n-3}2}\ssf f(s)\,.
\nonumber\end{align*}
Then the following commutative diagram holds,
let say in the Schwartz space setting\,:
\vspace{1mm}
\begin{equation*}\label{CommutativeDiagramHn}\begin{aligned}
&\hspace{-2.5mm}\mathcal{S}_{\ssf\text{even}}(\R)\\
\mathcal{H}\nearrow\hspace{2mm}
&\hspace{10mm}\nwarrow\mathcal{F}\\
\mathcal{S}_{\ssf\text{rad}}(\Hn)\hspace{5mm}
&\underset{\textstyle\mathcal{A}}\longrightarrow\hspace{5mm}
\mathcal{S}_{\ssf\text{even}}(\R)
\end{aligned}\end{equation*}
Here \ssf$\mathcal{S}_{\ssf\text{rad}}(\Hn)$
denotes the $L^2$ radial Schwartz space on \ssf$\Hn$,
which can be identified with
\,$(\cosh r)^{-\rho}\,\mathcal{S}_{\ssf\text{even}}(\R)$\ssf,
\,$\mathcal{F}$ \ssf the Euclidean Fourier transform on \ssf$\R$
\ssf and each arrow is an isomorphism.
Thus the inversion of the spherical Fourier transform \ssf$\mathcal{H}$
\ssf boilds down to the inversion of the Abel transform \ssf$\mathcal{A}$\ssf.
In odd dimension,
\begin{equation*}
\mathcal{A}^{-1}\ssb g\ssf(r)
=(2\ssf\pi)^{-\frac{n-1}2}\,\bigl(-\ssf\tfrac1{\sinh r}\ssf
\tfrac\partial{\partial\ssf r}\ssf\bigr)^{\ssb\frac{n-1}2}\ssf g(r)
\end{equation*}
while, in even dimension,
\begin{equation*}
\mathcal{A}^{-1}\ssb g\ssf(r)
=\tfrac{1\vphantom{\frac00}}{2^{\frac{n-1}2}\pi^{\frac n2}\vphantom{\frac|0}}
\int_{\,|r|}^{+\infty}\!\tfrac{ds}{\sqrt{\ssf\cosh s\,-\,\cosh r\ssf}}\;
\bigl(-\tfrac{\partial}{\partial s}\bigr)\ssf
\bigl(-\tfrac1{\sinh s}\ssf\tfrac\partial{\partial s}\bigr)^{\ssb\frac n2-1}
\ssf g\ssf(s)\,.
\end{equation*}
Consider finally the transform 
\begin{equation*}
\mathcal{A}^*g(x)=\int_K\!dk\;e^{\,\rho\,h(k\ssf x)}\,g(h(k\ssf x))\,,
\end{equation*}
which is dual to the Abel transform, i.e.,
\begin{equation*}
\int_{\ssf\Hn}\!dx\;f(x)\,\mathcal{A}^*g(x)\,
=\int_{-\infty}^{+\infty}\hspace{-1mm}d\ssf r\;\mathcal{A}f(r)\,g(r)\,,
\end{equation*}
and which is an isomorphism
between \ssf$\mathcal{C}_{\ssf\text{even}}^\infty(\R)$
and \ssf$\mathcal{C}_{\ssf\text{rad}}^{\ssf\infty}(\Hn)$\ssf.
It is given explicitly by
\begin{equation*}\label{DualAbelTransformHn}
\mathcal{A}^*\ssb g\ssf(r)
=\ssf\tfrac{2^{\frac{n-1}2}\,\Gamma(\frac n2)\vphantom{\frac0|}}
{\sqrt{\pi\ssf}\;\Gamma(\frac{n-1}2)\vphantom{\frac|0}}\;
(\sinh r)^{2-n}\int_{\,0}^{\ssf r}\!ds\,
(\ssf\cosh r\ssb-\ssb\cosh s\ssf)^{\ssb\frac{n-3}2}\ssf g\ssf(s)
\end{equation*}
and its inverse by
\begin{equation*}
(\mathcal{A}^*)^{-1}\ssb f\ssf(r)
=\ssf\tfrac{\sqrt{\pi\,}\vphantom{\frac0|}}
{2^{\frac{n-1}2}\,\Gamma(\frac n2)\vphantom{\frac0|}}\,
\tfrac\partial{\partial\ssf r}\,
\bigl(\ssf\tfrac1{\sinh r}\ssf\tfrac\partial{\sinh r}\ssf\bigr)^{\!\frac{n-3}2}
\ssf\bigl\{(\ssf\sinh r)^{\ssf n-2}\ssf f(r)\bigr\}
\end{equation*}
in odd dimension and by
\begin{equation*}
(\mathcal{A}^*)^{-1}\ssb f(r)
=\tfrac{1\vphantom{\frac00}}
{2^{\frac{n-1}2}\ssf(\frac n2-1)\ssf!\vphantom{\frac|0}}\,
\tfrac\partial{\partial\ssf r}\,
\bigl(\ssf\tfrac1{\sinh r}\ssf\tfrac\partial{\partial\ssf r}\ssf\bigr)^{\!\frac n2-1}\!
\int_{\,0}^{\ssf r}\ssb
\tfrac{ds}{\sqrt{\ssf\cosh r\,-\,\cosh s\ssf}}\,
(\ssf\sinh s)^{n-1}\ssf f(s)
\end{equation*}
in even dimension.

\begin{remark}
Notice that the spherical function \,$\phi_\lambda(r)$
is the dual Abel transform of the cosine function \,$\cos\lambda\ssf s$\ssf.
\end{remark}

\paragraph{\bf Applications.}
Spherical Fourier analysis is an efficient tool
for solving invariant PDEs on \hspace{.5mm}$\Hn$.
Here are some examples of evolution equations.
\smallskip

$\bullet$
\,The heat equation
\begin{equation*}\label{HeatEquationHn}\begin{cases}
\,\partial_{\ssf t\ssf}u(x,t)=\Delta_{\ssf x\ssf}u(x,t)\\
\,u(x,0)=f(x)
\end{cases}\end{equation*}
can be solved explicitly by means of the inverse Abel transform.
Specifically,
\begin{equation*}
u(x,t)=f\!*\ssb h_{\ssf t\ssf}(x)\,,
\end{equation*}
where the heat kernel is given by
\begin{equation*}
h_{\ssf t\ssf}(r)
=\tfrac{1\vphantom{\frac00}}{2^{\frac{n+1}2}\pi^{\frac n2}}\,
t^{\ssf-\ssf\frac12}\,e^{-\ssf\rho^{\ssf2\ssf}t}\,
\bigl(-\ssf\tfrac1{\sinh r}\ssf
\tfrac\partial{\partial\ssf r}\ssf\bigr)^{\ssb\frac{n-1}2}\ssf
e^{-\ssf\frac{r^2}{4\ssf t}}
\end{equation*}
in odd dimension and by
\begin{equation*}
h_{\ssf t\ssf}(r)
=(2\ssf\pi)^{-\frac{n+1}2}\,
t^{\ssf-\ssf\frac12}\,e^{-\ssf\rho^{\ssf2\ssf}t}
\int_{\,|r|}^{+\infty}\!\tfrac{ds}{\sqrt{\ssf\cosh s\,-\,\cosh r\ssf}}\;
\bigl(-\tfrac{\partial}{\partial s}\bigr)\ssf
\bigl(-\tfrac1{\sinh s}\ssf\tfrac\partial{\partial s}\bigr)^{\ssb\frac n2-1}
\ssf e^{-\ssf\frac{s^2}{4\ssf t}}
\end{equation*}
in even dimension.
Moreover, the following global estimate holds\,:
\begin{equation}\label{HeatKernelEstimateHn}
h_{\ssf t\ssf}(r)\,
\asymp\,e^{-\ssf\rho^{\ssf2\ssf}t}\,e^{-\ssf\rho\,r}\,
e^{-\frac{r^2}{4\ssf t}}\ssf\times\ssf\begin{cases}
\,t^{\ssf-\frac32}\,(1\!+\ssb r)
&\text{if \,}t\ssb\ge\!1\!+\ssb r\ssf,\\
\,t^{\ssf-\frac n2}\,(1\!+\ssb r)^{\frac{n-1}2}
&\text{if \,}0\!<\hspace{-.5mm}t\ssb\le\!1\!+\ssb r\ssf.\\
\end{cases}\end{equation}

$\bullet$
\,Similarly for the Schr\"odinger equation
\begin{equation*}\label{SchroedingerEquationHn}\begin{cases}
\,i\,\partial_{\ssf t\ssf}u(x,t)=\Delta_{\ssf x\ssf}u(x,t)\ssf,\\
\,u(x,0)=f(x)\ssf.
\end{cases}\end{equation*}
In this case
\hspace{.5mm}$u(x,t)\ssb=\hspace{-.5mm}f\!*\ssb h_{-i\ssf t\ssf}(x)$\ssf,
where
\begin{equation*}
h_{-i\ssf t\ssf}(r)
=\tfrac{1\vphantom{\frac00}}{2^{\frac{n+1}2}\pi^{\frac n2}}\,
e^{\,i\ssf\frac\pi4\sign(t)}\,|t|^{\ssf-\ssf\frac12}\,
e^{-\ssf i\,\rho^{\ssf2\ssf}t}\,
\bigl(-\ssf\tfrac1{\sinh r}\ssf
\tfrac\partial{\partial\ssf r}\ssf\bigr)^{\ssb\frac{n-1}2}\ssf
e^{-\ssf i\ssf\frac{r^2}{4\ssf t}}
\end{equation*}
in odd dimension,
\begin{align*}
h_{-i\ssf t\ssf}(r)\ssf
&=\ssf(2\ssf\pi)^{-\frac{n+1}2}\,
e^{\,i\ssf\frac\pi4\sign(t)}\,|t|^{\ssf-\ssf\frac12}\,
e^{-\ssf i\,\rho^{\ssf2\ssf}t}\\
&\,\times\ssf
\int_{\,|r|}^{+\infty}\!\tfrac{ds}{\sqrt{\ssf\cosh s\,-\,\cosh r\ssf}}\;
\bigl(-\tfrac{\partial}{\partial s}\bigr)\ssf
\bigl(-\tfrac1{\sinh s}\ssf\tfrac\partial{\partial s}\bigr)^{\ssb\frac n2-1}
\ssf e^{-\ssf i\ssf\frac{s^2}{4\ssf t}}
\end{align*}
in even dimension, and
\begin{equation*}\begin{aligned}
|\ssf h_{-i\ssf t\ssf}(r)|\,
&\lesssim\,e^{-\ssf\rho\,r}\ssf\times\ssf\begin{cases}
\,|t|^{\ssf-\frac32}\,(1\!+\ssb r)
&\text{if \,}|t|\!\ge\!1\!+\ssb r\\
\,|t|^{\ssf-\frac n2}\,(1\!+\ssb r)^{\frac{n-1}2}
&\text{if \,}0\!<\!|t|\!\le\!1\!+\ssb r\\
\end{cases}\end{aligned}\end{equation*}
in all dimensions.
\smallskip

$\bullet$
\,The shifted wave equation
\begin{equation*}\label{WaveEquationHn}\begin{cases}
\,\partial_{\ssf t}^{\ssf2}u(x,t)
=\bigl(\Delta_{\ssf x}\!+\ssb\rho^{\ssf2\ssf}\bigr)\ssf u(x,t)\\
\,u(x,0)=f(x)\ssf,\,
\partial_{\ssf t}|_{\ssf t=0}\,u(x,t)=g(x)
\end{cases}\end{equation*}
can be solved explicitly by means of the inverse dual Abel transform.
Specifically,
\begin{align*}
u(t,x)
&=\tfrac{1\vphantom{\frac00}}{2^{\frac{n+1}2}\pi^{\frac{n-1}2}}\,
\tfrac\partial{\partial\ssf t}\,\bigl(\ssf\tfrac1{\sinh t}\ssf
\tfrac\partial{\partial\ssf t}\ssf\bigr)^{\!\frac{n-3}2}\ssf
\Bigl\{\ssf\tfrac1{\sinh t}\int_{S(x,\ssf|t|)}\!dy\hspace{.75mm}f(y)\Bigr\}\\
&\ssf+\tfrac{1\vphantom{\frac00}}{2^{\frac{n+1}2}\pi^{\frac{n-1}2}}\,
\bigl(\ssf\tfrac1{\sinh t}\ssf
\tfrac\partial{\partial\ssf t}\ssf\bigr)^{\!\frac{n-3}2}\ssf
\Bigl\{\ssf\tfrac1{\sinh t}\int_{S(x,\ssf|t|)}\!dy\;g(y)\Bigr\}
\end{align*}
in odd dimension and
\begin{align*}
u(t,x)\ssf
&=\tfrac{1\vphantom{\frac00}}{2^{\frac{n+1}2}\pi^{\frac n2}}\,
\tfrac{\partial}{\partial|t|}\,
\bigl(\ssf\tfrac{\partial}{\partial\ssf(\cosh t)}\bigr)^{\!\frac n2-1}
\int_{B(x,\ssf|t|)}\!dy\;
\tfrac{f(y)}{\sqrt{\ssf\cosh t\,-\,\cosh d\ssf(y,\ssf x)\ssf}}\\
&\ssf+\tfrac{1\vphantom{\frac00}}{2^{\frac{n+1}2}\pi^{\frac n2}}\,
\operatorname{sign}(t)\,
\bigl(\ssf\tfrac{\partial}{\partial\ssf(\cosh t)}\bigr)^{\!\frac n2-1}\!
\int_{B(x,\ssf|t|)}\!dy\;
\tfrac{g(y)}{\sqrt{\ssf\cosh t\,-\,\cosh d\ssf(y,\ssf x)\ssf}}
\end{align*}
in even dimension.

\subsection{Spherical Fourier analysis on homogeneous trees}
\label{SphericalFourierAnalysisTree}

\begin{figure}[hb]
\psfrag{0}[c]{$0$}
\psfrag{1}[c]{$1$}
\psfrag{2}[c]{$2$}
\psfrag{-1}[c]{$-1$}
\psfrag{h}[c]{$h$}
\psfrag{omega}[c]{$\omega$}
\psfrag{geodesic}[l]{\color{green}geodesic}
\psfrag{horocycle}[l]{\color{blue}horocycle}
\includegraphics[height=70mm]{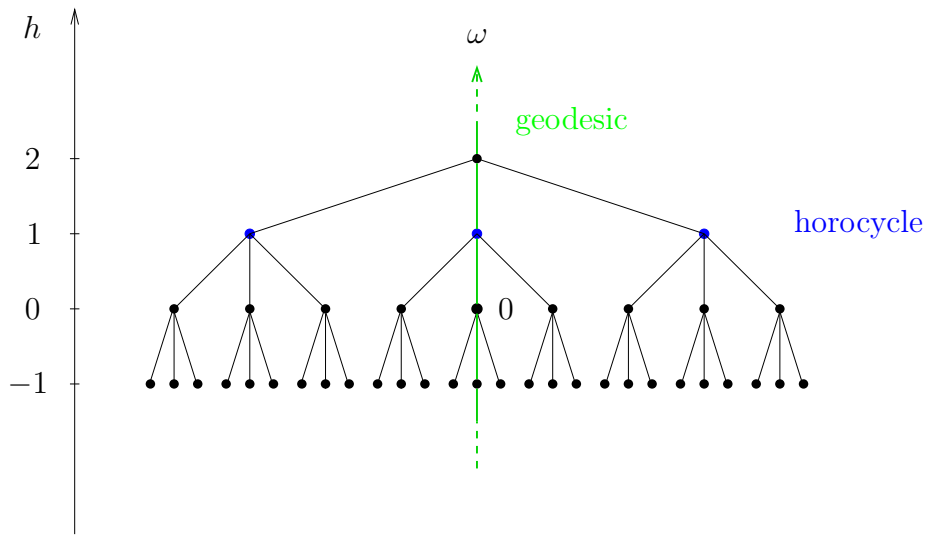}
\caption{Upper half--space picture of \,$\mathbb{T}_3$}
\label{UpperHalfSpacePictureTree}
\end{figure}

A homogeneous tree is a connected graph,
with no loop and with the same number of edges at each vertex.
Let us denote by \ssf$\Tq$
the set of vertices of the homogeneous tree
with \ssf$q\ssb+\!1\!>\!2$ \ssf edges.
It is equipped with the counting measure and with the geodesic distance,
given by the number of edges between two points.
The volume of any sphere \ssf$S(x,r)$ of radius \ssf$r\hspace{-.75mm}\in\!\N$
\ssf is given by
\vspace{-2mm}
\begin{equation*}
\delta(r)\ssf=\,\begin{cases}
\;1
&\text{if \,}r\!=\ssb0\ssf,\\
\,(\ssf q\ssb+\!1)\,q^{\,r-1}
&\text{if \,}r\hspace{-.75mm}\in\!\N^*.
\end{cases}\end{equation*}
Once we have chosen an origin \ssf$0\!\in\!\Tq$ and
an oriented geodesic \ssf$\omega\ssb:\ssb\Z\ssb\to\!\Tq$ \ssf through \ssf$0$,
let us denote by $|x|\hspace{-.7mm}\in\!\N$
\ssf the distance of a vertex \ssf$x\!\in\!\Tq$ to the origin
and by \ssf$h(x)\hspace{-.6mm}\in\!\Z$ \ssf its horocyclic height
(see Figure \ref{UpperHalfSpacePictureTree}).
Let \ssf$G$ \ssf be the isometry group of \ssf$\mathbb{T}_q$
and let \ssf$K$ \ssf be the stabilizer of \ssf$0$\ssf.
Then \ssf$G$ \ssf is a locally compact group,
\ssf$K$ \ssf is a compact open subgroup,
and \ssf$\mathbb{T}_q\approx G/K$\ssf.

\begin{figure}[ht]
\psfrag{0}[c]{$0$}
\psfrag{sphere}[l]{\color{red}\text{sphere}}
\psfrag{blank}[l]{}
\includegraphics[height=70mm]{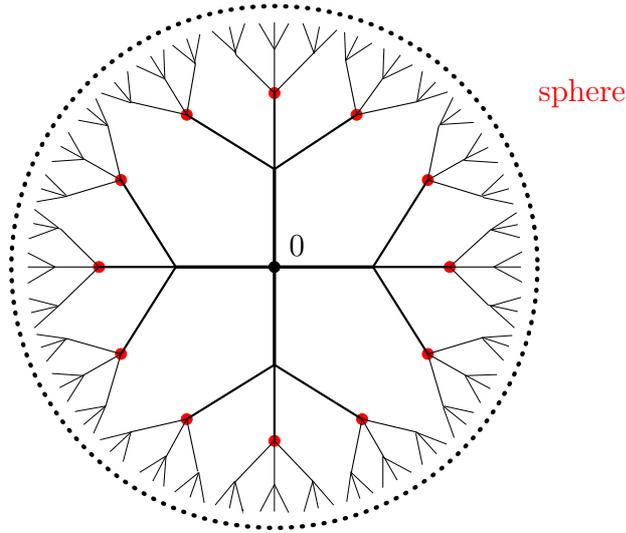}
\caption{Ball picture of \,$\mathbb{T}_3$}
\end{figure}

\begin{remark}
If \ssf$q$ \ssf is a prime number,
then \ssf$\Tq\ssb\approx\text{\rm PGL}(2,\Q_{\ssf q})
\ssf/\,\text{\rm PGL}(2,\Z_{\ssf q})$\ssf.
\end{remark}

The combinatorial Laplacian \ssf$\Delta$ \ssf on \ssf$\Tq$
is defined by \ssf$\Delta f\ssb=\ssb Af\ssb-\ssb f$\ssf,
where \ssf$A$ \ssf denotes the average operator
\vspace{-.5mm}
\begin{equation*}\label{AverageOperatorTq}
A f(x)=\tfrac1{q\,+\ssf1}\,\sum\nolimits_{\,y\in S(x,1)}f(y)\,.
\end{equation*}

\begin{remark}
Notice that $f$ is harmonic,
i.e., $\Delta f\!=\ssb0$
\ssf if and only if $f$ has the mean value property,
i.e., $f\!=\!Af$.
\end{remark}

The spherical functions \ssf$\phi_\lambda$ on \ssf$\Tq$
are the normalized radial eigenfunctions of \ssf$\Delta$\ssf,
or equivalently \ssf$A\ssf$.
Specifically,
\begin{equation*}\begin{cases}
\;A\,\phi_\lambda=\gamma(\lambda)\,\phi_\lambda\,,\\
\;\phi_\lambda(0)=1\,,
\end{cases}\end{equation*}
\vspace{-2.5mm}

\noindent
where \,$\gamma(\lambda)\ssb
=\ssb\frac{q^{\ssf i\lambda}\ssf+\,q^{-i\lambda}}{q^{\ssf1/2}\ssf+\,q^{-1/2}}$\,.

\begin{remark}
The spherical functions on \,$\Tq$ can be characterized again in several other ways.
Notably,
\begin{itemize}

\item
{\rm Explicit expression\,:}
\begin{equation}\label{SphericalFunctionsTq}
\phi_\lambda(r)=\begin{cases}
\,\mathbf{c}\vsf(\lambda)\,q^{\,(-1/2\ssf+\ssf i\ssf\lambda)\ssf r}\ssb
+\mathbf{c}\vsf(\ssb-\lambda)\,q^{\,(-1/2\ssf-\ssf i\ssf\lambda)\ssf r}
&\text{if \,}\lambda\!\in\hspace{-.5mm}\C\!\smallsetminus\!\frac\tau2\ssf\Z\ssf,\\
\,(-1)^{\ssf j\ssf r}\,\bigl(1\ssb
+\ssb\frac{q^{\ssf1/2}\ssf-\,q^{-1/2}}{q^{\ssf1/2}\ssf+\,q^{-1/2}}\,r\bigr)
\,q^{\ssf-\ssf r/2}
&\text{if \,}\lambda\ssb=\ssb\frac\tau2\,j\ssf,\\
\end{cases}\end{equation}
where \,$\mathbf{c}\vsf(\lambda)\ssb=\ssb
\frac1{q^{\ssf1/2}\ssf+\,q^{-1/2}}
\,\frac{q^{\ssf1/2\ssf+\ssf i\ssf\lambda}\ssf-\,q^{-1/2\ssf-\ssf i\ssf\lambda}}
{q^{\,i\ssf\lambda}\ssf-\,q^{-i\ssf\lambda}}$
\,and \,$\tau\ssb=\ssb\frac{2\pi}{\log q}$\,.
Notice that \eqref{SphericalFunctionsTq} is even and
\linebreak
\,$\tau$\ssf--\ssf periodic in \ssf$\lambda$\ssf.

\item
Same functional equations
as \eqref{FunctionalEquation1Sn}
and \eqref{FunctionalEquation2Sn}.

\item
The spherical functions are
the bi\ssf--\ssf$K$\!--\ssf invariant matrix coefficients of
the spherical principal series representations of \,$G$
on \ssf$L^2(\partial\,\Tq)$\ssf.

\end{itemize}
\end{remark}

The spherical Fourier transform of radial functions on \ssf$\Tq$
is defined by
\begin{equation}\label{SphericalFourierTransformTq}
\mathcal{H}f(\lambda)
=\ssf\sum\nolimits_{\,x\in\Tq}\ssb f(x)\,\phi_\lambda(x)
=f(0)+\tfrac{\,q^{1/2}\,+\;q^{-1/2}}{q^{1/2}}\,
\sum\nolimits_{\,r=1}^{+\infty}q^{\,r}\ssf f(r)\,\phi_\lambda(r)
\end{equation}
and the inversion formula reads
\begin{equation}\label{InverseSphericalFourierTransformTq}
f(r)\ssf=\ssf
\tfrac{q^{\ssf1/2}}{q^{\ssf1/2}\ssf+\,q^{-1/2}}\,\tfrac1\tau
\int_{\,0}^{\hspace{.5mm}\tau/2}\hspace{-1.5mm}
d\lambda\hspace{1mm}|\ssf\mathbf{c}\vsf(\lambda)\vsf|^{-2}\hspace{1mm}
\mathcal{H}f(\lambda)\,\phi_\lambda(r)\,.
\end{equation}
Consider the Abel transform
\vspace{1mm}
\begin{equation*}\label{AbelTransformTq}\begin{aligned}
\mathcal{A}f(h)
&=\ssf q^{\frac h2}\ssf
\sum\nolimits_{\hspace{-1.mm}\substack{
\vphantom{o}\\x\in\Tq\\h(x)=h}}\hspace{-1.5mm}
f(|x|)\\
&=\,q^{\frac{|h|}2}\ssf f(|h|)
+\tfrac{q\,-\ssf1}q\,
\smash{\sum\nolimits_{\ssf j=1}^{+\infty}}\ssf
q^{\frac{|h|}2+\ssf j}\,
f(|h|\!+\ssb2\ssf j\ssf)\,,
\end{aligned}\end{equation*}
\vspace{-1mm}

which is essentially the horocyclic Radon transform restricted to radial functions.
Then the following commutative diagram holds,
let say in the Schwartz space setting\,:
\vspace{1mm}
\begin{equation*}\label{CommutativeDiagramTree}\begin{aligned}
&\hspace{-6.5mm}C_{\text{even}}^{\ssf\infty}\ssb(\R\ssf/\tau\ssf\Z)\\
\mathcal{H}\nearrow\hspace{3mm}
&\hspace{11mm}\nwarrow\mathcal{F}\\
\mathcal{S}_{\ssf\text{rad}}(\Tq)\hspace{6mm}
&\underset{\textstyle\mathcal{A}}\longrightarrow\hspace{6mm}
\mathcal{S}_{\ssf\text{even}}(\Z)
\end{aligned}\end{equation*}
Here \ssf$\mathcal{S}_{\text{even}}(\Z)$ denotes
the space of even functions on \ssf$\Z$ \ssf such that
\vspace{-1mm}
\begin{equation*}
\sup\nolimits_{\ssf r\in\N^*}\ssf
r^{\ssf m}\,|f(r)|<+\infty
\qquad\forall\;m\!\in\!\N\ssf,
\end{equation*}
\ssf$\mathcal{S}_{\ssf\text{rad}}(\Tq)$
the space of radial functions on \ssf$\Tq$\ssf,
whose radial part coincides with
\,$q^{-\frac r2}\ssf\mathcal{S}(\N)$\ssf,
\begin{equation*}
\mathcal{F}\ssb f(\lambda)=\sum\nolimits_{\ssf h\in\Z}
q^{\,i\ssf\lambda\ssf h}\,f(h)
\end{equation*}
is a variant of the Fourier transform on \ssf$\Z$\ssf,
and each arrow is an isomorphism.
The inverse Abel transform is given by
\begin{equation*}\label{InverseAbelTransformTq}
\mathcal{A}^{-1}\ssb f(r)\ssf
=\ssf\sum\nolimits_{\ssf j=0}^{+\infty}\ssf q^{-\frac r2-\ssf j}\ssf
\bigl\{\ssf f(r\ssb+\ssb2\ssf j)-f(r\ssb+\ssb2\ssf j\ssb+\ssb2)\bigr\}\,.
\end{equation*}
Consider finally the transform
\begin{equation*}
\mathcal{A}^*\ssb g\ssf(r)\ssf
=\ssf\tfrac1{\delta(r)}\,
\sum\nolimits_{\ssf\substack{
\vphantom{o}\\x\in\Tq\\|x|=\ssf r}}\ssb
q^{\frac{h(x)}2}\ssf g\bigl(h(x)\bigr)\,,
\end{equation*}
\vspace{-4mm}

which is dual to the Abel transform, i.e.,
\begin{equation*}
\sum\nolimits_{\ssf x\in\Tq}f(x)\,\mathcal{A}^*g(|x|)
=\sum\nolimits_{\ssf h\in\Z}\mathcal{A}f(h)\,g(h)\,,
\end{equation*}
\vspace{-4mm}

and which is an isomorphism between the space of all even functions on \ssf$\Z$
\ssf and the space of all radial functions on \ssf$\Tq$\ssf.
It is given explicitly by
\ssf$\mathcal{A}^*\ssb g\ssf(0)\ssb=\ssb0$\ssf,
\ssf and
\vspace{.5mm}
\begin{equation*}\label{DualAbelTransformTq}
\mathcal{A}^*\ssb g\ssf(r)
=\tfrac{2\,q}{q\,+\ssf1}\,q^{-\frac r2}f(r)
+\tfrac{q\,-\ssf1}{q\,+\ssf1}\,q^{-\frac r2}\,
\sum\nolimits_{\hspace{-6mm}\substack{
\vphantom{o}\\-r\,<\,h\,<\,r\\
\text{\rm$h$ \ssf has same parity as $r$}}}
\hspace{-6mm}f(h)
\end{equation*}
\vspace{-3mm}

if \ssf$r\hspace{-.75mm}\in\!\N^*$.
Its inverse is given by
\ssf$(\mathcal{A}^*)^{-1}\ssb f\ssf(0)=f\ssf(0)$\ssf,
\begin{equation*}\begin{aligned}
(\mathcal{A}^*)^{-1}\ssb f(h)\ssf
&=\tfrac{q^{\ssf1/2}\ssf+\,q^{-1/2}}2\,q^{\frac{h-1}2}f(h)\\
&\,-\ssf\tfrac{q\,-\,q^{-1}}2\,q^{-\frac h2}\ssf
\sum\nolimits_{\ssf0\ssf<\,r\;\text{\rm odd}\,<\,h}
\ssf q^{\,r}\ssf f(r)
\end{aligned}\end{equation*}
\vspace{-3mm}

if \,$h\!\in\!\N$ \ssf is odd, and
\begin{equation*}\begin{aligned}
(\mathcal{A}^*)^{-1}\ssb f(h)\ssf
&=\ssf\tfrac{q^{1/2}\ssf+\,q^{-1/2}}2\;q^{\frac{h-1}2}f(h)
-\tfrac{q^{1/2}\ssf-\,q^{-1/2}}2\;q^{-\frac{h-1}2}f(0)\\
&\,-\ssf\tfrac{q\,-\,q^{-1}}2\,q^{-\frac h2}\ssf
\sum\nolimits_{\ssf0\ssf<\,r\;\text{\rm even}\,<\,h}
\ssf q^{\,r}\ssf f(r)
\end{aligned}\end{equation*}
\vspace{-3mm}

if \,$h\!\in\!\N^*$ is even.
\smallskip

\paragraph{\bf Applications.}
Let us use spherical Fourier analysis
to study discrete evolution equations on \ssf$\Tq$\ssf,
as we did in the differential setting on \ssf$\Hn$\ssf.
\smallskip

$\bullet$
\,Consider the heat equation with discrete time \ssf$t\!\in\!\N$
\begin{equation*}\label{HeatEquationTq}\begin{cases}
\,u(x,t\hspace{-.4mm}+\!1)-u(x,t)=\Delta_{\ssf x\ssf}u(x,t)\\
\,u(x,0)=f(x)
\end{cases}\end{equation*}
or, equivalently, the simple random walk
\begin{equation*}
u(x,t)=A^{\ssf t\ssb}f(x)=f\!*\ssb h_{\ssf t}\ssf(x)\,.
\end{equation*}
Its density \ssf$h_{\ssf t}(x)$ vanishes
unless \ssf$|x|\!\le\ssb t$ \ssf have the same parity.
In this case,
\begin{equation}\label{HeatKernelEstimateTq}
h_{\ssf t\ssf}(x)\asymp
\tfrac{1\ssf+\ssf|x|\vphantom{\frac00}}
{(1\ssf+\ssf t)\,\sqrt{\ssf1\ssf+\ssf t\ssf-\ssf|x|\ssf}}\;
\gamma_{\ssf0}^{\ssf t}\;q^{-\frac{|x|}2}\,
e^{-\ssf t\,\psi\ssf\left(\frac{1\ssf+\ssf|x|}{1\ssf+\,t\;}\right)}\,,
\end{equation}
\vspace{-4mm}

\noindent
where \,$\psi(z)\ssb
=\ssb\frac{1\ssf+\ssf z}2\ssf\log(1\!+\ssb z)\ssb
-\ssb\frac{1\ssf-\ssf z}2\ssf\log(1\!-\ssb z)$
\hspace{.5mm}and \,$\gamma_{\ssf0}\ssb=\ssb\gamma(0)\ssb=\ssb
\frac{2\vphantom{|}}{q^{\ssf1/2}\ssf+\,q^{-1/2}}\ssb<\ssb1$
\ssf is the spectral radius of \ssf$A$\ssf.
\smallskip

$\bullet$
\,The shifted wave equation with discrete time \ssf$t\!\in\!\Z$
\begin{equation*}\label{ShiftedWaveEquationTree}\begin{cases}
\,\gamma_{\ssf0}\,\Delta_{\ssf t}^{\Z}\ssf u(x,t)
=\bigl(\Delta_{\ssf x}^{\ssb\Tq}\!
+\!1\!-\ssb\gamma_{\ssf0}\bigr)\ssf u(x,t)\\
\,u(x,0)\ssb=\ssb f(x)\ssf,\;
\{\ssf u(x,1)-u(x,-1)\}/\ssf2\ssb=\ssb g(x)\\
\end{cases}\end{equation*}
can be solved explicitly by using the inverse dual Abel transform.
Specifically,
\begin{equation*}
u(x,t)=C_{\ssf t\ssf}f(x)+S_{\ssf t\ssf}\ssf g\ssf(x)\,,
\end{equation*}
where
\begin{equation*}\begin{cases}
\,C_{\ssf t\ssf}=\frac{M_{\ssf|t|}\ssf-\,M_{\ssf|t|-2}}{2\vphantom{|}}\\
\,S_{\ssf t\ssf}=\sign(t)\,M_{\ssf|t|-1}\\
\end{cases}\end{equation*}
and
\begin{equation*}
M_{\ssf t\ssf}f\ssf(x)
=\ssf q^{-\frac t2}\ssf
\sum\nolimits_{\hspace{-2.75mm}\substack{\vphantom{o}\\
d\ssf(y,\ssf x)\ssf\le\,t\\
t\ssf-\ssf d\ssf(y,\ssf x)
\hspace{1mm}\text{\rm even}}}
\hspace{-2.5mm}f(y)
\end{equation*}
\vspace{-4mm}

\noindent
if \hspace{.5mm}$t\!\ge\!0$\ssf,
while \hspace{.5mm}$M_{\ssf t}\!=\ssb0$
\hspace{.5mm}if \hspace{.5mm}$t\!<\!0$\ssf.

\subsection{Comments, references and further results}
\label{CommentsRankOne}

\begin{itemize}

\item
Our main reference for classical special functions is \cite{DLMF}.
\smallskip

\item
The spherical Fourier analysis presented in this section
takes place on homogeneous spaces \ssf$G/K$
\ssf associated with Gelfand pairs \ssf$(G,K)$\ssf.
\smallskip

\item
The sphere \ssf$\Sn$ and the hyperbolic space \ssf$\Hn$ are dual symmetric spaces.
By letting their curvature tend to \ssf$0$\ssf,
the Euclidean space \ssf$\Rn$ is obtained as a limit.
At the level of spherical functions,
the duality is reflected by the relation
\begin{equation*}
\phi_{\,\ell}^{\ssf\Sn}\!(\cos\theta_1)
=\phi_{\,\lambda}^{\ssf\Hn}\!(r)
\end{equation*}
between \eqref{SphericalFunctionsSn} and \eqref{SphericalFunctionsHn},
when we specialize
\ssf$\lambda\ssb=\pm\,i\,(\ssf\rho\ssb+\ssb\ell\,)$
\ssf(\ssf$\ell\!\in\!\N$\ssf)
and take \,$r\ssb=\ssb\pm\, i\,\theta_1$\ssf.
And \eqref{SphericalFunctionsRn} is a limit
of \eqref{SphericalFunctionsHn} and \eqref{SphericalFunctionsSn}\,:
\begin{equation*}
\phi_{\,\lambda}^{\ssf\Rn}\!(r)=
\lim\nolimits_{\,\epsilon\to0}\ssf
\phi_{\ssf\frac\lambda\epsilon}^{\ssf\Hn}\ssb(\epsilon\ssf r)
=\lim\nolimits_{\,\ell\to+\infty}\ssf
\phi_{\,\ell}^{\ssf\Sn}\!(\cos\tfrac{\lambda\ssf r\vphantom{|}}{\ell})\,.
\end{equation*}
Homogeous trees may be considered as discrete analogs of hyperbolic spaces.
This is for instance justified by the structural similarity between
\,$\mathbb{H}^2\ssb\approx\text{\rm PSL}(2,\R)\ssf/\,\text{\rm PSO}(2)$
\ssf and
\,$\Tq\approx\text{PGL}(2,\Q_{\ssf q})\ssf/\,\text{PGL}(2,\Z_{\ssf q})$
\ssf when \ssf$q$ \ssf is a prime number.
At the analytic level, an actual relation is provided
by the meta--theory developed by Cherednik \cite{Cherednik2005},
which includes as limit cases
spherical Fourier analysis on \ssf$\Tq$ and on \ssf$\Hn$,
as well as on \ssf$\Rn$ or on \ssf$\Sn$.
\smallskip

\item
The material in Subsection \ref{SphericalFourierAnalysisHn} generalizes to the
class of Riemannian symmetric spaces of noncompact type and of rank \ssf$1$\ssf,
which consist of all hyperbolic spaces 
\begin{equation*}
\Hn\!=\ssb\text{\rm H}^n(\R)\ssf,
\,\text{\rm H}^n(\C)\ssf,
\,\text{\rm H}^n(\mathbb{H})\ssf,
\,\text{\rm H}^2(\mathbb{O})\ssf,
\end{equation*}
and further to the class of Damek--Ricci spaces.
One obtains this way a group theoretical interpretation of Jacobi functions
\ssf$\phi_\lambda^{\alpha,\beta}$
for infinitely many discrete parameters
\begin{equation*}
\alpha=\tfrac{m_1\ssf+\,m_2\ssf-\ssf1}2\ssf,
\,\beta=\tfrac{m_2\ssf-\ssf1}2\ssf.
\end{equation*}
Our main references are
\cite{Faraut1982},
\cite{Koornwinder1984},
\cite{AnkerDamekYacoub1996},
\cite{Rouviere1999}
and \cite{AnkerMartinotPedonSetti2013}.
\smallskip

\item
Our main references for Subsection \ref{SphericalFourierAnalysisTree} are
\cite{FigatalamancaPicardello1983},
\cite{CowlingMedaSetti1998},
and \cite{AnkerMartinotPedonSetti2013}.
Evolution equations (heat, Schr\"odinger, wave) with continuous time
were also considered on homogeneous trees
(see \cite{Setti1998},
\cite{MedollaSetti1999},
\cite{CowlingMedaSetti2000}
and \cite{Jamaleddine2013}).
\smallskip

\item
Spherical Fourier analysis generalizes to higher rank
(see \cite{Helgason1984},
\cite{GangolliVaradarajan1988}
for Riemannian symmetric spaces
and \cite{Macdonald1971},
\cite{Parkinson2006},
\cite{ManteroZappa2011}
for affine buildings).
\vspace{1mm}

\noindent
\textit{Classification}\,:
\vspace{2mm}

\begin{center}
\begin{tabular}{|c||c|c|c|}
\hline
$\vphantom{\Big|}$
\hspace{11mm} type \hspace{11mm}
& \ssf constant curvature\ssf
& \hspace{10mm} rank $1$ \hspace{10mm}
& \hspace{5mm} general case \hspace{5mm} \\
\hline\hline
$\vphantom{\Big|}$
Euclidean
& $\color{red}\mathbb{R}^n$
& 
& $\mathfrak{p}\!\rtimes\!K/K$ \\
\hline
$\vphantom{\Big|}$
compact
& $\color{red}\Sn\!=\text{S}(\mathbb{R}^{n+1})$
& $\text{S}(\mathbb{F}^n)$
& $U/K$ \\
\hline
$\vphantom{\Big|}$
non compact
& $\color{red}\mathbb{H}^n\!=\text{H}^n(\mathbb{R})$
& $\text{H}^n(\mathbb{F})$
& $G/K$ \\
\hline
$\vphantom{\Big|}$
$p$\ssf--\ssf adic
&
& {\color{red}homogeneous trees}
& affine buildings \\
\hline
\end{tabular}
\end{center}
\vspace{2mm}

\noindent
\textit{Notation}\,:
\begin{itemize}

\item[$\circ$]
$\ssf\mathfrak{g}_{\vsf\mathbb{C}}$
\ssf is a complex semisimple Lie algebra,

\item[$\circ$]
$\ssf\mathfrak{g}$ \,is a noncompact real form
of \ssf$\mathfrak{g}_{\vsf\mathbb{C}}$
and \ssf$\mathfrak{g}\ssb=\ssb\mathfrak{k}\oplus\mathfrak{p}$
\ssf a Cartan decomposition of \ssf$\mathfrak{g}$\ssf,

\item[$\circ$]
$\ssf\mathfrak{u}\ssb=\ssb\mathfrak{k}\oplus i\ssf\mathfrak{p}$
\,is the compact dual form of \ssf$\mathfrak{g}$\ssf,

\item[$\circ$]
$\ssf\apartment$ \,is a Cartan subspace of \ssf$\mathfrak{p}$\ssf,

\item[$\circ$]
$\ssf\mathfrak{g}\ssb=\ssb\apartment\oplus\mathfrak{m}\oplus\ssb
\bigl(\oplus_{\vsf\alpha\ssf\in\vsf\mathcal{R}}\,\mathfrak{g}_{\vsf\alpha}\bigr)$
\ssf is the root space decomposition
of \ssf$(\mathfrak{g}\ssf,\ssb\apartment)$\ssf,

\item[$\circ$]
$\ssf\mathcal{R}^{\ssb+}$ is a positive root subsystem
and \ssf$\chamber$ the corresponding positive Weyl chamber
in \ssf$\apartment$\ssf,

\item[$\circ$]
$\ssf\mathfrak{n}\ssb=\ssb
\oplus_{\vsf\alpha\ssf\in\vsf\mathcal{R}^{\ssb+}}\ssf\mathfrak{g}_{\vsf\alpha}$
is the corresponding nilpotent Lie subalgebra,

\item[$\circ$]
$\ssf m_\alpha\!=\hspace{-.4mm}\dim\mathfrak{g}_{\vsf\alpha}$
is the multiplicity of the root \ssf$\alpha$\ssf,

\item[$\circ$]
$\ssf\rho\ssb=\ssb\frac12\displaystyle\sum\nolimits_{\ssf\alpha\ssf\in\vsf\mathcal{R}^{\ssb+}}\!m_\alpha$\ssf,

\item[$\circ$]
$\ssf G_{\vsb\mathbb{C}}$
is a complex Lie group with finite center and Lie algebra
\ssf$\mathfrak{g}_{\vsf\mathbb{C}}$\ssf,

\item[$\circ$]
$\ssf G$, $U$\ssb, $K$ \ssf and \ssf$N$
are the Lie subgroups of \ssf$G_{\vsb\mathbb{C}}$
corresponding to the Lie subalgebras 
\ssf$\mathfrak{g}$\ssf, $\mathfrak{u}$\ssf,
$\mathfrak{k}$ \ssf and \ssf$\mathfrak{n}$\ssf.

\end{itemize}
\vspace{1mm}

\noindent
\textit{Special functions}\,:
\vspace{1mm}

\begin{itemize}

\item[$\circ$]
Bessel functions on \ssf$\mathfrak{p}$\,:
\begin{equation}\label{BesselFunctionP}
\phi_{\ssf\lambda}^{\,\mathfrak{p}}(x)\ssf=\int_K\!dk\;
e^{\,i\,\langle\ssf\lambda\ssf,\ssf(\Ad k)\ssf x\ssf\rangle}
\qquad\forall\;\lambda\!\in\hspace{-.5mm}\mathfrak{p}_{\vsf\C}\ssf,
\;\forall\;x\!\in\hspace{-.5mm}\mathfrak{p}_{\vsf\C}\ssf.
\end{equation}

\item[$\circ$]
Spherical functions on \ssf$G/K$\,:
\begin{equation}\label{SphericalFunctionGK}
\phi_{\ssf\lambda}^{\ssf G}(x)\ssf=\int_K\!dk\;
e^{\,\langle\,i\ssf\lambda\ssf+\ssf\rho\ssf,\,\apartment\ssf(k\ssf x)\,\rangle}
\qquad\forall\;\lambda\!\in\!\apartment_{\ssf\C}\ssf,
\;\forall\;x\!\in\!G\ssf,
\end{equation}
where \ssf$\apartment\vsf(y)$ denotes
the \ssf$\apartment$\ssf--\ssf component of \ssf$y\!\in\!G$
\ssf in the Iwasawa decomposition
\,$G\ssb=\ssb N(\exp\apartment)K$\ssf.
\vspace{1mm}

\item[$\circ$]
Spherical functions on affine buildings are classical Macdonald polynomials
(i.e. Hall--Littlewood polynomials for the type
\ssf$\smash{\widetilde{A}}$\hspace{.5mm}).
\end{itemize}
\smallskip

\item
The global heat kernel estimate \eqref{HeatKernelEstimateHn}
generalizes as follows to \ssf$G/K$\,:
\begin{equation}\label{HeatKernelEstimateGK}\begin{aligned}
h_{\ssf t\ssf}(x)\ssf&\asymp\,t^{\ssf-\frac n2}\,
\Bigl\{\,\prod\nolimits_{\ssf\alpha\ssf\in\vsf R^+}
\bigl(\ssf1\!+\ssb\langle\alpha\ssf,\ssb x^{\ssb+}\rangle\vsf\bigr)
\bigl(\ssf1\!+\ssb t\ssb+\ssb\langle\alpha\ssf,\ssb x^{\ssb+}\rangle\vsf
\bigr)^{\frac{m_\alpha+\ssf m_{2\alpha}}2}\Bigr\}\\
&\ssf\times\,e^{-\ssf|\rho|^2\ssf t\,
-\,\langle\ssf\rho\vsf,\ssf x^+\rangle\,-\,\frac{|x^{\ssb+}\ssb|^2}{4\,t}}\,,
\end{aligned}\end{equation}
where \ssf$R^{\ssf+}$ is the set of positive indivisible roots in \ssf$\mathcal{R}$
\ssf and \ssf$x^+$ denotes the
\ssf$\overline{\smash{\chamber}\vphantom{X}}$--\,component
of \ssf$x\!\in\!G$ \ssf in the Cartan decomposition
\,$G\ssb=\ssb K(\vsf\exp\overline{\smash{\chamber}\vphantom{X}}\vsf)\ssf K$
(see \cite{AnkerOstellari2004} and the references therein).
\smallskip

\item
Random walks are harder to analyze on affine buildings.
A global estimate similar to
\eqref{HeatKernelEstimateTq} and \eqref{HeatKernelEstimateGK}
was established in \cite{AnkerSchapiraTrojan2006}
for the simple random walk on affine buildings of type
\ssf$\smash{\widetilde{\text{\rm A}}_{\ssf2}}$\ssf.
In general, the main asymptotics of random walks were obtained in \cite{Trojan2013}.

\end{itemize}

\subsection{Epilogue}

This section outlines spherical Fourier analysis around 1980
(except for the later applications to evolution equations).
In the 1980s, Heckman and Opdam addressed the following problem
(which goes back to Koornwinder for the root system \ssf$\text{BC}_{\ssf2}$)\,:
for any root system \ssf$R$\ssf,
construct a continuous family of special functions on~\ssf$\apartment$ \ssf
generalizing spherical functions on the corresponding symmetric spaces \ssf$G/K$\ssf,
as Jacobi functions (or equivalently the Gauss hypergeometric function)
generalize spherical functions on hyperbolic spaces.
This problem was solved during the 1990s,
mainly by Cherednik, Dunkl, Heckman, Macdonald, and Opdam,
and has actually given rise to a large theory
of special functions associated to root system,
which is nowadays often referred to as Dunkl theory.

\begin{remark}
A different generalization of hypergeometric functions to Grassmannians
was developed by Aomoto and Gelfand at the end of the 20th century.
\end{remark}

\section{Rational Dunkl theory}
\label{RationalDunklTheory}

Rational Dunkl theory originates from the seminal paper \cite{Dunkl1989}.
This theory of special functions in several variables encompasses
\begin{itemize}

\item
Euclidean Fourier analysis
(which corresponds to the multiplicity \ssf$k\ssb=\ssb0$\ssf)\ssf,

\item
classical Bessel functions in dimension \ssf$1$\ssf,

\item
generalized Bessel functions associated with
Riemannian symmetric spaces of Euclidean type
(which correspond to a discrete set of multiplicities \ssf$k$\ssf)\ssf.

\end{itemize}
In this subsection,
we use \cite{Roesler2002}
(or alternately \cite{Graczyk2008}, pp.\;1--69)
as our primary reference
and quote only later works.
Our notation goes as follows
(see the appendix for more details)\,:
\begin{itemize}

\item
$\,\apartment$ \,is a Euclidean vector space of dimension \ssf$n$\ssf,
which we identify with its dual space,

\item
$\ssf R$ \ssf is a root system,
which is reduced but not necessarily crystallographic,

\item
$\ssf W$ is the associated reflection group,

\item
$\,\chamber$ is a positive Weyl chamber in \ssf$\apartment$\ssf,
$\overline{\smash{\chamber}\vphantom{X}}$ its closure,
$R^+$ the corresponding positive root subsystem,
and $S$ the subset of simple roots,

\item
$\,\overline{\apartment_+\vphantom{0}}$
denotes the closed cone generated by $R^+$,
which is the dual cone of \ssf$\overline{\smash{\chamber}\vphantom{X}}$\ssf,

\item
for every \ssf$x\!\in\!\apartment$\ssf,
\,$x^+$ denotes the element of the orbit \ssf$W\ssb x$ \ssf
which lies in \ssf$\overline{\smash{\chamber}\vphantom{X}}$,

\item
$\ssf w_{\ssf0}$ \ssf denotes the longest element in \ssf$W$\ssb,
which interchanges \ssf$\smash{\chamber}$ and \ssf$\smash{-\ssf\chamber}$,
respectively \ssf$R^+$ and \ssf$R^-\!=\ssb-\ssf R^+$,

\item
$\,k$ \,is a multiplicity,
which will remain implicit in most formulae
and which will be assumed to be nonnegative after a while,

\item
$\ssf\gamma
=\ssb{\displaystyle\sum\nolimits_{\ssf\alpha\in R^+}}\ssb k_\alpha$\ssf,
\vspace{.5mm}

\item
$\ssf\delta(x)\ssb=\ssb\displaystyle\prod\nolimits_{\ssf\alpha\in R^+}\!
|\ssf\langle\ssf\alpha\ssf,x\ssf\rangle\ssf|^{\ssf2\ssf k_\alpha}$
\ssf is the reference density in the case \ssf$k\ssb\ge\ssb0$\ssf.

\end{itemize}
 
\subsection{Dunkl operators} 

\begin{definition}\label{DefinitionRationalDunklOperators}
The rational Dunkl operators, which are often simply called Dunkl operators,
are the differential--difference operators defined by
\begin{equation}\label{DunklOperators}
D_\xi f(x)=\ssf\partial_{\ssf\xi}f(x)
+\ssf\sum\nolimits_{\ssf\alpha\in R^+}\!
k_\alpha\,
\tfrac{\langle\ssf\alpha\ssf,\,\xi\ssf\rangle}{\langle\ssf\alpha\ssf,\,x\ssf\rangle}\,
\bigl\{\ssf f(x)\ssb-\ssb f(r_{\ssb\alpha\ssf}x)\ssf\bigr\}
\end{equation}
for every \,$\xi\!\in\!\mathfrak a$\ssf.
\end{definition}

\begin{remark}\label{RemarksDunklOperators}
${}$
\begin{itemize}

\item
Notice that Dunkl operators \,$D_\xi$
reduce to partial derivatives \,$\partial_{\ssf\xi}$
when \,$k\ssb=\ssb0$\ssf.

\item
The choice of \,$R^{\ssf+}$ plays no role
in Definition \ref{DefinitionRationalDunklOperators},
as
\begin{equation*}
\sum\nolimits_{\ssf\alpha\in R^+}\!k_\alpha\,
\tfrac{\langle\ssf\alpha\ssf,\,\xi\ssf\rangle}{\langle\ssf\alpha\ssf,\,x\ssf\rangle}\,
\bigl\{\ssf f(x)\ssb-\ssb f(r_{\ssb\alpha\ssf}x)\ssf\bigr\}
=\sum\nolimits_{\ssf\alpha\in R}\!\tfrac{k_\alpha}2\,
\tfrac{\langle\ssf\alpha\ssf,\,\xi\ssf\rangle}{\langle\ssf\alpha\ssf,\,x\ssf\rangle}\,
\bigl\{\ssf f(x)\ssb-\ssb f(r_{\ssb\alpha\ssf}x)\ssf\bigr\}\,.
\end{equation*}

\item
Dividing by \,$\langle\alpha,x\rangle$
\,produces no actual singularity in \eqref{DunklOperators},
as
\begin{equation*}
\tfrac{f(x)\ssf-\ssf f(r_{\alpha\ssf}x)}{\langle\ssf\alpha\ssf,\,x\ssf\rangle}
=-\,\tfrac1{\langle\ssf\alpha\ssf,\,x\ssf\rangle}\int_{\,0}^{\ssf1}\!dt\;
\tfrac{\partial}{\partial\ssf t}\ssf
f(x\ssb-\ssb t\,\langle\ssf\alpha^{\!\vee}\!,x\ssf\rangle\,\alpha)
=\tfrac2{|\alpha|^2}\int_0^1\!dt\,\partial_\alpha f(x\ssb
-\ssb t\,\langle\ssf\alpha^{\!\vee}\!,x\ssf\rangle\,\alpha)\,.
\end{equation*}

\end{itemize}
\end{remark}

Commutativity is a remarkable property of Dunkl operators.

\begin{theorem}
Fix a multiplicity \,$k$\ssf. Then
\begin{equation*}
D_\xi\circ D_\eta=D_\eta\circ D_\xi
\qquad\forall\;\xi,\eta\ssb\in\ssb\apartment\ssf.
\end{equation*}
\end{theorem}

This result leads to the notion of Dunkl operators \ssf$D_p$\ssf,
for every polynomial \ssf$p\hspace{-.5mm}\in\!\mathcal{P}(\apartment)$\ssf,
and of their symmetric part \ssf$\smash{\widetilde{D}}_p$
on symmetric (\ssf i.e.~$W$\ssb--\ssf invariant) functions.

\begin{example}
${}$
\begin{itemize}

\item
The Dunkl Laplacian is given by
\begin{align*}
\Delta f(x)
=\sum\nolimits_{\ssf j=1}^{\,n}\!D_j^{\ssf2}f(x)
&=\!\overbrace{\vphantom{\frac||}
\sum\nolimits_{\ssf j=1}^{\,n}\!\partial_j^{\ssf2}f(x)
+\ssf\sum\nolimits_{\ssf\alpha\in R^+}\!
\tfrac{2\ssf k_\alpha}{\langle\ssf\alpha,\ssf x\ssf\rangle}\,
\partial_\alpha f(x)}^{\textstyle\text{\rm differential part
\ssf$\smash{\widetilde{\Delta}f(x)}$}}\\
&\ssf-\underbrace{\sum\nolimits_{\ssf\alpha\in R^+}\!
\tfrac{k_\alpha\ssf|\alpha|^2}{\langle\ssf\alpha\ssf,\,x\ssf\rangle^{\ssf2}}\,
\bigl\{\ssf f(x)\ssb-\ssb f(r_{\ssb\alpha\ssf}x)\ssf\bigr\}
\vphantom{\frac||}}_{\textstyle\text{\rm difference part}}\,,
\end{align*}
where \,$D_j$\ssf, respectively \,$\partial_j$ denote
the Dunkl operators, respectively the partial derivatives
with respect to an orthonormal basis of \,$\apartment$\ssf.

\item
In dimension \ssf$1$,
the Dunkl operator is given by
\begin{equation*}
Df(x)=\bigl(\tfrac\partial{\partial x}\bigr)f(x)
+\tfrac kx\,\bigl\{\ssf f(x)\!-\!f(-x)\ssf\bigr\}
\end{equation*}
and the Dunkl Laplacian by
\begin{equation*}
Lf(x)=\bigl(\tfrac\partial{\partial x}\bigr)^{\ssb2}f(x)
+\tfrac{2\ssf k}x\,\bigl(\tfrac\partial{\partial x}\bigr)f(x)
-\tfrac k{x^2}\,\bigl\{\ssf f(x)\!-\!f(-x)\ssf\bigr\}\,.
\end{equation*}

\end{itemize}
\end{example}

Here are some other properties of Dunkl operators.

\begin{proposition}
${}$
\begin{itemize}

\item
The Dunkl operators map the following function spaces into themselves\,:
\begin{equation*}
\mathcal{P}(\apartment)\ssf,\;
\mathcal{C}^{\ssf\infty}(\apartment)\ssf,\;
\mathcal{C}_c^{\ssf\infty}(\apartment)\ssf,\;
\mathcal{S}(\apartment)\ssf,\;\dots
\end{equation*}
More precisely, the Dunkl operators \,$D_\xi$\ssf,
with \,$\xi\!\in\!\apartment$\ssf,
are homogeneous operators of degree \ssf$-1$ on polynomials.

\item
\text{\rm $W$\ssb--\ssf equivariance\,:}
For every \,$w\!\in\!W$ \ssb and
\,$p\hspace{-.5mm}\in\!\mathcal{P}(\apartment)$\ssf,
we have
\begin{equation*}
w\circ D_p\circ w^{-1}\ssb=D_{w\ssf p}\,.
\end{equation*}
Hence \,$\widetilde{D}_{p\ssf q}\hspace{-.5mm}
=\ssb\widetilde{D}_p\ssb\circ\ssb\widetilde{D}_q$\ssf,
for all symmetric $($i.e.~$W$\hspace{-.5mm}--\ssf invariant\,$)$ polynomials
\,$p\ssf,\ssb q\hspace{-.5mm}\in\!\mathcal{P}(\apartment)^W$.

\item
\text{\rm Skew--adjointness\,:}
Assume that \,$k\ssb\ge\ssb0$\ssf.
Then, for every \,$\xi\!\in\!\apartment$\ssf,
\begin{equation*}
\int_{\ssf\apartment}\ssb dx\;\delta(x)\,D_\xi f(x)\hspace{1mm}g(x)
=\ssf-\int_{\ssf\apartment}\ssb dx\;\delta(x)\,f(x)\,D_\xi\vsf g(x)\,.
\end{equation*}

\end{itemize}
\end{proposition}

\subsection{Dunkl kernel}

\begin{theorem}
For generic multiplicities \,$k$
and for every \,$\lambda\!\in\!\apartment_{\ssf\C}$\ssf,
the system
\begin{equation*}\begin{cases}
\,D_\xi\ssf E_\lambda\ssb=\langle\ssf\lambda\ssf,\xi\ssf\rangle\,E_\lambda
\quad\forall\;\xi\!\in\!\apartment\ssf,\\
\,E_\lambda(0)\ssb=\ssb1\ssf,
\end{cases}\end{equation*}
has a unique smooth solution on \,$\apartment$\ssf,
which is called the Dunkl kernel.
\end{theorem}

\begin{remark}
In this statement, generic means that
\,$k$ belongs to a dense open subset \,$K_{\text{\rm reg}}$ of \,$K$,
whose complement is a countable union of algebraic sets.
The set \,$K_{\text{\rm reg}}$ is known explicitly
and it contains in particular
\,$\{\,k\!\in\hspace{-1mm}K\,|\,\Re k\ssb\ge\ssb0\,\}$\ssf.
\end{remark}

\begin{definition}
The generalized Bessel function is the average
\begin{equation}\label{BesselFunctionDunkl}
J_\lambda(x)
=\tfrac{1\vphantom{\frac00}}{|W|\vphantom{\frac00}}
\sum\nolimits_{\ssf w\in W}E_\lambda(w\ssf x)
=\tfrac{1\vphantom{\frac00}}{|W|\vphantom{\frac00}}
\sum\nolimits_{\ssf w\in W}E_{\ssf w\vsf \lambda}(x)\,.
\end{equation}
\end{definition}

\begin{remark}
${}$
\begin{itemize}

\item
Mind the possible formal confusion between \eqref{BesselFunctionDunkl}
and the classical Bessel function of the first kind \,$J_\nu$\ssf.

\item
Conversely, the Dunkl kernel \,$E_\lambda(x)$
can be recovered by applying to the generalized Bessel function \,$J_\lambda(x)$
 a linear differential operator in \,$x$
 whose coefficients are rational functions of \ssf$\lambda$
$($see \cite[proposition 6.8.(4)\vsf]{Opdam1993}$)$.

\item
In dimension \ssf$1$,
\ssf$K_{\text{\rm reg}}$ is the complement
of \ssf$-\ssf\N\ssb-\ssb\frac12$ in \,$\C$\ssf.
The generalized Bessel function \eqref{BesselFunctionDunkl} reduces to
the modified Bessel function encountered in Subsection~\ref{Hankel}$\,:$
\begin{equation*}
J_\lambda(x)=j_{\ssf k-\frac12}(\lambda\ssf x)\,.
\end{equation*}
The Dunkl kernel is a combination of two such functions$\,:$
\begin{equation*}
E_\lambda(x)=\underbrace{
j_{\ssf k-\frac12}(\lambda\ssf x)\vphantom{\tfrac00}
}_{\text{\rm even}}+\underbrace{
\tfrac{\lambda\ssf x}{2\ssf k+1}\,j_{\ssf k+\frac12}(\lambda\ssf x)
}_{\text{\rm odd}}
\end{equation*}
It can be also expressed in terms of the confluent hypergeometric function\,$:$
\begin{equation*}\begin{aligned}
E_\lambda(x)
&=\tfrac{\Gamma(k\ssf+\frac12)}{\sqrt{\pi\ssf}\,\Gamma(k)}\ssf
\int_{-1}^{+1}\hspace{-1mm}du\,
(1\!-\ssb u)^{k-1}\ssf(1\!+\ssb u)^k\,e^{\,\lambda\ssf x\ssf u}\\
&=\ssf e^{\,\lambda\ssf x}\underbrace{
\tfrac{\Gamma(2\ssf k\ssf+1)}{\Gamma(k)\,\Gamma(k\ssf+1)}\ssf
\int_{\ssf0}^{\,1}\!dv\hspace{.7mm}
v^{\ssf k-1}\ssf(1\!-\ssb v)^k\,e^{-2\ssf\lambda\ssf x\ssf v}
}_{\textstyle{}_1\vsf\text{\rm F}_{\vsb1}
(k\ssf;\ssb 2\ssf k\!+\!1\ssf;\ssb-\ssf2\ssf\lambda\ssf x)}\ssf.
\end{aligned}\end{equation*}

\item
When \,$k\ssb=\ssb0$\ssf,
the Dunkl kernel \,$E_\lambda(x)$ reduces to the exponential
\,$e^{\,\langle\ssf\lambda\ssf,\,x\ssf\rangle}$
and the generalized Bessel function \,$J_\lambda(x)$ to
\begin{equation}\label{DefinitionCosh}
\Cosh_{\vsf\lambda}(y)
=\tfrac{1\vphantom{\frac00}}{|W|\vphantom{\frac00}}
\sum\nolimits_{\ssf w\in W}e^{\,\langle\ssf w\lambda\ssf,\,y\ssf\rangle}\,.
\end{equation}

\item
As far as we know,
the non--symmetric Dunkl kernel had not occurred previously
in special functions, group theory or geometric analysis.

\item
Bessel functions \eqref{BesselFunctionP}
on Riemannian symmetric spaces of Euclidean type
\,$\mathfrak{p}\!\rtimes\!K/K$
are special cases of \eqref{BesselFunctionDunkl},
corresponding to crystallographic root systems
and to certain discrete sets of mutiplicities.
More precisely, if
\begin{itemize}

\item[$\circ$]
$\,\mathfrak{g}\hspace{-.4mm}=\ssb\mathfrak{p}\ssb\oplus\ssb\mathfrak{k}$
\ssf is the associated semisimple Lie algebra,

\item[$\circ$]
$\,\apartment$ is a Cartan subspace of \,$\mathfrak{p}$\ssf,

\item[$\circ$]
$\ssf\mathcal{R}$ is the root system of \,$(\mathfrak{g},\apartment)$\ssf,

\item[$\circ$]
$\,m_{\ssf\alpha}$ is the multiplicity of \,$\alpha\!\in\!\mathcal{R}$\ssf,

\end{itemize}
and
\begin{itemize}

\item[$\circ$]
$\ssf R$ is the subsystem of indivisible roots in \ssf$\mathcal{R}$\ssf,

\item[$\circ$]
$\,k_\alpha\!=\ssb\frac{m_{\ssf\alpha}+\,m_{\ssf2\ssf\alpha}\vphantom{|}}2$
\;$\forall\,\alpha\!\in\hspace{-1mm}R$\ssf,

\end{itemize}
then
\vspace{-1mm}
\begin{equation*}
\phi_{\ssf\lambda}^{\,\mathfrak{p}}(x)=J_{\ssf i\ssf\lambda}(x)
\qquad\forall\;\lambda\!\in\!\apartment_{\ssf\C}\ssf,
\;\forall\;x\!\in\!\apartment\ssf.
\end{equation*}

\end{itemize}
\end{remark}

In next proposition, we collect some properties of the Dunkl kernel.

\begin{proposition}
${}$
\begin{itemize}

\item
{\rm Regularity\,:}
\ssf$E_\lambda(x)$ extends to a holomorphic function
in \,$\lambda\!\in\!\apartment_{\ssf\C}$\ssf,
$x\!\in\!\apartment_{\ssf\C}$
and \,$k\!\in\hspace{-1mm}K_{\text{\rm reg}}$\ssf.

\item
{\rm Symmetries\,:}
\vspace{-1mm}
\begin{equation*}\begin{cases}
\,E_\lambda(x)\ssb=\ssb E_x(\lambda)\ssf,\\
\,E_{\ssf w\vsf\lambda}(w\ssf x)\ssb=\ssb E_\lambda(x)
&\forall\;w\!\in\!W\ssb,\\
\,E_\lambda(t\ssf x)\ssb=\ssb E_{\ssf t\vsf\lambda}(x)
&\forall\;t\ssb\in\ssb\C\ssf,\\
\,\overline{E_\lambda(x)}\ssb
=\ssb E_{\ssf\overline{\lambda}\ssf}(\overline{x})
&\text{when \,}k\ssb\ge\ssb0\ssf.
\end{cases}\end{equation*}
\vspace{-1mm}

\item
{\rm Positivity\,:}
Assume that \,$k\ssb\ge\ssb0$\ssf.
Then,
\begin{equation*}
0<E_\lambda(x)\le
e^{\,\langle\ssf\lambda^{\ssb+}\!,\,x^+\rangle}
\qquad\forall\;\lambda\!\in\!\apartment\ssf,
\,\forall\;x\!\in\!\apartment\ssf.
\end{equation*}

\item
{\rm Global estimate\,:}
Assume that \,$k\ssb\ge\ssb0$\ssf.
Then, for every \,$\xi_1,\dots,\xi_N\!\in\!\apartment$\ssf,
\begin{equation*}
\bigl|\ssf\partial_{\xi_1}\ssb\dots\,\partial_{\xi_N}\ssf E_\lambda(x)\ssf\bigr|
\le|\vsf\xi_1|\ssf\dots\ssf|\vsf\xi_N|\,|\lambda|^N\,
e^{\,\langle\ssf(\Re\ssb\lambda)^{\ssb+}\ssb,\,(\Re x)^{\ssb+}\ssf\rangle}\qquad\forall\;\lambda\!\in\!\apartment_{\ssf\C}\ssf,
\,\forall\;x\!\in\!\apartment_{\ssf\C}\ssf.
\end{equation*}

\end{itemize}
\end{proposition}

\subsection{Dunkl transform}

From now on, we assume that \ssf$k\ssb\ge\ssb0$\ssf.

\begin{definition}
The Dunkl transform is defined by
\begin{equation}\label{DunklTransform}
\mathcal{H}f(\lambda)
=\int_{\ssf\apartment}\ssb dx\;
\delta(x)\,f(x)\,E_{-i\lambda}(x)\,.
\end{equation}
\end{definition}

In next theorem, we collect the main properties of the Dunkl transform.

\pagebreak

\begin{theorem}
${}$
\begin{itemize}

\item
The Dunkl transform is an automorphism
of the Schwartz space \,$\mathcal{S}(\apartment)$\ssf.

\item
We have
\begin{equation*}\begin{cases}
\,\mathcal{H}(D_\xi f)(\lambda)\ssb
=i\,\langle\ssf\xi,\lambda\ssf\rangle\,\mathcal{H}f(\lambda)\ssb
&\forall\;\xi\!\in\!\apartment\ssf,\\
\,\mathcal{H}\ssf(\ssf\langle\ssf\xi,\ssf.\,\rangle\,f\ssf)\ssb
=i\,D_\xi\ssf\mathcal{H}f
&\forall\;\xi\!\in\!\apartment\ssf,\\
\,\mathcal{H}(wf)(\lambda)\ssb
=\mathcal{H}f(w\lambda)\ssb
&\forall\;w\!\in\!W,\\
\,\mathcal{H}\ssf\bigl[\ssf f(\vsf t\,.\,)\ssf\bigr](\lambda)
=|t|^{-\ssf n-2\ssf\gamma}\,(\mathcal{H}f)(\vsf t^{-1}\ssb\lambda)
&\forall\;t\!\in\!\R^*.\\
\end{cases}\end{equation*}
\vspace{-1mm}

\item
{\rm Inversion formula\,:}
\begin{equation}\label{InverseDunklTransform}
f(x)=\vsf\crat^{-2}\vsf\int_{\ssf\apartment}\ssb d\lambda\;
\delta(\lambda)\,\mathcal{H}f(\lambda)\,E_{\,i\lambda}(x)\,,
\end{equation}
\vspace{-4mm}

\noindent
where
\vspace{-1mm}
\begin{equation}\label{MehtaConstant}
\crat=\ssb\int_{\ssf\apartment}\ssb dx\;\delta(x)\,e^{-\frac{|x|^2}2}
\end{equation}
\vspace{-2mm}

\noindent
is the so--called Mehta--Macdonald integral.

\item
{\rm Plancherel identity\,:}
The Dunkl transform extends to an isometric
automorphism of \,$L^2(\apartment,\ssb\delta(x)\ssf dx)$\ssf,
up to a positive constant. Specifically,
\begin{equation*}\label{PlancherelFormulaDunkl}
\int_{\ssf\apartment}\ssb d\lambda\;
\delta(\lambda)\,|\ssf\mathcal{H}f(\lambda)\vsf|^{\ssf2}
=\vsf\crat^{\hspace{.5mm}2}\int_{\ssf\apartment}\ssb
dx\;\delta(x)\,|\vsf f(x)\vsf|^{\ssf2}\ssf.
\end{equation*}

\item
{\rm Riemann--Lebesgue Lemma\,:}
The Dunkl transform maps \,$L^1(\apartment,\ssb\delta(x)\ssf dx)$
into \,$\mathcal{C}_{\ssf0}(\apartment)$\ssf.

\item
{\rm Paley--Wiener Theorem\,:}
The Dunkl transform is an isomorphism
between \,$\mathcal{C}_{\ssf c}^{\ssf\infty}(\apartment)$
and the Paley--Wiener space \,$\mathcal{PW}(\apartment_{\ssf\C})$\ssf,
which consists of all holomorphic functions
\linebreak
\,$h\ssb:\ssb\apartment_{\ssf\C}\ssb\longrightarrow\ssb\C$
such that
\begin{equation}\label{PaleyWienerCondition}
\exists\;R\ssb>\ssb0\ssf,
\;\forall\,N\hspace{-1mm}\in\!\N\ssf,
\;\sup\nolimits_{\ssf\lambda\ssf\in\ssf\apartment_{\ssf\C}}\ssf
(1\!+\!|\lambda|)^N\ssf e^{\ssf-\ssf R\,|\ssb\Im\lambda\ssf|}\,
|\ssf h(\lambda)\ssf|<+\infty\,.
\end{equation}
More precisely, the support of
\,$f\!\in\ssb\mathcal{C}_{\ssf c}^{\ssf\infty}(\apartment)$
is contained in the closed ball \,$\smash{\overline{B(0,\ssb R)}}$
if and only if \,$h\ssb=\ssb\mathcal{H}f$ satisfies \eqref{PaleyWienerCondition}.

\end{itemize}
\end{theorem}

\begin{remark}
${}$
\begin{itemize}

\item
Notice that \eqref{DunklTransform}
and \eqref{InverseDunklTransform} are symmetric,
as the Euclidean Fourier transform \eqref{FourierTransformRn}
and its inverse \eqref{InverseFourierTransformRn},
or the Hankel transform \eqref{HankelTransform}
and its inverse \eqref{InverseHankelTransform}.

\item
In the \ssf$W$\hspace{-.5mm}--\ssf invariant case,
\,$E_{\pm\ssf i\vsf\lambda}(x)$ is replaced by \,$J_{\pm\ssf i\vsf\lambda}(x)$
in \eqref{DunklTransform} and \eqref{InverseDunklTransform}.

\item
The Dunkl transform of a radial function is again a radial function.

\item
The following sharper version of the Paley--Wiener Theorem
was proved in \cite{AmriAnkerSifi2010},
as a consequence of the corresponding result in the trigonometric setting
$($\ssb see Theorem \ref{PropertiesCherednikTransform}\hspace{.5mm}$)$
and thus under the assumption that \ssf$R$ is crystallographic.
Given a \ssf$W$\hspace{-.5mm}--\ssf invariant convex compact
neighborhood \,$C$ of the origin in \,$\apartment$\ssf,
consider the gauge \,$\chi(\lambda)\ssb
=\ssb\max_{\ssf x\vsf\in\vsf C}\,\langle\ssf\lambda\ssf,x\ssf\rangle\ssf$.
Then the support of \,$f\!\in\ssb\mathcal{C}_{\ssf c}^{\ssf\infty}(\apartment)$
is contained in \,$C$ if and only if its Dunkl transform \,$h\ssb=\ssb\mathcal{H}f$ satisfies the condition
\vspace{-.5mm}
\begin{equation}\label{GeometricPaleyWienerCondition}
\forall\,N\hspace{-1mm}\in\!\N\ssf,
\;\sup\nolimits_{\ssf\lambda\ssf\in\ssf\apartment_{\ssf\C}}\ssf
(1\!+\!|\lambda|)^N\ssf e^{\ssf-\ssf\chi(\Im\lambda)}\,
|\ssf h(\lambda)\ssf|<+\infty\,.
\end{equation}

\end{itemize}
\end{remark}

\begin{problem}
Extend the latter result to the non--crystallographic case.
\end{problem}

\subsection{Heat kernel}

The heat equation
\begin{equation*}\label{HeatEquationDunkl}\begin{cases}
\,\partial_{\ssf t\ssf}u(x,t)=\Delta_{\ssf x\ssf}u(x,t)\\
\,u(x,0)=f(x)
\end{cases}\end{equation*}
can be solved via the Dunkl transform
(under suitable assumptions).
This way, one obtains
\begin{equation*}\label{SolutionHeatEquationDunkl}
u(x,t)=\int_{\ssf\apartment}\ssb dy\;\delta(y)\,f(y)\,h_{\ssf t\ssf}(x,y)\,,
\end{equation*}
where the heat kernel is given by
\begin{equation*}
h_{\ssf t\ssf}(x,y)
=\vsf\crat^{-2}\int_{\ssf\apartment}\ssb d\lambda\;
\delta(\lambda)\,e^{-\ssf t\,|\lambda|^2}
E_{\ssf i\lambda}(x)\,E_{-i\lambda}(y)
\qquad\forall\;t\ssb>\ssb0\,,
\;\forall\;x,y\ssb\in\ssb\apartment\,.
\end{equation*}
In next proposition, we collect
some properties of the heat kernel established by R\"osler.

\begin{proposition}
${}$
\begin{itemize}

\item
$\ssf h_{\ssf t\ssf}(x,y)$ is an smooth symmetric probability density.
More precisely,
\begin{itemize}

\item[$\circ$]
$\ssf h_{\ssf t\ssf}(x,y)$ is an analytic function in
\ssf$(t,x,y)\!\in\!(0,+\infty)\!\times\!\apartment\!\times\!\apartment$\ssf,

\item[$\circ$]
$\ssf h_{\ssf t\ssf}(x,y)\ssb=\ssb h_{\ssf t\ssf}(y,x)$\ssf,

\item[$\circ$]
$\ssf h_{\ssf t\ssf}(x,y)\ssb>\ssb0$ and
\ssf$\int_{\ssf\apartment}\ssb dy\,\delta(y)\,h_{\ssf t\ssf}(x,y)\ssb=\ssb1$\ssf.

\end{itemize}

\item
{\rm Semigroup property\,:}
\vspace{-1mm}
\begin{equation*}
h_{\ssf s+t\ssf}(x,y)=\int_{\apartment}dz\;\delta(z)\,h_{\ssf s\ssf}(x,z)\,h_{\ssf t\ssf}(z,y)\,.
\end{equation*}

\item
{\rm Expression by means of the Dunkl kernel\,:}
\vspace{-1mm}
\begin{equation}\label{ExpressionHeatKernelDunkl}
h_{\ssf t\ssf}(x,y)=\vsf\crat^{-1}\,(2\ssf t)^{-\frac n2-\ssf\gamma}\,
e^{-\frac{|x|^2}{4\ssf t}-\frac{|y|^2}{4\ssf t}}\,
E_{\smash{\frac{x\vphantom{|}}{\sqrt{2\ssf t\ssf}}}}\ssb
\bigl(\tfrac{y\vphantom{|}}{\sqrt{2\ssf t\ssf}}\bigr)
\qquad\forall\;t\ssb>\ssb0\ssf,
\;\forall\;x\ssf,y\ssb\in\ssb\apartment\,.
\end{equation}

\item
{\rm Upper estimate\,:}
\vspace{-2mm}
\begin{equation*}\label{EstimateHeatKernelDunkl}
h_{\ssf t\ssf}(x,y)
\le\vsf\crat^{-1}\hspace{.5mm}(2\ssf t)^{-\frac n2-\ssf\gamma}\,
\max\nolimits_{\,w\ssf\in\ssf W}
e^{-\frac{|\ssf w\ssf x\ssf-\ssf y\ssf|^2}{4\ssf t}}\,.
\end{equation*}

\end{itemize}
\end{proposition}

\begin{remark}
${}$
\begin{itemize}

\item
In \cite{AnkerBensalemDziubanskiHamda2015},
the following sharp heat kernel estimates were obtained in dimension \ssf$1$
$($and also in the product case\ssf$):$
\vspace{-1mm}
\begin{equation}\label{EstimateHeatKernelDunkl1D}
h_{\ssf t}(x,y)\,\asymp\,\begin{cases}
\;t^{-k-\frac12}\,
e^{-\frac{x^2\ssb+\ssf y^2}{4\ssf t}}
&\text{if \;}|\ssf x\ssf y\ssf|\!\le\ssb t\ssf,\\
\;t^{-\frac12}\,(x\ssf y)^{-k}\,
e^{-\frac{(x-y)^2}{4\ssf t}}
&\text{if \;}x\ssf y\ssb\ge\ssb t\ssf,\\
\;t^{\ssf\frac12}\,(-\ssf x\ssf y)^{-k-1}\,
e^{-\frac{(x+y)^2}{4\ssf t}}
&\text{if \,}-\ssb x\ssf y\ssb\ge\ssb t\ssf.\\
\end{cases}
\end{equation}
Notice the lack of Gaussian behavior when \,$-\,x\ssf y\ssb\ge\ssb t$\ssf,
in particular when \,$y\ssb=\ssb-\ssf x$ tends to infinity
faster then \ssf$\sqrt{\ssf t\ssf}$.

\item
The Dunkl Laplacian is the infinitesimal generator
of a Feller--Markov process on \ssf$\apartment\ssf$,
which has remarkable features
$($Brownian motion with jumps\ssf$)$
and which has drawn a lot of attention in the 2000s.
We refer to \cite{Graczyk2008} and \cite{Demni2009}
for probabilistic aspects of Dunkl theory.

\end{itemize}
\end{remark}

\begin{problem}
Prove in general heat kernel estimates similar to \eqref{EstimateHeatKernelDunkl1D}.
\end{problem}

\subsection{Intertwining operator and (dual) Abel transform}

Consider the Abel transform
\vspace{-.5mm}
\begin{equation*}
\mathcal{A}\ssf=\ssf\mathcal{F}^{-1}\ssb\circ\mathcal{H}\,,
\end{equation*}
which is obtained by composing the Dunkl transform \ssf$\mathcal{H}$
\ssf with the inverse Euclidean Fourier transform \ssf$\mathcal{F}^{-1}$
\ssb on \ssf$\apartment$\ssf,
and the dual Abel transform \,$\mathcal{A}^*$,
which satisfies
\begin{equation*}
\int_{\ssf\apartment}dx\;\delta(x)\,f(x)\,\mathcal{A}^*\ssb g(x)\,
=\int_{\ssf\apartment}dy\;\mathcal{A}f(y)\,g(y)\,.
\end{equation*}

\begin{theorem}
\label{DualAbelTransformDunkl}
${}$
\begin{itemize}

\item
The dual Abel transform \,$\mathcal{A}^*$ \ssb coincides
with the intertwining operator~\,$\mathcal{V}$
defined on polynomials by Dunkl
and extended to smooth functions by Trim\`eche.

\item
{\rm Intertwining property\,:} for every \,$\xi\ssb\in\ssb\apartment$\ssf,
\begin{equation}\label{InterwiningPropertyDunkl}
\mathcal{A}\circ D_\xi\ssb=\partial_\xi\circ\mathcal{A}
\qquad\text{and}\qquad
\mathcal{V}\circ\partial_\xi=D_\xi\ssb\circ\mathcal{V}\,.
\end{equation}

\item
{\rm Symmetries\,:}
\vspace{1mm}

\centerline{\begin{tabular}{clcl}
$\mathcal{A}\ssf(wf\vsf)=w\ssf(\mathcal{A}f\vsf)$
&and
&$\mathcal{V}\ssf(w\ssf g)=w\ssf(\mathcal{V}g)$
&$\forall\;w\!\in\!W$,
$\vphantom{\frac||}$\\
$\mathcal{A}\ssf\bigl[\ssf f(\ssf t\,.\,)\ssf\bigr](y)
=|t|^{-2\ssf\gamma}(\mathcal{A}f\vsf)(\vsf t\ssf y)$
&and
&$\mathcal{V}\ssf\bigl[\ssf g\ssf(\ssf t\,.\,)\ssf\bigr](x)
=(\mathcal{V}g\vsf)(\vsf t\ssf x)$
&$\forall\;t\!\in\!\R^*$.
$\vphantom{\frac||}$\\
\end{tabular}}
\vspace{1mm}

\item
For every \,$x\ssb\in\ssb\apartment$\ssf,
there is a unique Borel probability measure \,$\mu_{\ssf x}$ on \,$\apartment$
such that
\vspace{-.5mm}
\begin{equation}\label{MeasureMuDunkl}
\mathcal{V}g(x)=\int_{\ssf\apartment}\ssb d\mu_{\ssf x}(y)\,g(y)\,.
\end{equation}
\vspace{-3.5mm}

\noindent
The support of \,$\mu_{\ssf x}$ is contained in
the convex hull of \,$W\ssb x$\ssf.
Moreover, if \,$k\!>\!0$\ssf,
the support of \,$\mu_{\ssf x}$ is
\hspace{.5mm}$W$\hspace{-.5mm}--\hspace{.5mm}invariant
and contains \,$W\ssb x$\ssf.

\item
$\ssf\mathcal{A}$ \ssf is an automorphism of the spaces
\,$\mathcal{C}_c^{\ssf\infty}(\apartment)$ and \,$\mathcal{S}(\apartment)$\ssf,
while \ssf$\mathcal{V}$ is an automorphism
of \,$\mathcal{C}^{\ssf\infty}(\apartment)$,
with
\vspace{-.5mm}
\begin{equation*}
|\mathcal{V}g(x)|\le\ssf\max\nolimits_{\ssf y\ssf\in\ssf\co(W\ssb x)}|g(y)|
\qquad\forall\;x\!\in\!\apartment\ssf.
\end{equation*}

\end{itemize}
\end{theorem}

The following integral representations,
which follow from \eqref{BesselFunctionDunkl},
\eqref{InterwiningPropertyDunkl}
and \eqref{MeasureMuDunkl},
generalize \eqref{SphericalFunctionsRn}
and \eqref{BesselFunctionP}
in the present setting.

\begin{corollary}
For every \,$\lambda\ssb\in\ssb\apartment_{\ssf\C}$\ssf,
we have
\vspace{-.5mm}
\begin{equation*}\label{IntegralRepresentationDunklKernel}
E_\lambda(x)=\int_{\ssf\apartment}\ssb
d\mu_{\ssf x}(y)\,e^{\,\langle\ssf\lambda\ssf,\,y\ssf\rangle}
\end{equation*}
\vspace{-3.5mm}

\noindent
and
\vspace{-1mm}
\begin{equation*}\label{IntegralRepresentationBesselFunctionDunkl}
J_\lambda(x)=\int_{\ssf\apartment}\ssb
d\mu_{\ssf x}^{\vsb W}\!(y)\,\Cosh_{\vsf\lambda}(y)\,,
\end{equation*}
\vspace{-2mm}

\noindent
where \,$\Cosh_{\vsf\lambda}$ is defined in \eqref{DefinitionCosh} and
\vspace{-1mm}
\begin{equation*}
\mu_{\ssf x}^{\vsb W}\!
=\tfrac{1\vphantom{\frac00}}{|W|\vphantom{\frac00}}
\sum\nolimits_{\ssf w\vsf\in\vsf W}\mu_{\ssf w\ssf x}\,.
\end{equation*}
\end{corollary}

\begin{remark}
${}$
\begin{itemize}

\item
The first three items in Theorem \ref{DualAbelTransformDunkl}
hold for all multiplicities \,$k\hspace{-.5mm}\in\!K_{\text{\rm reg}}$\ssf.

\item
The following symmetries hold\;$:$
\vspace{-1mm}
\begin{equation*}\begin{cases}
\,d\vsf\mu_{\ssf w\vsf x}(w\ssf y)=d\vsf\mu_{\ssf x}(y)
&\forall\;w\!\in\!W\ssf,\\
\,d\vsf\mu_{\ssf t\vsf x}(t\ssf y)=d\vsf\mu_{\ssf x}(y)
&\forall\;t\!\in\!\R^*.
\end{cases}\end{equation*}
\vspace{-2.5mm}

\item
In \cite{DejeuRoesler2002},
it is conjectured that the measure \,$\mu_{\ssf x}$
is absolutely continuous with respect to the Lebesgue measure
under the following two assumptions$\,:$
\begin{itemize}

\item[$\circ$]
$\ssf x$ is regular
$($which means that
\,$\langle\alpha,x\rangle\ssb\ne\ssb0$\ssf,
for every \,$\alpha\!\in\!R$\ssf$)$\ssf,

\item[$\circ$]
$\ssf\apartment$ \ssf is spanned by the roots \,$\alpha$
with multiplicity \,$k_{\ssf\alpha}\!>\!0$\ssf.
\end{itemize}

\item
These conjectures hold in dimension $1$
$($hence in the product case\ssf$)$,
where
\vspace{-1mm}
\begin{equation*}
d\ssf\mu_{\ssf x}(y)=\tfrac
{\Gamma(k\ssf+\frac12)\vphantom{\frac0|}}
{\sqrt{\pi\ssf}\,\Gamma(k)\vphantom{\frac|0}}\,
(\ssf|x|\!+\ssb\sign(x)\,y\ssf)\,(\ssf x^2\!-\ssb y^2\ssf)^{\ssf k-1}\,
\1_{\ssf(-|x|,\ssf+|x|\ssf)}\vsf(y)\,dy
\end{equation*}
if \,$x\ssb\ne\ssb0$\ssf,
while \,$\mu_{\ssf0}$ is the Dirac measure at the origin.

\end{itemize}
\end{remark}

\subsection{Generalized translations, convolution and product formula}

\begin{definition}
${}$
\begin{itemize}

\item
The generalized convolution corresponds, via the Dunkl transform,
to point\-wise mul\-ti\-pli\-ca\-tion$\,:$
\vspace{-.5mm}
\begin{equation}\label{ConvolutionDunkl}
(\vsf f\hspace{-.4mm}*\ssb g\vsf)\ssf(x)=\vsf\crat^{-2}
\int_{\ssf\apartment}\ssb d\lambda\;\delta(\lambda)\,
\mathcal{H}f(\lambda)\,\mathcal{H}g(\lambda)\,E_{\,i\lambda}(x)\,.
\end{equation}
\vspace{-2mm}

\item
The generalized translations are defined by
\begin{equation}\label{TranslationDunkl}
(\vsf\tau_y\ssf f\ssf)\ssf(x)=\vsf\crat^{-2}
\int_{\ssf\apartment}\ssb d\lambda\;\delta(\lambda)\,
\mathcal{H}f(\lambda)\,E_{\,i\lambda}(x)\,E_{\,i\lambda}(y)
=(\vsf\tau_x\ssf f\ssf)\ssf(y)\,.
\end{equation}

\end{itemize}
\end{definition}

The key objects here are the tempered distributions
\begin{equation}\label{KernelProductDunkl}
f\ssf\longmapsto\,\langle\,\nu_{\ssf x,\ssf y}\,,f\,\rangle\,,
\end{equation}
which are defined by \eqref{TranslationDunkl}
and which enter the product formula
\begin{equation}\label{ProductFormulaDunkl}
E_\lambda(x)\,E_\lambda(y)=\langle\,\nu_{\ssf x,\ssf y}\,,E_\lambda\,\rangle\,.
\end{equation}

\begin{remark}
${}$
\begin{itemize}

\item
When \,$k\ssb=\ssb0$\ssf,
\eqref{ConvolutionDunkl} reduces to the usual convolution on \,$\apartment$\ssf,
\eqref{TranslationDunkl} to \,$(\tau_yf)(x)\ssb=\ssb f(x\ssb+\ssb y)$\ssf,
and \,$\nu_{\ssf x,\ssf y}\ssb=\delta_{\ssf x+y}$\ssf.

\item
In the \,$W$\!--\ssf invariant case,
\eqref{ConvolutionDunkl} becomes
\vspace{-.5mm}
\begin{equation*}
(\vsf f\hspace{-.4mm}*\ssb g\vsf)\ssf(x)=\vsf\crat^{-2}
\int_{\ssf\apartment}\ssb d\lambda\;\delta(\lambda)\,
\mathcal{H}f(\lambda)\,\mathcal{H}g(\lambda)\,J_{\,i\lambda}(x)
\end{equation*}
\vspace{-4.5mm}

\noindent
and \eqref{ProductFormulaDunkl}
\vspace{-.5mm}
\begin{equation}\label{ProductFormulaDunklSymmetric}
J_\lambda(x)\,J_\lambda(y)=\langle\,\nu_{\ssf x,\ssf y}^W\,,J_\lambda\,\rangle\,,
\end{equation}
\vspace{-4.5mm}

\noindent
where
\vspace{-.5mm}
\begin{equation*}\label{KernelProductDunklSymmetric}
\nu_{\ssf x,\ssf y}^{\,W}=\tfrac{1\vphantom{\frac00}}{|W|\vphantom{\frac00}}
\sum\nolimits_{\ssf w\in W}\nu_{\ssf w\ssf x,\ssf w\ssf y\vphantom{\frac00}}\,.
\end{equation*}

\end{itemize}
\end{remark}

\begin{lemma}\label{SupportNu}
The distributions \eqref{KernelProductDunkl} are compactly supported.
\begin{itemize}

\item
Specifically, \,$\nu_{\ssf x,\ssf y}$ is supported in the spherical shell
\begin{equation*}
\bigl\{\,z\!\in\!\apartment\bigm|
\bigl|\ssf|x|\!-\!|y|\ssf\bigr|\ssb\le\ssb|z|\ssb\le\ssb|x|\!+\!|y|\,\bigr\}\,.
\end{equation*}

\item
Assume that \,$W$ is crystallographic.
Then \,$\nu_{\ssf x,\ssf y}$ is actually supported in
\begin{equation}\label{SupportProductDunkl}
\bigl\{\,z\!\in\!\apartment\bigm|
z^+\hspace{-1mm}\preccurlyeq\hspace{-.5mm}
x^+\hspace{-1mm}+\hspace{-.5mm}y^+,\,
z^+\hspace{-.25mm}\succcurlyeq
y^+\hspace{-1mm}+\ssb w_{\ssf0}\ssf x^+
\hspace{1mm}\text{and }\hspace{1mm}
x^+\hspace{-1mm}+\ssb w_{\ssf0}\ssf y^+
\,\bigr\}\,,
\end{equation}
where \,$\preccurlyeq$ denotes the partial order on \,$\apartment$
associated with the cone \,$\overline{\apartment_+\vphantom{0}}$.

\end{itemize}
\end{lemma}

\begin{figure}[b]
\begin{center}
\psfrag{x+y}[c]{\color{red}$x+y$}
\psfrag{x-y}[c]{\color{red}$x-y$}
\includegraphics[height=70mm]{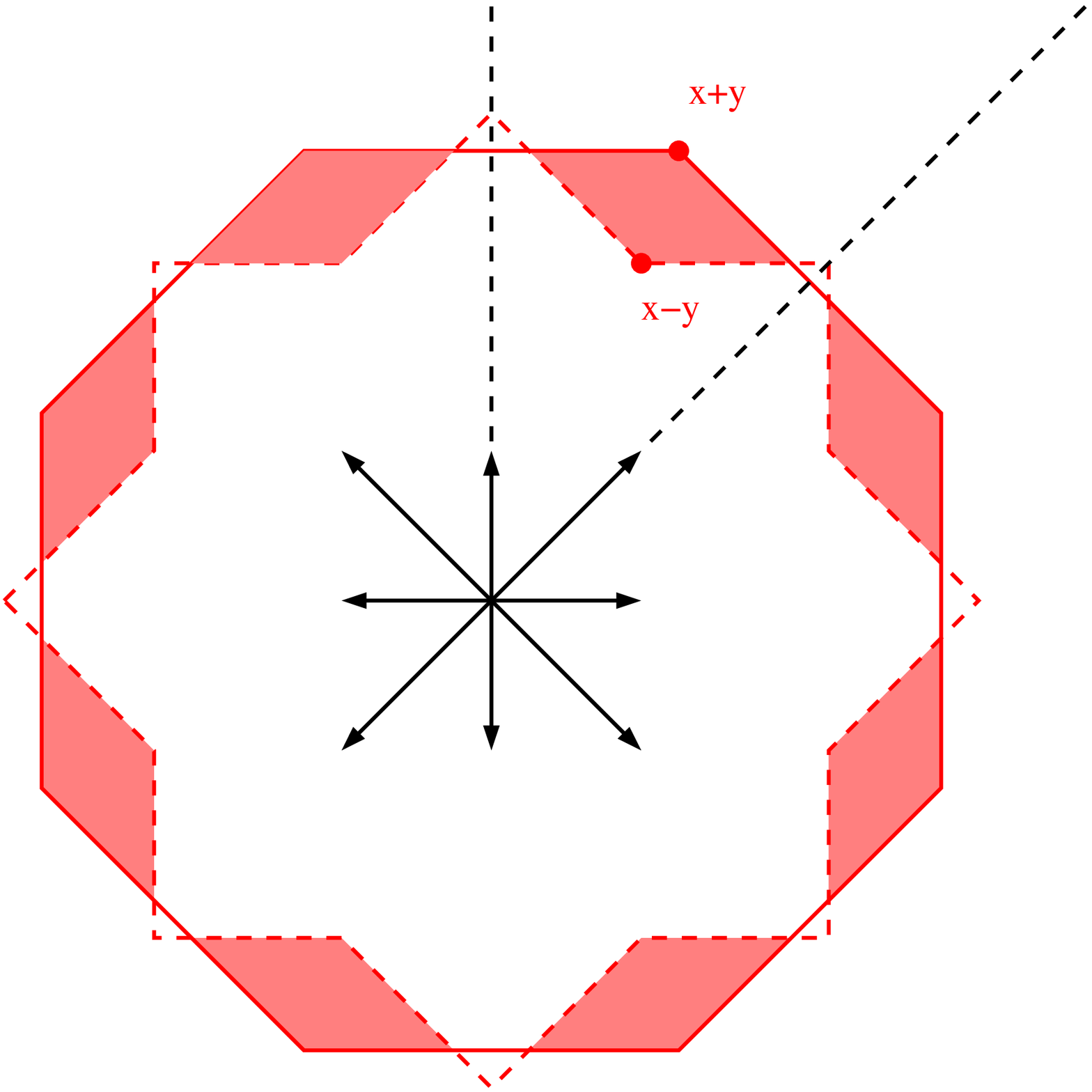}
\caption{Picture of the set \eqref{SupportProductDunkl}
for the root system \ssf$B_{\ssf2}$}
\end{center}
\end{figure}

\begin{example}
In dimension \ssf$1$\ssf,
\,$\nu_{\ssf x,\ssf y}$ is a bounded signed measure.
Specifically,
\vspace{-.5mm}
\begin{equation*}
d\ssf\nu_{\ssf x,\ssf y}\ssf(z)=\begin{cases}
\,\nu(x,y,z)\,|z|^{\ssf2\ssf k}\,dz
&\text{if }\,x,y\ssb\in\ssb\R^*,\\
\,d\ssf\delta_y(z)
&\text{if }\,x\ssb=\ssb0\,,\\
\,d\ssf\delta_x(z)
&\text{if }\,y\ssb=\ssb0\,,\\
\end{cases}\end{equation*}
\vspace{-4mm}

\pagebreak

\noindent
where
\vspace{-.5mm}
\begin{equation*}
\nu(x,y,z)=\tfrac
{\Gamma(k\ssf+\frac12)\vphantom{\frac0|}}
{\sqrt{\pi\ssf}\,\Gamma(k)\vphantom{\frac|0}}\,
\tfrac{\vphantom{\frac0|}
(z\ssf+\ssf x\ssf+\ssf y)\ssf(z\ssf+\ssf x\ssf-\ssf y)\ssf(z\ssf-\ssf x\ssf+\ssf y)}
{2\,x\ssf y\ssf z\vphantom{\frac|0}}\,
\tfrac{\vphantom{\frac0|}
\{\ssf|z|^2-\ssf(\vsf|x|-|y|\vsf)^2\}^{\vsf k\vsf-1}\ssf
\{(\vsf|x|+|y|\vsf)^2-\ssf|z|^2\}^{\vsf k\vsf-1}}
{\{\ssf2\ssf|x|\ssf|y|\ssf|z|\ssf\}^{\ssf2\ssf k\ssf-1}\vphantom{\frac|0}}
\end{equation*}
\vspace{-3.5mm}

\noindent
if \,$x,\ssb y,\ssb z\!\in\!\R^*$ satisfy the triangular inequality
\,$\bigl||x|\!-\!|y|\bigr|\!<\!|z|\!<\!|x|\!+\!|y|$
and \,$\nu(x,y,z)\ssb=\ssb0$ otherwise.
Moreover, let
\vspace{-1mm}
\begin{equation*}
M=\ssf\sup\nolimits_{\ssf x,y\ssf\in\ssf\R}\ssf
\int_{\ssf\R}\ssb d\ssf|\nu_{\ssf x,\ssf y}|(z)
=\sqrt{2\ssf}\,\tfrac
{\left[\ssf\Gamma(k\ssf+\frac12)\ssf\right]^2\vphantom{\frac0|}}
{\Gamma(k\ssf+\frac14)\,\Gamma(k\ssf+\frac34)\vphantom{\frac|0}}\,.
\end{equation*}
\vspace{-3mm}

\noindent
Then \,$M\!\ge\!1$ and
\,$M\hspace{-1mm}\nearrow\!\sqrt{2\ssf}$
\ssf as \,$k\hspace{-.75mm}\nearrow\!+\infty$\ssf.
\end{example}

In general there is a lack information about \eqref{KernelProductDunkl}
and the following facts are conjectured \cite{Roesler2003}.

\begin{problem}
\label{Problem1DunklProductKernel}
${}$

\noindent{\rm(a)}
The distributions \,$\nu_{\ssf x,\ssf y}$ are bounded signed Borel measures.

\noindent{\rm(b)}
They are uniformly bounded in \,$x$ and \,$y$\ssf.

{\rm(c)}
The measures \,$\nu_{\ssf x,\ssf y}^{\,W}$ are positive.
\end{problem}

If \ssf$\nu_{\,x,\ssf y}^{(W\vsb)}$ is a measure, notice that it is normalized by
\begin{equation*}
\int_{\ssf\apartment}\ssb d\ssf\nu_{\,x,\ssf y}^{(W\vsb)}\ssb(z)=1\,.
\end{equation*}
Problem \ref{Problem1DunklProductKernel} and especially item (b)
is important for harmonic analysis.
It implies indeed the following facts,
for the reference measure \ssf$\delta(x)\ssf dx$ on $\apartment$\ssf.

\begin{problem}
\label{Problem2DunklProductKernel}
${}$

\noindent{\rm(d)}
The generalized translations \eqref{TranslationDunkl}
are uniformly bounded on \ssf$L^1$ and hence on \ssf$L^{\ssf p}$,
for every $1\!\le\hspace{-.5mm}p\hspace{-.4mm}\le\!\infty$\ssf.

\noindent{\rm(e)}
{\rm Young's inequality\,:}
For all \hspace{.5mm}$1\!\le\hspace{-.5mm}p,\ssb q,\ssb r\!\le\!\infty$
satisfying \,$\frac1p\!+\!\frac1q\!-\!\frac1r\!=\!1$\ssf,
there exists a con\-stant \,$C\!\ge\hspace{-.5mm}0$ such that
\begin{equation}\label{Young}
\|\ssf f\hspace{-.4mm}*\ssb g\hspace{.4mm}\|_{L^r}
\le C\,\|\ssf f\ssf\|_{L^p}\,\|\hspace{.4mm}g\hspace{.4mm}\|_{L^q}\,.
\end{equation}
\end{problem}

Beside the trivial \ssf$L^2$ setting
and the one--dimensional case (hence the product case),
here are two more situations where
Problems \ref{Problem1DunklProductKernel} and
\ref{Problem2DunklProductKernel} have been solved.
\medskip

\begin{itemize}

\item
\textbf{Radial case \cite{Roesler2003}.}
Translations of radial functions are positive.
Specifically,
for radial functions \,$f(z)\hspace{-.5mm}=\ssb{}^{\prime}\!f(|z|)$
\ssf and nonzero \,$y\!\in\!\apartment$\ssf,
we have
\begin{equation}\label{TranslationDunklRadial}
(\vsf\tau_{\vsf y}\vsf f\vsf)(x)
=\ssb\int_{\ssf\apartment}\ssb d\vsf\mu_{\vsb\frac y{|y|}\ssb}(z)\;
{}^{\prime}\!f
\bigl(\sqrt{\ssf|x|^2\!+\ssb|y|^2\!+\ssb2\,\langle x,z\rangle\,|y|\ssf}\,\bigr)\,.
\end{equation}
Hence \eqref{Young} holds if \ssf$f$ \ssb or \ssf$g$ is radial.
\medskip

\item
\textbf{Symmetric space case.}
Assume that the multiplicity \ssf$k$ \ssf corresponds to
a Riemannian symmetric space of Euclidean type.
Then \,$\nu_{\ssf x,\ssf y}^W$ \ssf is a positive measure
and \eqref{Young} holds for \ssf$W$\ssb--\ssf invariant functions.

\end{itemize}

\subsection{Comments, references and further results}
\label{CommentsRationalDunkl}

\begin{itemize}

\item
The computation of the integral \eqref{MehtaConstant} has a long history.
A closed form was conjectured by Mehta
for the root systems of type \ssf$\text{\rm A}$
\ssf and by Macdonald for general root systems.
For the four infinite families of classical root systems,
it can be actually deduced from an earlier integral formula of Selberg.
In general, the Mehta--Macdonald formula was proved by Opdam,
first for crystallographic root systems \cite{Opdam1989}
and next for all root systems \cite{Opdam1993}.
His proof was simplified by Etingof \cite{Etingof2010},
who removed in particular the computer--assisted calculations used in the last cases.
\smallskip

\item
We have not discussed the shift operators,
which move the multiplicity \ssf$k$ \ssf by integers
and which have proven useful in the $W$\hspace{-.5mm}--\ssf invariant setting
(see \cite{Opdam1993}).
\smallskip

\item
The following asymptotics hold for the Dunkl kernel,
under the assumption that \,$k\ssb\ge\ssb0$ $($see \cite{DejeuRoesler2002}$):$
there exists \,$\upsilon\hspace{-.5mm}:\hspace{-.5mm}W\!\longrightarrow\ssb\C$
such that, for every \,$w\!\in\!W$ \ssb and
for every \,$\lambda\vsf,\ssb x\!\in\!\chamber\ssb,$
\begin{equation}\label{AsymptoticsDunklKernel}
\lim\nolimits_{\,t\ssf\to+\infty}\ssf
(i\ssf t\ssf)^\gamma\,
e^{-\ssf i\ssf t\,\langle\lambda\ssf,\ssf w\ssf x\ssf\rangle}\,
E_{\ssf i\ssf t\vsf\lambda}(w\ssf x)
=\upsilon(w)\,\delta(\lambda)^{-\frac12}\,\delta(x)^{-\frac12}\,\,.
\end{equation}
If \,$w\ssb=\ssb\Id$\ssf, we have
\,$\upsilon(\Id)\ssb=\ssb(2\pi)^{-\frac n2}\ssf c$\ssf,
where \ssf$c$ \ssf is defined by \eqref{MehtaConstant},
and \eqref{AsymptoticsDunklKernel} holds more generally
when \,$\widetilde{t}\ssb=\ssb i\ssf t$ \,tends to infinity
in the half complex space
\,$\{\,\widetilde{t}\!\in\!\C\,|\Re\widetilde{t}\ssb\ge\ssb0\ssf\}$\ssf.
If \,$w\ssb\ne\ssb\Id$\ssf,
the Dunkl kernel is expected to have a different asymptotic behavior,
when \,$\Re\widetilde{t}$ \,becomes positive.
In dimension \ssf$1$\ssf, we have indeed, for any \ssf$0\!<\!\epsilon\!\le\!\frac\pi2$\ssf,
\vspace{-.5mm}
\begin{equation*}
\lim\nolimits_{\hspace{-2.5mm}\substack{
\widetilde{t}\ssf\to+\infty\\
|\ssb\arg\widetilde{t}\,|\ssf\le\frac\pi2-\ssf\epsilon}
\hspace{-3mm}}
\widetilde{t}^{\hspace{.75mm}k\vsf+1}\ssf
e^{-\ssf\widetilde{t}\ssf\lambda\ssf x}\ssf
E_{\ssf\widetilde{t}\ssf\lambda}(-x)
=\tfrac{2^{\ssf k-1}\,k\;\Gamma(\ssf k+\frac12)}{\sqrt{\pi\ssf}}\,
\lambda^{\ssb-k-1}\,x^{-k-1}\,.
\end{equation*}
\vspace{-3mm}

\noindent
This discrepancy plays an important role
in \cite{AnkerBensalemDziubanskiHamda2015}.
\smallskip

\item
In the fourth item of Theorem \ref{DualAbelTransformDunkl},
the sharper results about the support of \ssf$\mu_{\ssf x}$
when \ssf$k\!>\!0$ \ssf
were obtained in \cite{GallardoRejeb2016b}.
\smallskip

\item
Specific information is available
in the \ssf$W$\hspace{-.5mm}--\ssf invariant setting
for the root systems \ssf$\text{\rm A}_{\ssf n}\ssf$.
In this case, an integral recurrence formula over \ssf$n$ \ssf
was obtained in \cite{Amri2014} and \cite{Sawyer2016} for $J_\lambda$\vsf,
by taking rational limits of corresponding formulae in the trigonometric case
(see the tenth item in Subsection 4.7).
Moreover, an explicit expression of \ssf$\mu_{\ssf x}^W$
\ssb is deduced in~\cite{Sawyer2016}.
As a consequence, the support of \ssf$\mu_{\ssf x}$
is shown to be equal to
the convex hull of \hspace{.5mm}$W\ssb x$\ssf, when \ssf$k\!>\!0$\ssf.
\smallskip

\item
The asymmetric setting is harder.
Beyond the one--dimensional case (and the product case),
explicit expressions of the measure \,$\mu_{\ssf x}$
are presently available in some two--dimensional cases.
For the root system \ssf$\text{\rm A}_{\vsf2}$\ssf,
two closely related expressions were obtained,
first in \cite{Dunkl1995} and recently in \cite{Amri2016}.
For the root system \ssf$\text{\rm B}_{\hspace{.15mm}2}\ssf$,
a complicated formula was obtained in \cite{Dunkl2007}
and a simpler one recently in \cite{AmriDemni2016}.
The case of dihedral root systems \,$\text{\rm I}_{\ssf2}(m)$
is currently investigated $($see \cite{DeleavalDemniYoussfi2015},
\cite{ConstalesDebieLian2016} and the references therein\ssf$)$.
\smallskip

\item
In Lemma \ref{SupportNu},
the sharper result in the crystallographic case
was obtained in \cite{AmriAnkerSifi2010}.
\smallskip

\item
An explicit product formula
was obtained in \cite{Roesler2007}
for generalized Bessel functions
associated with root systems of type \ssf$\text{\rm B}$
and for three one--dimensional families of multiplicities
(which are two--dimensional in this case).
The method consists in computing a product formula
in the symmetric space case,
which corresponds to a discrete set of multiplicities~\ssf$k\ssf$,
and in extending it holomorphically in \ssf$k$\ssf.
The resulting measure lives in a matrix cone,
which projects continuously onto \ssf$\smash{\overline{\chamber}}$,
and its image \ssf$|W|\hspace{.5mm}\nu_{\ssf x,\ssf y}^{\,W}$
\ssf is a probability measure
if \ssf$k\hspace{-.5mm}\ge\hspace{-.5mm}0$\ssf.
\smallskip

\item
Potential theory in the rational Dunkl setting has been studied in
\cite{MejjaoliTrimeche2001},
\cite{GallardoGodefroy2003},
\cite{MaslouhiYoussfi2007},
\cite{GallardoRejeb2015},
\cite{GallardoRejeb2016a},
\cite{GallardoRejeb2016b},
\cite{GallardoRejebSifi2016},
\cite{GraczykLuksRoesler2016}
(see also \cite{Rejeb2015} and the references therein).
\smallskip

\item
Many current works deal with generalizations of results in Euclidean harmonic analysis
to the rational Dunkl setting. Among others, let us mention
\begin{itemize}

\item[$\circ$]
\cite{ThangaveluXu2005}
about the Hardy--Littlewood and the Poisson maximal functions,

\item[$\circ$]
\cite{AmriSifi2012a} and \cite{AmriSifi2012b}
about singular integrals and Calderon--Zygmund theory,

\item[$\circ$]
\cite{AnkerBensalemDziubanskiHamda2015}
and \cite{Dziubanski2015}
about the Hardy space \ssf$H^1$.

\end{itemize}
\smallskip

\item[$\bullet$]
A further interesting deformation of Euclidean Fourier analysis,
encompassing rational Dunkl theory and the Laguerre semigroup,
was introduced and studied in \cite{BensaidKobayashiOrsted2012}.

\end{itemize}

\section{Trigonometric Dunkl theory}

Trigonometric Dunkl theory was developed
in the symmetric case by Heckman and Opdam in the 1980s,
and in the non--symmetric case by Opdam and Cherednik in the 1990s.
This theory of special functions in several variables encompasses
\begin{itemize}

\item
Euclidean Fourier analysis
(which corresponds to the multiplicity \ssf$k\ssb=\ssb0$\ssf)\ssf,

\item
Jacobi functions in dimension \ssf$1$\ssf,

\item
spherical functions associated with
Riemannian symmetric spaces of noncompact type
(which correspond to a discrete set of multiplicities \ssf$k$\ssf)\ssf.

\end{itemize}
In this subsection,
we use \cite{Opdam2000} as our primary reference
and quote only later works.
We resume the notation of Section \ref{RationalDunklTheory},
with some modifications\,:
\begin{itemize}

\item
the root system \ssf$R$ \ssf is now assumed to be
crystallographic but not necessarily reduced,

\item
$\widetilde{R}$ \ssf denotes the subsystem of non--multipliable roots,

\item
the reference density in the case \ssf$k\ssb\ge\ssb0$ \ssf is now
\,$\delta(x)\ssb=\ssb\smash{\displaystyle\prod\nolimits_{\ssf\alpha\vsf\in R^+}}
\bigl|\,2\ssf\sinh\frac{\langle\ssf\alpha\ssf,x\ssf\rangle}2\,\bigr|
{\vphantom{|}}^{\ssf2\ssf k_\alpha}$\,,

\end{itemize}
and some addenda\,:

\begin{itemize}

\item
$\ssf Q$ \ssf denotes the root lattice,
$Q^\vee$ the coroot lattice
and $P$ the weight lattice,

\item
$\ssf\rho\ssf
=\ssb{\displaystyle\sum\nolimits_{\ssf\alpha\in\R^+}}
\ssb\tfrac{k_\alpha}2\,\alpha$\,.

\end{itemize}

\subsection{Cherednik operators}

\begin{definition}
The trigonometric Dunkl operators, which are often called Cherednik operators,
are the differential--difference operators defined by
\begin{equation}\label{CherednikOperators}
D_\xi f(x)=\ssf\partial_{\ssf\xi}f(x)
+\ssf\sum\nolimits_{\ssf\alpha\in R^+}\!k_\alpha\,
\tfrac{\langle\ssf\alpha\ssf,\,\xi\ssf\rangle}
{1\ssf-\,e^{-\ssf\langle\alpha\ssb,\ssf x\rangle}}\,
\bigl\{\ssf f(x)\ssb-\ssb f(r_{\ssb\alpha\ssf}x)\ssf\bigr\}\ssf
-\ssf\langle\ssf\rho\ssf,\xi\ssf\rangle\,f(x)
\end{equation}
for every \,$\xi\!\in\!\mathfrak a$\ssf.
\end{definition}

Notice that the counterpart of Remark \ref{RemarksDunklOperators}
holds in the present setting.
In next theorem,
we collect properties of Cherednik operators.
The main one is again commutativity,
which leads to Cherednik operators \ssf$D_p$\ssf,
for every polynomial \ssf$p\hspace{-.5mm}\in\!\mathcal{P}(\apartment)$\ssf,
and to their symmetric parts \ssf$\smash{\widetilde{D}}_p$
on $W$\hspace{-.5mm}--\ssf invariant functions.

\begin{theorem}
${}$
\begin{itemize}

\item
For any fixed multiplicity \,$k$\ssf,
the Cherednik operators \eqref{CherednikOperators} commute pairwise.

\item
The Cherednik operators map the following function spaces into themselves$\,:$
\begin{equation*}
\C\ssf\bigl[e^{\ssf P\ssf}\bigr]\hspace{.4mm},\;
\mathcal{P}(\apartment)\ssf,\;
\mathcal{C}^{\ssf\infty}(\apartment)\ssf,\;
\mathcal{C}_c^{\ssf\infty}(\apartment)\ssf,\;
\mathcal{S}^{\ssf2}(\apartment)\hspace{-.4mm}
=\hspace{-.5mm}\smash{\bigl(\Cosh_{\vsf\rho}\bigr)^{\ssb-1}}
\mathcal{S}(\apartment)\ssf,\;\dots
\end{equation*}
where \,$\C\ssf\bigl[e^{\ssf P\ssf}\bigr]$ denotes
the algebra of polynomials in \,$e^{\ssf\lambda}$
$(\ssf\lambda\hspace{-.4mm}\in\!P\ssf)$ and
\,$\Cosh_{\vsf\rho}$ is defined in \eqref{DefinitionCosh}.

\item
\text{\rm $W$\hspace{-.5mm}--\ssf equivariance\,:}
For every \,$w\!\in\!W$ \hspace{-.5mm}and \,$\xi\!\in\!\apartment$\ssf,
we have
\begin{equation*}
\bigl(\ssf w\circ D_\xi\circ w^{-1}\bigr)\ssf f\ssf(x)
=D_{w\ssf\xi}f(x)+\sum\nolimits_{\ssf\alpha\ssf\in R^+\ssb\cap\ssf wR^-}
k_\alpha\,\langle\ssf\alpha\ssf,w\ssf\xi\ssf\rangle\,f(r_{\ssb\alpha\ssf}x)\,.
\end{equation*}
\vspace{-4mm}

\noindent
Hence \,$\widetilde{D}_{p\ssf q}\hspace{-.5mm}
=\ssb\widetilde{D}_p\ssb\circ\ssb\widetilde{D}_q$
\ssf for all symmetric polynomials
\,$p\ssf,q\ssb\in\ssb\mathcal{P}(\apartment)^W$.

\item
\text{\rm Adjointness\,:}
Assume that \,$k\ssb\ge\ssb0$\ssf.
Then, for every \,$\xi\!\in\!\apartment$\ssf,
\begin{equation*}
\int_{\ssf\apartment}\ssb dx\;\delta(x)\,D_\xi f(x)\hspace{1mm}g(-x)
=\int_{\ssf\apartment}\ssb dx\;\delta(x)\,f(x)\,D_\xi\vsf g(-x)\,.
\end{equation*}

\end{itemize}
\end{theorem}

\pagebreak

\begin{example}
${}$
\begin{itemize}

\item
The Heckman--Opdam Laplacian is given by
\begin{align*}
\Delta f(x)
=\sum\nolimits_{\ssf j=1}^{\,n}\!D_j^{\ssf2}f(x)
&=\!\overbrace{\vphantom{\frac||}
\sum\nolimits_{\ssf j=1}^{\,n}\!\partial_j^{\ssf2}f(x)
+\ssf\sum\nolimits_{\ssf\alpha\in R^+}\!
k_\alpha\ssf\coth\tfrac{\langle\ssf\alpha\ssf,\,x\ssf\rangle}2\;\partial_\alpha f(x)
+|\rho|^2f(x)}^{\textstyle\text{\rm differential part
\ssf$\smash{\widetilde{\Delta}f(x)}$}}\\
&\ssf-\underbrace{\sum\nolimits_{\ssf\alpha\in R^+}\!
k_\alpha\,
\tfrac{|\alpha|^2\vphantom{\frac00}}
{4\ssf\sinh^2\!\frac{\langle\alpha,\ssf x\rangle}2}\,
\bigl\{\ssf f(x)\ssb-\ssb f(r_{\ssb\alpha\ssf}x)\ssf\bigr\}
\vphantom{\frac||}}_{\textstyle\text{\rm difference part}}\,,
\end{align*}
where \,$D_j$\ssf, respectively \,$\partial_j$ denote
the Cherednik operators, respectively the par\-tial de\-riv\-a\-tives
with respect to an orthonormal basis of \,$\apartment$\ssf.

\item
In dimension \ssf$1$,
the Cherednik operator is given by
\begin{align*}
Df(x)
&=\bigl(\tfrac\partial{\partial x}\bigr)f(x)
+\bigl\{\tfrac{k_1\vphantom{|}}{1\ssf-\,e^{-x}}\ssb
+\ssb\tfrac{2\,k_2\vphantom{|}}{1\ssf-\,e^{-2\ssf x}}\bigr\}\,
\bigl\{\ssf f(x)\!-\!f(-x)\ssf\bigr\}
-\rho\,f(x)\\
&=\bigl(\tfrac\partial{\partial x}\bigr)f(x)
+\bigl\{\ssf\tfrac{k_1}2\ssf\coth\tfrac x2\hspace{-.4mm}
+\ssb k_2\vsf\coth x\ssf\bigr\}\,
\bigl\{\ssf f(x)\!-\!f(-x)\ssf\bigr\}
-\rho\,f(-x)
\end{align*}
and the Heckman--Opdam Laplacian by
\begin{align*}
\Delta f(x)
&=\bigl(\tfrac\partial{\partial x}\bigr)^{\ssb2}f(x)
+\bigl\{\ssf k_1\coth\tfrac x2\hspace{-.4mm}
+\ssb2\,k_2\vsf\coth x\ssf\bigr\}\,
\bigl(\tfrac\partial{\partial x}\bigr)f(x)
+\rho^{\ssf2}f(x)\\
&\ssf-\ssf\Bigl\{\ssf
\tfrac{k_1\vphantom{\frac o0}}{4\ssf\sinh^2\!\frac x2\vphantom{\frac|0}}\ssb
+\ssb\tfrac{k_2\vphantom{\frac o0}}{\sinh^2\!x\vphantom{\frac|0}}
\Bigr\}\,\bigl\{\ssf f(x)\!-\!f(-x)\ssf\bigr\}\,,
\end{align*}
where \,$\rho\ssb=\ssb\frac{k_1}2\!+\ssb k_2$\ssf.
\end{itemize}
\end{example}

\subsection{Hypergeometric functions}

\begin{theorem}
Assume that \,$k\ssb\ge\ssb0$\ssf.
Then, for every \,$\lambda\!\in\!\apartment_{\ssf\C}$\ssf,
the system
\begin{equation*}\label{SystemOpdamHypergeometricFunction}\begin{cases}
\,D_\xi\ssf G_\lambda\ssb=\langle\ssf\lambda\ssf,\xi\ssf\rangle\,G_\lambda
\quad\forall\;\xi\!\in\!\apartment\ssf,\\
\,G_\lambda(0)\ssb=\ssb1\ssf.
\end{cases}\end{equation*}
has a unique smooth solution on \,$\apartment$\ssf,
which is called the Opdam hypergeometric function.
\end{theorem}

\begin{definition}
The Heckman--Opdam hypergeometric function is the average
\begin{equation}\label{HeckmanOpdamHypergeometricfunction}
F_\lambda(x)
=\tfrac{1\vphantom{\frac00}}{|W|\vphantom{\frac00}}
\sum\nolimits_{\ssf w\in W}G_\lambda(w\ssf x)\,.
\end{equation}
\end{definition}

\begin{remark}
${}$
\begin{itemize}

\item
Conversely, \,$G_\lambda(x)$
can be recovered by applying to \,$F_\lambda(x)$
 a linear differential operator in \,$x$
 whose coefficients are rational functions of \ssf$\lambda$\ssf.

\item
The expression \,$G_\lambda(x)$ extends to a holomorphic function of
\,$\lambda\ssb\in\ssb\apartment_{\ssf\mathbb{C}}$\ssf,
\ssf$x\ssb\in\ssb\apartment\ssb+\ssb i\,U$
and \,$k\ssb\in\hspace{-.5mm}V$\ssb,
where \,$U$ is a \,$W$\hspace{-.5mm}--\ssf invariant open neighborhood of \,$0$ in \,$\apartment$
and \,$V$ is a \,$W$\hspace{-.5mm}--\ssf invariant open neighborhood of
\,$\{\ssf k\ssb\in\!K\,|\,k\ssb\ge\ssb0\ssf\}$\ssf.

\item
The Heckman--Opdam hypergeometric function
\eqref{HeckmanOpdamHypergeometricfunction}
is characterized by the system
\begin{equation*}\label{SystemHeckmanOpdamHypergeometricFunction}\begin{cases}
\,\widetilde{D}_p\,F_\lambda\ssb=p\vsf(\lambda)\,F_\lambda
\quad\forall\;p\ssb\in\ssb\mathcal{P}(\apartment)^W,\\
\,F_\lambda(0)\ssb=\ssb1\ssf.
\end{cases}\end{equation*}

\item
In dimension \ssf$1$,
the Heckman--Opdam hypergeometric function reduces to
the Gauss hypergeometric function \,${}_2\ssf\text{\rm F}_{\ssb1}$
or, equivalently, to the Jacobi functions \,$\phi_\lambda^{\alpha,\beta}:$
\vspace{-1mm}
\begin{equation*}
F_\lambda(x)={}_2\ssf\text{\rm F}_{\ssb1}(
\rho\hspace{-.4mm}+\!\lambda\ssf,\rho\hspace{-.4mm}-\!\lambda\ssf;
k_1\hspace{-.8mm}+\hspace{-.4mm}k_2\!+\!\tfrac12\ssf;-\sinh^2\ssb\tfrac x2\ssf)
=\phi_{\,i\,2\ssf\lambda}^{\ssf k_1\ssb+\ssf k_2-\frac12,\,k_2-\frac12}(\tfrac x2)\,,
\end{equation*}
and the Opdam hypergeometric function to a combination of two such functions$\,:$
\vspace{-1mm}
\begin{equation*}
G_\lambda(x)
=\phi_{\,i\,2\ssf\lambda}^{\ssf k_1\ssb+\ssf k_2-\frac12,\,k_2-\frac12}(\tfrac x2)
+\tfrac{\rho\,+\ssf\lambda}{2\ssf k_1+\ssf2\ssf k_2\vsf+1}\,(\sinh x)\,
\phi_{\,i\,2\ssf\lambda}^{\ssf k_1\ssb+\ssf k_2+\frac12,\,k_2+\frac12}(\tfrac x2)\,.
\end{equation*}

\item
When \,$k\ssb=\ssb0$\ssf,
\,$G_\lambda(x)$ \ssf reduces to the exponential
\,$e^{\,\langle\ssf\lambda\ssf,\,x\ssf\rangle}$
\ssf and \,$F_\lambda(x)$ \ssf to \,$\Cosh_{\vsf\lambda}(x)$\ssf.

\item
The functions \,$\ssf G_{-\rho}$ and \,$F_{-\rho}$ are equal to \ssf$1$\ssf.

\item
Spherical functions \,$\phi_{\ssf\lambda}^{\ssf G}$
on Riemannian symmetric space \ssf$G/K$ \ssf of noncompact type
are Heckman--Opdam hypergeometric functions.
Specifically, if
\vspace{.5mm}

\centerline{$\begin{cases}
\,\text{$\mathcal{R}$ \ssf is the root system
of \ssf$(\mathfrak{g},\apartment)$\ssf,}\\
\,m_{\vsf\alpha}\!=\ssb\dim\mathfrak{g}_{\vsf\alpha}\ssf,
\end{cases}$}

\noindent
set
\vspace{-.5mm}

\centerline{$\begin{cases}
\,R\ssb=\ssb2\,\mathcal{R}\ssf,\\
\,k_{\ssf2\alpha}\!=\ssb\frac12\ssf m_{\vsf\alpha}\ssf.
\end{cases}$}

\noindent
Then
\vspace{-.5mm}

\centerline{$
\phi_\lambda^{\ssf G}(\exp x)=F_{i\frac\lambda2}\ssb(2\ssf x)\,.
$}

\end{itemize}
\end{remark}

We collect in the next two propositions
asymptotics and estimates of the hypergeometric functions.

\begin{proposition}\label{HarishChandraExpansionTrigonometric}
The following Harish--Chandra type expansions hold\,$:$
\begin{gather*}\label{HarishChandraExpansionF}
F_\lambda(x)=\sum\nolimits_{\ssf w\vsf\in W}
\mathbf{c}\vsf(w\lambda)\,\Phi_{w\lambda}(x)\,,\\
G_\lambda(x)=
\tfrac{1\vphantom{\frac0|}}
{{\textstyle\prod\nolimits_{\ssf\alpha\vsf\in\widetilde{R}^+}^{\vphantom{0}}}
\left(\langle\ssf\lambda\vsf,\ssf\alpha^{\ssb\vee\vsb}\rangle)\hspace{.5mm}
-\ssf\frac12\ssf k_{\alpha/2}\ssf-\hspace{.5mm}k_\alpha\right)}
\,\sum\nolimits_{\ssf w\in W}
\mathbf{c}\vsf(w\lambda)\,\Psi_{w,\ssf\lambda}(x)\,.
\nonumber\end{gather*}
\vspace{-2.5mm}

\noindent
Here
\vspace{-1.5mm}
\begin{equation*}
\mathbf{c}\vsf(\lambda)=\ssf c_{\ssf0}\ssf
\prod\nolimits_{\ssf\alpha\vsf\in R^{\vsf+}}\!\tfrac
{\Gamma\vsf(\ssf\langle\ssf\lambda\vsf,\ssf\alpha^{\ssb\vee\vsb}\rangle\ssf
+\frac12\ssf k_{\alpha/2}\ssf)\vphantom{\frac0|}}
{\Gamma\vsf(\ssf\langle\ssf\lambda\vsf,\ssf\alpha^{\ssb\vee\vsb}\rangle\ssf
+\frac12\ssf k_{\alpha/2}\vsf+\ssf k_\alpha\ssf)\vphantom{\frac|0}}\,,
\end{equation*}
where \,$c_{\ssf0}$ is a positive constant
such that \,$\mathbf{c}(\rho)\ssb=\ssb1$\ssf,
and
\begin{equation*}
\Phi_\lambda(x)
=\sum\nolimits_{\ssf\ell\ssf\in\vsf Q^+}
\Gamma_{\ssb\ell\ssf}(\lambda)\,
e^{\,\langle\ssf\lambda\ssf-\vsf\rho\ssf-\ssf\ell\ssf,\,x\ssf\rangle}\,,
\hspace{2mm}
\Psi_{w,\ssf\lambda}(x)
=\sum\nolimits_{\ssf\ell\ssf\in\vsf Q^+}
\Gamma_{\ssb\ell\ssf}(w,\lambda)\,
e^{\,\langle\ssf w\lambda\ssf-\vsf\rho\ssf-\ssf\ell\ssf,\,x\ssf\rangle}
\end{equation*}
\vspace{-4mm}

\noindent
are converging series,
for generic \,$\lambda\!\in\!\apartment_{\ssf\mathbb{C}}$
and for every \,$x\!\in\!\chamber$.
\end{proposition}

\begin{proposition}\label{EstimatesG}
Assume that \,$k\ssb\ge\ssb0$\ssf.
\begin{itemize}

\item
All functions \,$G_\lambda$
with \ssf$\lambda\!\in\!\apartment$
are strictly positive.
\vspace{1mm}

\item
The ground function \,$G_0$ has the following behavior$\,:$
\begin{equation*}\label{EstimateGzero}
G_0(x)\ssf\asymp\,\Bigl\{\,\prod\nolimits_{\hspace{-1mm}
\substack{\vphantom{o}\\\alpha\ssf\in\widetilde{R}^+\\
\langle\alpha,\ssf x\rangle\ssf\ge\ssf0}}\hspace{-1.5mm}
\bigl(\ssf1\!+\ssb\langle\alpha,x\ssf\rangle\ssf\bigr)\ssf\Bigr\}\;
e^{-\,\langle\ssf\rho\vsf,\ssf x^{\vsb+}\rangle}
\qquad\forall\;x\!\in\!\apartment\ssf.
\end{equation*}
\vspace{-4.5mm}

\noindent
In particular,
\begin{equation*}
G_0(x)\ssf\asymp\,\Bigl\{\,\prod\nolimits_{\ssf\alpha\ssf\in\widetilde{R}^+}\ssb
\bigl(\ssf1\!+\ssb\langle\alpha,x\ssf\rangle\ssf\bigr)\ssf\Bigr\}\;
e^{-\,\langle\ssf\rho\vsf,\ssf x^{\vsb+}\rangle}
\end{equation*}
\vspace{-4.5mm}

\noindent
if \,$x\ssb\in\ssb\overline{\smash{\chamber}\vphantom{X}}$\ssf, while
\begin{equation*}
G_0(x)\asymp\ssf e^{-\,\langle\ssf\rho\vsf,\ssf x^{\vsb+}\rangle}
\end{equation*}
\vspace{-6mm}

\noindent
if \,$x\ssb\in\ssb-\,\overline{\smash{\chamber}\vphantom{X}}$.
\vspace{1mm}

\item
For every \,$\lambda\!\in\!\apartment_{\ssf\C}$\ssf,
\,$\mu\!\in\!\apartment$ and \,$x\!\in\!\apartment$\ssf,
we have
\begin{equation*}
|\ssf G_{\lambda+\mu}(x)\ssf|\le
e^{\,\langle\ssf(\Re\lambda)^{\ssb+}\!,\,x^{\vsb+}\rangle}\,G_\mu(x)\,.
\end{equation*}
In particular, the following estimates hold,
for every \,$\lambda\!\in\!\apartment_{\ssf\C}$
and \,$x\!\in\!\apartment$\ssf,
\begin{equation*}
|\ssf G_\lambda(x)\ssf|
\le G_{\ssf\Re\ssb\lambda}(x)
\le G_{\vsf0}(x)\,e^{\,\langle\ssf(\Re\lambda)^{\ssb+}}\,.
\end{equation*}
\end{itemize}
\end{proposition}

\subsection{Cherednik transform}

From now on, we assume that \ssf$k\ssb\ge\ssb0$\ssf.

\begin{definition}
The Cherednik transform is defined by
\begin{equation}\label{CherednikTransform}
\mathcal{H}f(\lambda)=\int_{\ssf\apartment}\ssb
dx\;\delta\vsf(x)\,f(x)\,G_{\vsf i\vsf\lambda}(-x)\,.
\end{equation}
\end{definition}

In next theorem, we collect the main properties of the Cherednik transform.

\begin{theorem}\label{PropertiesCherednikTransform}
${}$
\begin{itemize}

\item
The Cherednik transform is an isomorphism
between the \ssf$L^2$ Schwartz space
\vspace{.5mm}

\centerline{$
\mathcal{S}^{\ssf2}(\mathfrak{a})
=\bigl(\Cosh_{\vsf\rho}\bigr)^{\ssb-1}\,
\mathcal{S}(\mathfrak{a})
$}\vspace{.75mm}

\noindent
and the Euclidean Schwartz space
\,$\mathcal{S}(\mathfrak{a})$\ssf.

\item
{\rm Paley--Wiener Theorem\,:}
The Cherednik transform is an isomorphism
between \,$\mathcal{C}_{\ssf c}^{\ssf\infty}(\apartment)$
and the Paley--Wiener space \,$\mathcal{PW}(\apartment_{\ssf\C})$\ssf.
More precisely,
let \,$C$ be a \,$W$\hspace{-.5mm}--\ssf invariant convex compact
neighborhood  of the origin in \,$\apartment$
and let \,$\chi(\lambda)\ssb
=\ssb\max_{\ssf x\vsf\in\vsf C}\,\langle\ssf\lambda\ssf,x\ssf\rangle$
\ssf be the associated gauge.
Then the support of
\,$f\!\in\ssb\mathcal{C}_{\ssf c}^{\ssf\infty}(\apartment)$
is contained in \,$C$ if and only if \,$h\ssb=\ssb\mathcal{H}f$
satisfies \eqref{GeometricPaleyWienerCondition}.

\item
{\rm Inversion formula\,:}
\begin{equation}\label{InverseCherednikTransform}
f(x)=\vsf\ctrig^{\ssf-2}\int_{\ssf\apartment}\ssb
d\lambda\,\widetilde{\delta}\vsf(\lambda)\,
\mathcal{H}f(\lambda)\,G_{\vsf i\vsf\lambda}(x)\,,
\end{equation}
where
\begin{equation}\begin{aligned}\label{AsymmetricPlancherelMeasure}
\widetilde{\delta}\vsf(\lambda)
&=\,\tfrac{c_{\ssf0}^{\ssf2}\vphantom{\frac0|}}
{|\ssf\mathbf{c}\vsf(i\vsf\lambda)\vsf|^{\ssf2}\vphantom{\frac|0}}\;
\prod\nolimits_{\ssf\alpha\ssf\in\widetilde{R}^{\vsf+}}\!\tfrac
{-\ssf i\ssf\langle\ssf\lambda\vsf,\ssf\alpha^{\ssb\vee\vsb}\rangle\ssf
+\frac12\ssf k_{\alpha/2}\vsf+\ssf k_\alpha\vphantom{\frac0|}}
{-\ssf i\ssf\langle\ssf\lambda\vsf,\ssf\alpha^{\ssb\vee\vsb}\rangle
\vphantom{\frac|0}}\\
&=\,\prod\nolimits_{\ssf\alpha\ssf\in R^{\vsf+}}\!\tfrac
{\Gamma\vsf(\ssf i\ssf\langle\ssf\lambda\vsf,\ssf\alpha^{\ssb\vee\vsb}\rangle\ssf
+\frac12\ssf k_{\alpha/2}\vsf+\ssf k_\alpha\ssf)\vphantom{\frac0|}}
{\Gamma\vsf(\ssf i\ssf\langle\ssf\lambda\vsf,\ssf\alpha^{\ssb\vee\vsb}\rangle\ssf
+\frac12\ssf k_{\alpha/2}\ssf)\vphantom{\frac|0}}
\,\tfrac
{\Gamma\vsf(-\ssf i\ssf\langle\ssf\lambda\vsf,\ssf\alpha^{\ssb\vee\vsb}\rangle\ssf
+\frac12\ssf k_{\alpha/2}\vsf+\ssf k_\alpha+\vsf1)\vphantom{\frac0|}}
{\Gamma\vsf(-\ssf i\ssf\langle\ssf\lambda\vsf,\ssf\alpha^{\ssb\vee\vsb}\rangle\ssf
+\frac12\ssf k_{\alpha/2}\vsf+\vsf1)\vphantom{\frac|0}}
\end{aligned}\end{equation}
and \,$\ctrig$ is a positive constant\/.
\end{itemize}
\end{theorem}

\begin{remark}
${}$
\begin{itemize}

\item
In the \ssf$W$\hspace{-.5mm}--\ssf invariant case,
the Cherednik transform \eqref{CherednikTransform}
reduces to
\vspace{-1mm}
\begin{equation}\label{SymmetricCherednikTransform}
\mathcal{H}f(\lambda)=\int_{\ssf\apartment}\ssb
dx\;\delta\vsf(x)\,f(x)\,F_{\vsf i\vsf\lambda}(-x)
\end{equation}
\vspace{-4mm}

\noindent
and its inverse \eqref{InverseCherednikTransform}
to
\vspace{-.5mm}
\begin{equation}\label{InverseSymmetricCherednikTransform}
f(x)=\vsf\ctrig^{\ssf-2}
\int_{\ssf\apartment}\ssb d\lambda\,
\tfrac{c_{\ssf0}^{\ssf2}\vphantom{\frac0|}}
{|\ssf\mathbf{c}\vsf(i\vsf\lambda)\vsf|^{\ssf2}\vphantom{\frac|0}}\,
\mathcal{H}f(\lambda)\,F_{\vsf i\vsf\lambda}(x)\,.
\end{equation}
It is an isomorphism between 
\,$\mathcal{S}^{\ssf2}(\mathfrak{a})^W\hspace{-1mm}
=\hspace{-.5mm}\bigl(\Cosh_{\vsf\rho}\bigr)^{\ssb-1}\ssf
\mathcal{S}(\mathfrak{a})^W$
\ssb and \,$\mathcal{S}(\mathfrak{a})^W$,
which extends to an isometric isomorphism, up to a positive constant,
between \,$L^2(\apartment\ssf,\ssb\delta(x)\ssf dx)^W$ and
\,$L^2(\apartment\ssf,\ssb|\vsf\mathbf{c}(i\lambda)|^{-2}\ssf d\lambda)^W$.

\item
Formulae \eqref{SymmetricCherednikTransform}
and \eqref{InverseSymmetricCherednikTransform}
are not symmetric,
as the spherical Fourier transform \eqref{SphericalFourierTransformHn}
and its inverse \eqref{InverseSphericalFourierTransformHn}
on hyperbolic spaces \,$\Hn$,
or their counterparts \eqref{SphericalFourierTransformTq}
and \eqref{InverseSphericalFourierTransformTq}
on homogeneous trees \,$\Tq$\ssf.
The asymmetry is even greater between
\eqref{CherednikTransform} and \eqref{InverseCherednikTransform},
where the density \eqref{AsymmetricPlancherelMeasure} is complex--valued.

\item
There is no straightforward Plancherel identity
for the full Cherednik transform \eqref{CherednikTransform}.
Opdam has defined in \cite{Opdam1995} a vector--valued transform
leading to a Plancherel identity in the non\ssf--\,$W$\hspace{-.5mm}--\ssf invariant case.
\end{itemize}
\end{remark}

\subsection{Rational limit}\label{RationalLimit}
Rational Dunk theory (in the crystallographic case)
is a suitable limit of trigonometric Dunk theory,
as Hankel analysis on \ssf$\Rn$ is a limit of spherical Fourier analysis on \ssf$\Hn$.
More precisely, assume that
the root system \ssf$R$ \ssf is both crystallographic and reduced.
Then,
\begin{itemize}

\item
the Dunkl kernel is the following limit of Opdam hypergeometric functions\,:
\begin{equation*}
E_{\ssf\lambda}(x)
=\ssf\lim\nolimits_{\,\epsilon\ssf\to\ssf0}\ssf
G_{\epsilon^{-1}\lambda}(\epsilon\ssf x)\,,
\end{equation*}

\item
the Dunkl transform \ssf$\mathcal{H}_{\ssf\text{\rm rat}}$
is a limit case of the Cherednik transform
\ssf$\mathcal{H}_{\ssf\text{\rm trig}}$\,:
\begin{equation*}
(\vsf\mathcal{H}_{\ssf\text{\rm rat}}f\vsf)(\lambda)
=\ssf\lim\nolimits_{\,\epsilon\ssf\to\ssf0}\,
\epsilon^{-n-2\vsf\gamma}\,\bigl\{\vsf
\mathcal{H}_{\ssf\text{\rm trig}}\ssf[\ssf f(\epsilon^{-1}\ssf.\,)\ssf]
\vsf\bigr\}\ssf(\epsilon^{-1}\lambda)\,,
\end{equation*}

\item
likewise for the inversion formulae \eqref{InverseDunklTransform}
and \eqref{InverseCherednikTransform}\,:
\begin{equation*}
(\vsf\mathcal{H}_{\ssf\text{\rm rat}}^{\ssf-1}\vsf f\vsf)(x)
=\ssf\const\,\lim\nolimits_{\,\epsilon\ssf\to\ssf0}\,
\epsilon^{\,n+2\vsf\gamma}\,\bigl\{\vsf
\mathcal{H}_{\hspace{.5mm}\text{\rm trig}}^{\ssf-1}\ssf[\ssf f(\epsilon\,.\,)\ssf]
\vsf\bigr\}\ssf(\epsilon\ssf x)\,.
\end{equation*}
\end{itemize}

\subsection{Intertwining operator and (dual) Abel transform}

In the trigonometric setting,
consider again the Abel transform
\vspace{-.5mm}
\begin{equation*}
\mathcal{A}\ssf=\ssf\mathcal{F}^{-1}\ssb\circ\mathcal{H}\,,
\end{equation*}
which is obtained by composing the Cherednik transform \ssf$\mathcal{H}$
\ssf with the inverse Euclidean Fourier transform \ssf$\mathcal{F}^{-1}$
\ssb on \ssf$\apartment$\ssf,
and the dual Abel transform \,$\mathcal{A}^*$,
which satisfies
\begin{equation*}
\int_{\ssf\apartment}dx\;\delta(x)\,f(x)\,\mathcal{A}^*\ssb g(x)\,
=\int_{\ssf\apartment}dy\;\mathcal{A}f(y)\,g(y)\,.
\end{equation*}

\begin{remark}
In \cite{Trimeche2010},
\,$\mathcal{A}^*\!=\hspace{-.5mm}\mathcal{V}$
\,is called the trigonometric Dunkl intertwining operator

and \,$\mathcal{A}\ssb=\hspace{-.5mm}\mathcal{V}^*$ \ssf the dual operator.
\end{remark}

\begin{proposition}
${}$
\begin{itemize}

\item
For every \,$\xi\ssb\in\ssb\apartment$\ssf,
\vspace{1mm}

\centerline{$
\mathcal{A}\circ D_\xi\ssb=\partial_\xi\circ\mathcal{A}
\qquad\text{and}\qquad
\mathcal{V}\circ\partial_\xi=D_\xi\ssb\circ\mathcal{V}\,.
$}\vspace{.5mm}

\item
For every \,$x\ssb\in\ssb\apartment$\ssf,
there is a unique tempered distribution \,$\mu_{\ssf x}$ on \,$\apartment$
such that
\vspace{1mm}

\centerline{$
\mathcal{V}g(x)=\langle\,\mu_{\ssf x}\ssf,g\,\rangle\,.
$}\vspace{1mm}

\noindent
Moreover,
the support of \,$\mu_{\ssf x}$ is contained
in the convex hull of \,$W\ssb x$\ssf.

\end{itemize}
\end{proposition}

\begin{corollary}
For every \,$\lambda\ssb\in\ssb\apartment_{\ssf\C}$\ssf,
we have
\vspace{-.5mm}
\begin{equation}\label{RepresentationHypergeometricFunctions}
G_\lambda(x)=\langle\,\mu_{\ssf x}\ssf,e^{\,\lambda}\,\rangle
\quad\text{and}\quad
F_\lambda(x)=\langle\,\mu_{\ssf x}^{\vsb W},\Cosh_{\vsf\lambda}\,\rangle\,,
\end{equation}
\vspace{-5mm}

\noindent
where \,$\Cosh_{\vsf\lambda}$ is defined in \eqref{DefinitionCosh} and
\,$\displaystyle
\mu_{\ssf x}^{\vsb W}\!
=\tfrac{1\vphantom{\frac00}}{|W|\vphantom{\frac00}}
\sum\nolimits_{\ssf w\in W}\mu_{\ssf w\ssf x}$\,.

\end{corollary}

\begin{remark}
${}$
\begin{itemize}

\item
\,The distribution \,$\mu_{\ssf x}$ is most likely a probability measure,
as in the rational setting.

\item
This is true in dimension $1$
$($hence in the product case\ssf$)$,
where
\begin{equation*}
d\ssf\mu_{\ssf x}(y)=\begin{cases}
\hspace{3mm}d\vsf\delta_{\vsf x}(y)
&\text{if \,$x\hspace{-.5mm}=\hspace{-.5mm}0$ or if
\,$k_1\hspace{-1mm}=\hspace{-.5mm}k_2\!=\hspace{-.5mm}0$\ssf.}\\
\,\mu(x,y)\hspace{.5mm}dy
&\text{otherwise}.
\end{cases}
\end{equation*}
As far as it is concerned,
the density \,$\mu(x,y)$ vanishes unless \,$|y|\!<\!|x|$\ssf.
In the generic case,
where \,$k_1\hspace{-.8mm}>\hspace{-.5mm}0$ and \,$k_2\!>\hspace{-.5mm}0$\ssf,
it is given explicitly by
\begin{equation}\label{MuTrigonometricGeneric1D}\begin{aligned}
\mu(x,y)
&=\ssf2^{\ssf k_1\vsb+\ssf k_2-\vsf2}\hspace{.5mm}
\smash{\tfrac{\Gamma(k_1\vsb+\ssf k_2+\frac12)\vphantom{\frac0|}}
{\sqrt{\pi\ssf}\,\Gamma(k_1)\,\Gamma(k_2)\vphantom{\frac|0}}}\,
|\ssb\sinh\tfrac x2\ssf|^{-2\ssf k_1}\hspace{.5mm}
|\ssb\sinh x\ssf|^{-2\ssf k_2}\hspace{.5mm}
\vphantom{\Big|}\\
&\ssf\times\smash{\int_{\ssf|y|}^{\,|x|}}dz\,(\vsf\sinh\tfrac z2)\,
(\cosh\tfrac z2\hspace{-.5mm}-\ssb\cosh\tfrac y2\vsf)^{k_1\vsb-1}\hspace{.5mm}
(\cosh x\ssb-\ssb\cosh z)^{k_2-1}
\vphantom{\Big|}\\
&\hspace{12.4mm}\times(\vsf\sign x)\,
\bigl\{\ssf e^{\ssf\frac x2}\hspace{.5mm}(\vsf2\cosh\tfrac x2\vsf)\ssb
-e^{-\frac y2}\hspace{.5mm}(\vsf2\cosh\tfrac z2\vsf)\bigr\}\,.
\vphantom{\Big|}
\end{aligned}\end{equation}
In the limit case, 
where \,$k_1\hspace{-.8mm}=\hspace{-.5mm}0$ and \,$k_2\!>\hspace{-.5mm}0$\ssf,
\begin{equation}\label{DensityMuLimitCase1D}
\mu(x,y)=2^{\ssf k_2-1}\hspace{.5mm}
\tfrac{\Gamma(k_2+\frac12)\vphantom{\frac0|}}
{\sqrt{\pi\ssf}\,\Gamma(k_2)\vphantom{\frac|0}}\,
|\ssb\sinh x\ssf|^{-2\ssf k_2}\hspace{.5mm}
(\cosh x\ssb-\ssb\cosh y)^{k_2-1}\hspace{.5mm}
(\vsf\sign x)\hspace{.5mm}(e^{\ssf x}\!-\ssb e^{-\ssf y}\ssf)\hspace{.5mm}.
\end{equation}
In the other limit case,
where \,$k_1\hspace{-.8mm}>\hspace{-.5mm}0$ and \,$k_2\!=\hspace{-.5mm}0$\ssf,
the density is half of \eqref{DensityMuLimitCase1D},
with \,$k_2$\ssf, $x$\vsf, $y$ replaced respectively
by \,$k_1$\vsf, $\frac x2$\vsf, $\frac y2$\ssf.
\end{itemize}
\end{remark}

\subsection{Generalized translations, convolution and product formula}

\begin{definition}
${}$
\begin{itemize}

\item
The generalized convolution corresponds, via the Cherednik transform,
to pointwise multiplication$\,:$
\vspace{-.5mm}
\begin{equation*}\label{ConvolutionCherednik}
(\vsf f\hspace{-.4mm}*\ssb g\vsf)\ssf(x)
=\vsf\ctrig^{\ssf-2}\int_{\ssf\mathfrak{a}}\ssb
d\lambda\;\widetilde{\delta}\vsf(\lambda)\,
\mathcal{H}f(\lambda)\,\mathcal{H}g(\lambda)\,G_{\ssf i\lambda}(x)\,.
\end{equation*}
\vspace{-2mm}

\item
The generalized translations are defined by
\begin{equation}\label{TranslationCherednik}
(\vsf\tau_y\ssf f\ssf)\ssf(x)
=\vsf\ctrig^{\ssf-2}\int_{\ssf\mathfrak{a}}\ssb
d\lambda\;\widetilde{\delta}\vsf(\lambda)\,
\mathcal{H}f(\lambda)\,G_{\ssf i\lambda}(x)\,G_{\ssf i\lambda}(y)\,.
\end{equation}

\end{itemize}
\end{definition}

The key objects are again the tempered distributions
\begin{equation*}\label{KernelProductOpdam}
f\ssf\longmapsto\,\langle\,\nu_{\ssf x,\ssf y}\,,f\,\rangle
\end{equation*}
defined by \eqref{TranslationCherednik}
and their averages
\begin{equation*}\label{SymmetricKernelProductOpdam}
\nu_{\ssf x,\ssf y}^{\,W}=\tfrac{1\vphantom{\frac00}}{|W|\vphantom{\frac00}}
\sum\nolimits_{\ssf w\in W}\nu_{w\ssf x,\ssf w\ssf y\vphantom{\frac00}}\,,
\end{equation*}
which enter the product formulae
\begin{equation*}\label{ProductFormulaOpdam}
G_\lambda(x)\,G_\lambda(y)
=\langle\,\nu_{\ssf x,\ssf y}\ssf,G_\lambda\,\rangle
\end{equation*}
\vspace{-6mm}

\noindent
and
\vspace{-1mm}
\begin{equation}\label{ProductFormulaHeckmanOpdam}
F_\lambda(x)\,F_\lambda(y)
=\langle\,\nu_{\ssf x,\ssf y}^{\,W}\ssf,F_\lambda\,\rangle\,.
\end{equation}

\begin{example}\label{ProductOpdam1D}
In dimension \ssf$1$\ssf, 
the distributions \,$\nu_{\ssf x,\ssf y}$ are signed measures,
which are uniformly bounded in \,$x$ and \,$y$\ssf.
Explicitly \cite{AnkerAyadiSifi2012},
\vspace{1mm}

\centerline{$
d\ssf\nu_{\ssf x,\ssf y}\ssf(z)=\begin{cases}
\,\nu(x,y,z)\,dz
&\text{if }\,x,y\ssb\in\ssb\R^*,\\
\,d\ssf\delta_y(z)
&\text{if }\,x\ssb=\ssb0\ssf,\\
\,d\ssf\delta_x(z)
&\text{if }\,y\ssb=\ssb0\ssf,\\
\end{cases}$}\vspace{1mm}

where the density \,$\nu(x,y,z)$ is given by the following formulae,
when \,$x,y,\ssb z\hspace{-.5mm}\in\hspace{-.5mm}\R^*$
satisfy the triangular inequality
\vspace{0mm}

\centerline{$
\bigl|\ssf|x|\ssb-\ssb|y|\ssf\bigr|\ssb<\ssb|z|\ssb<\ssb|x|\ssb+\ssb|y|\,,
$}\vspace{0mm}

and vanishes otherwise.

\begin{itemize}

\item
Assume that \,$k_1\!>\ssb0$ and \,$k_2\!>\ssb0$\ssf. Then
\vspace{1.5mm}
\begin{align*}
\nu(x,y,z)
&=\smash{2^{\ssf k_1\vsb-\vsf2}\hspace{.5mm}
\tfrac{\Gamma(k_1+\ssf k_2\vsf+\frac12)\vphantom{\frac0|}}
{\sqrt{\pi\ssf}\,\Gamma(k_1)\,\Gamma(k_2)\vphantom{\frac|0}}\,
\sign(x\ssf y\ssf z)\,
|\sinh\ssb\tfrac x2\sinh\ssb\tfrac y2\hspace{.4mm}|^{-2\ssf k_1-2\ssf k_2}
\ssf(\cosh\ssb\tfrac z2)^{\vsf2\ssf k_2}}\vphantom{\Big|}\\
&\ssf\times\smash{\int_{\,0}^{\ssf\pi\vphantom{|}}\hspace{-1mm}d\chi\,
(\sin\ssb\chi)^{2\ssf k_2-1}}\vphantom{\Big|}\\
&\hspace{5mm}\times\smash{
\bigl[\,\cosh\ssb\tfrac x2\cosh\ssb\tfrac y2\cosh\ssb\tfrac z2\cos\ssb\chi
-\ssb\tfrac{1\ssf+\hspace{.5mm}\cosh x\hspace{.5mm}
+\hspace{.5mm}\cosh y\hspace{.5mm}+\hspace{.5mm}\cosh z\vphantom{\frac00}}
{4\vphantom{\frac00}}\hspace{1mm}\bigr]_+^{\ssf k_1-1}}\vphantom{\Big|}\\
&\hspace{5mm}\times\smash{\bigl[\,
\sinh\tfrac{\vphantom{\frac00}x\hspace{.5mm}+\hspace{.5mm}y\hspace{.5mm}+\hspace{.5mm}z}{2\vphantom{\frac00}}
\ssb-\vsb2\ssf\cosh\ssb\tfrac x2\cosh\ssb\tfrac y2\sinh\ssb\tfrac z2}
\vphantom{\Big|}\\
&\hspace{11mm}+\ssb\smash{
\tfrac{k_1\vsb+\ssf 2\ssf k_2\vphantom{\frac00}}{k_2\vphantom{\frac00}}
}\ssf\cosh\ssb\tfrac x2\cosh\ssb\tfrac y2\ssf\cosh\ssb\tfrac z2
\hspace{.5mm}(\sin\ssb\chi)^2\vphantom{\Big|}\\
&\hspace{11mm}+\smash{\tfrac{\vphantom{\frac00}
\sinh z\hspace{.5mm}-\hspace{.5mm}\sinh x\hspace{.5mm}-\hspace{.5mm}\sinh y}
{2\vphantom{\frac00}}\hspace{.5mm}\cos\ssb\chi\,\bigr]}\,.\vphantom{\Big|}
\end{align*}

\item
Assume that  \,$k_1\!=\ssb0$ and \,$k_2\!>\ssb0$\ssf. Then
\begin{align*}
\nu(x,y,z)
&=\smash{2^{\ssf2\ssf k_2\vsf-1}\hspace{.5mm}
\tfrac{\Gamma(k_2\vsf+\frac12)\vphantom{\frac0|}}
{\sqrt{\pi\ssf}\,\Gamma(k_2)\vphantom{\frac|0}}\,
\sign(x\ssf y\ssf z)\,
|\vsf(\sinh x)\vsf(\sinh y)\vsf|^{\vsf-\vsf2\ssf k_2}}\vphantom{\Big|}\\
&\ssf\times\smash{\bigl[\hspace{.5mm}
\sinh\tfrac
{x\hspace{.5mm}+\hspace{.5mm}y\hspace{.5mm}+\hspace{.5mm}z\vphantom{\frac00}}
{2\vphantom{\frac00}}
\,\sinh\tfrac
{-\hspace{.5mm}x\hspace{.5mm}+\hspace{.5mm}y\hspace{.5mm}
+\hspace{.5mm}z\vphantom{\frac00}}
{2\vphantom{\frac00}}
\,\sinh\tfrac
{x\hspace{.5mm}-\hspace{.5mm}y\hspace{.5mm}+\hspace{.5mm}z\vphantom{\frac00}}
{2\vphantom{\frac00}}
\,\sinh\tfrac
{x\hspace{.5mm}+\hspace{.5mm}y\hspace{.5mm}-\hspace{.5mm}z\vphantom{\frac00}}
{2\vphantom{\frac00}}
\,\bigr]{\vphantom{|}}^{\ssf k_2}}\vphantom{\Big|}\\
&\ssf\times\smash{\bigl[\hspace{.5mm}\sinh\tfrac
{x\hspace{.5mm}+\hspace{.5mm}y\hspace{.5mm}-\hspace{.5mm}z\vphantom{\frac00}}
{2\vphantom{\frac00}}\,\bigr]{\vphantom{|}}^{-1}
\hspace{1mm}e^{\ssf\frac{x\ssf+\ssf y\ssf-\ssf z}2}}\,.\vphantom{\Big|}
\end{align*}

\item
In the other limit case,
where \,$k_1\!>\ssb0$ and \,$k_2\vspace{-.5mm}=\ssb0\vphantom{\frac00}$\ssf,
the density is again half of the previous one,
with \,$k_2$\ssf, $x$\vsf, $y$ replaced respectively
by \,$k_1$\vsf, $\frac x2$\vsf, $\frac y2$\ssf.

\end{itemize}
\end{example}

In higher dimension, we have the trigonometric counterparts of Problems
\ref{Problem1DunklProductKernel} \& \ref{Problem2DunklProductKernel}
but fewer results than in the rational case. In particular,
there is no formula like \eqref{TranslationDunklRadial} for radial functions.
A new property is the following Kunze--Stein phenomenon,
which is typical of the semisimple setting
and which was proved in \cite{AnkerAyadiSifi2012}
(see also \cite{Ayadi2011})
and \cite{Trimeche2012}.

\begin{proposition}\label{KunzeStein}
Let \,$1\hspace{-.5mm}\le\ssb p\hspace{-.4mm}<\ssb2$\ssf.
Then there exists a constant \,$C\hspace{-.5mm}>\hspace{-.5mm}0$ \ssf such that
\vspace{1.5mm}

\centerline{$
\|\ssf f\hspace{-.4mm}*\ssb g\ssf\|_{L^2}
\le C\,\|\ssf f\ssf\|_{L^p}\ssf\|\ssf g\ssf\|_{L^2}\,,
$}\vspace{2mm}

for every \,$f\!\in\!L^{\vsf p}(\apartment\ssf,\ssb\delta(x)\ssf dx)$
and \,$g\hspace{-.5mm}\in\!L^{\vsf2}(\apartment\ssf,\ssb\delta(x)\ssf dx)$\ssf.
\end{proposition}

\subsection{Comments, references and further results}
\label{CommentsTrigonometricDunkl}

\begin{itemize}

\item
The joint action of the Cherednik operators \ssf$D_p$\ssf,
with \ssf$p\hspace{-.5mm}\in\!\mathcal{P}(\apartment)$\ssf,
and of the Weyl group \ssf$W$ \ssb may look intricate.
It corresponds actually to a faithful representation of
a graded affine Hecke algebra \cite{Opdam1995}.
\smallskip

\item
Heckman \cite{Heckman1991} considered initially
the following trigonometric version
\begin{equation*}\label{HeckmanOperators}
{}^{\prime}\hspace{-.5mm}D_\xi f(x)=\ssf\partial_{\ssf\xi}f(x)
+\ssf\sum\nolimits_{\ssf\alpha\in R^+}\hspace{-1mm}
\tfrac{k_\alpha}2\,\langle\ssf\alpha\ssf,\xi\ssf\rangle\,
\coth\tfrac{\langle\ssf\alpha\ssf,\,x\ssf\rangle}2\,
\bigl\{\ssf f(x)\ssb-\ssb f(r_{\ssb\alpha\ssf}x)\ssf\bigr\}
\end{equation*}
of rational Dunkl operators,
which are closely connected to \eqref{CherednikOperators}$\,:$
\begin{equation*}
D_\xi f(x)={}^{\prime}\hspace{-.5mm}D_\xi f(x)
-\ssf\sum\nolimits_{\alpha\in R^+}\hskip-1mm
\tfrac{k_\alpha}2\,\langle\ssf\alpha\ssf,\xi\ssf\rangle\,
f(r_{\ssb\alpha\ssf}x)\,.
\end{equation*}
These operators are \,$W$\hspace{-.5mm}--\ssf equivariant$\,:$
\vspace{-1mm}
\begin{equation*}
w\circ{}^{\prime}\hspace{-.5mm}D_\xi\circ w^{-1}\ssb
={}^{\prime}\hspace{-.5mm}D_{w\ssf\xi}\,,
\end{equation*}
\vspace{-6.5mm}

\noindent
and skew--invariant$\,:$
\vspace{-.5mm}
\begin{equation*}
\int_{\ssf\apartment}\ssb dx\;\delta(x)\,
(\vsf{}^{\prime}\hspace{-.5mm}D_\xi f\ssf)(x)\hspace{1mm}g(x)
=\ssf-\int_{\ssf\apartment}\ssb dx\;\delta(x)\,
f(x)\,(\vsf{}^{\prime}\hspace{-.5mm}D_\xi\vsf g)(x)\,,
\end{equation*}
\vspace{-4.5mm}

\noindent
but they don't commute\,:
\vspace{-1mm}
\begin{equation*}
[\ssf{}^{\prime}\hspace{-.5mm}D_\xi\ssf,
\ssb{}^{\prime}\hspace{-.5mm}D_\eta\ssf]\,f(x)
=\ssf\sum\nolimits_{\ssf\alpha,\beta\in R^+}
\hspace{-1mm}\tfrac{k_\alpha\ssf k_\beta}4\,
\bigl\{\ssf\langle\ssf\alpha\ssf,\xi\ssf\rangle\,\langle\ssf\beta\ssf,\eta\ssf\rangle
-\langle\ssf\beta\ssf,\xi\ssf\rangle\,\langle\ssf\alpha\ssf,\eta\ssf\rangle\ssf\bigr\}
\,f(r_{\ssb\alpha}\ssf r_{\ssb\beta}\hspace{.4mm}x)\,.
\end{equation*}

\item
The hypergeometric functions \ssf$x\ssb\longmapsto\ssb G_\lambda(x)$
and \ssf$x\ssb\longmapsto\hspace{-.5mm}F_\lambda(x)$
extend holomophically to a tube \hspace{.5mm}$\apartment\ssb+\ssb i\,U$
\ssf in \ssf$\apartment_{\ssf\C}$\ssf.
The optimal width for \ssf$F_\lambda$ was investigated in \cite{KroetzOpdam2008}.
\smallskip

\item
Proposition \ref{HarishChandraExpansionTrigonometric}
was obtained in \cite{Opdam1995}.
The asymptotic behavior of \,$F_\lambda$
was fully determined in \cite{NaranayanPasqualePusti2014}.
This paper contains in particular a proof of the estimate
\vspace{1mm}\begin{equation*}
F_\lambda(x)\ssf\asymp\,\Bigl\{\,\prod\nolimits_{\hspace{-1mm}
\substack{\vphantom{o}\\\alpha\ssf\in\widetilde{R}^+\\
\langle\alpha,\ssf\lambda\rangle\ssf\ne\ssf0}}\hspace{-1.5mm}
\bigl(\ssf1\!+\ssb\langle\alpha,x\ssf\rangle\ssf\bigr)\ssf\Bigr\}\;
e^{\,\langle\ssf\lambda\ssf-\vsf\rho\vsf,\ssf x\ssf\rangle}
\qquad\forall\;\lambda\ssf,x\hspace{-.5mm}\in\hspace{-.5mm}\overline{\smash{\chamber}\vphantom{X}}\ssf,
\end{equation*}
\vspace{-2.5mm}

\noindent
which was stated in \cite{Schapira2008} (see also \cite{Schapira2006}),
and the following generalization of a celebrated result of Helgason \& Johnson
in the symmetric space case\,:
\vspace{1mm}

\centerline{
$F_\lambda$ is bounded if and only if
\ssf$\lambda$ \ssf belongs to the convex hull of \ssf$W\hspace{-.5mm}\rho$\ssf.
}\smallskip

\item
The sharp estimates in Proposition \ref{EstimatesG}
were obtained in \cite{Schapira2008} (see also \cite{Schapira2006})
and \cite{RoeslerKoornwinderVoit2013}.
\smallskip

\item
As in the rational case, we have not discussed the shift operators,
which move the multiplicity \ssf$k$ \ssf by integers
and which have proven useful in the $W$\hspace{-.5mm}--\ssf invariant setting
(see \cite[Part\;I, Ch.\;3]{HeckmanSchlichtkrull1994},
\cite[Section\;5]{Opdam2000}).
\smallskip

\item
Rational limits in Subsection \ref{RationalLimit} have a long prehistory.
In the Dunkl setting, they have been used for instance in
\cite{RoeslerVoit2004},
\cite{BensaidOrsted2005},
\cite{Dejeu2006},
\cite{AmriAnkerSifi2010},
\cite{Amri2014},
\cite{Sawyer2016},
\dots
\,(seemingly first and independently in preprint versions
of \cite{BensaidOrsted2005} and \cite{Dejeu2006}).
\smallskip

\item
There are other interesting limits between
special functions occurring in Dunkl theory.
For instance,
in \cite{RoeslerKoornwinderVoit2013}
and \cite{RoeslerVoit2016},
Heckman--Opdam hypergeometric functions
associated with the root system \ssf$\text{\rm A}_{\ssf n-1}$
are obtained as limits of Heckman--Opdam hypergeometric functions
associated with the root system \ssf$\text{\rm BC}_{\ssf n}$\ssf,
when some multiplicities tend to infinity.
See \cite{RoeslerVoit2008}
for a similar result about generalized Bessel functions.
\smallskip

\item
The expressions \eqref{RepresentationHypergeometricFunctions}
are substitutes for the integral representations
\eqref{SphericalFunctionsHn} and \eqref{SphericalFunctionGK}.
A different integral representation of \ssf$F_\lambda$
is established in \cite{Sun2016}.
\smallskip

\item
Formula \eqref{MuTrigonometricGeneric1D} was obtained in \cite{Anker2016}
and used there to prove the positivity of \ssf$\mu_{\ssf x}$
when \ssf$k_1\hspace{-1mm}>\hspace{-.5mm}0$
\ssf and \ssf $k_2\!>\hspace{-.5mm}0$\ssf.
A more complicated expression was obtained previously
in \cite{GallardoTrimeche2010} and in \cite{Ayadi2011}.
It was used in \cite{GallardoTrimeche2010}
to disprove mistakenly the positivity of \ssf$\mu_{\ssf x}$\vsf.
Later on, the positivity of \ssf$\mu_{\ssf x}$ in the general case
was investigated in \cite{Trimeche2015}, using the positivity of the heat kernel.
\smallskip

\item
It is natural to look for recurrence formulae over \ssf$n$
\ssf for the five families of classical crystallographic root systems
\ssf$\text{\rm A}_{\ssf n}$\ssf,
$\text{\rm B}_{\ssf n}$\ssf,
$\text{\rm C}_{\ssf n}$\ssf,
$\text{\rm BC}_{\ssf n}$\ssf,
$\text{\rm D}_{\ssf n}$
(see the appendix).
In the case of \ssf$\text{\rm A}_{\ssf n}$\ssf,
an integral recurrence formula for \ssf$F_\lambda$
(or for Jack polynomials)
was discovered independently by several authors
(see for instance
\cite{Sawyer1997},
\cite{OkounkovOlshanski1997},
\cite{HallnasRuijsenaars2015}).
An explicit expression of \ssf$\mu_{\ssf x}^W$ \ssb
is deduced in \cite{Sawyer1997} and \cite{Sawyer2016},
first for \ssf$x\!\in\!\chamber$
and next for any \ssf$x\!\in\!\overline{\chamber}$.
In particular,
if \ssf$k\hspace{-.5mm}>\hspace{-.5mm}0$\ssf,
then \ssf$\mu_{\ssf x}^W$ \ssb is a probability measure,
whose support is equal to the convex hull of \ssf$W\ssb x$
and which is absolutely continuous with respect to the Lebesgue measure,
except for \ssf$x\hspace{-.5mm}=\hspace{-.5mm}0$
\ssf where \ssf$\mu_{\ssf x}^W\hspace{-1.25mm}
=\hspace{-.5mm}\delta_{\vsf0}$\ssf.
\smallskip

\item
As in the rational case
(see the eighth item in Subsection \ref{CommentsRationalDunkl}),
an explicit product formula was obtained
in \cite{Roesler2010} and \cite{Voit2015}
for Heckman--Opdam hypergeometric functions
associated with root systems of type \ssf$\text{\rm BC}$
\ssf and for certain continuous families of multiplicities.
\smallskip

\item
Probabilistic aspects of trigonometric Dunkl theory
were studied in \cite{Schapira2008} and \cite{Schapira2007}
(see also \cite{Schapira2006}).
Regarding the heat kernel \ssf$h_{\ssf t}(x,y)$\vsf,
the estimate \eqref{HeatKernelEstimateGK}
was shown to hold for \ssf$h_{\ssf t}(x,0)$
and some asymptotics were obtained for \ssf$h_{\ssf t}(x,y)$\vsf.
But there is no trigonometric counterpart
of the expression \eqref{ExpressionHeatKernelDunkl},
neither precise information like \eqref{EstimateHeatKernelDunkl1D}
about the full behavior of \ssf$h_{\ssf t}(x,y)$\vsf.
\smallskip

\item
The bounded harmonic functions for the Heckman--Opdam Laplacian
were determined in \cite{Schapira2009b}.
\end{itemize}

\appendix
\section{Root systems}
\label{Appendix}

In this appendix, we collect some information about root systems and reflection groups.
More details can be found in classical textbooks
such as \cite{Humphreys1990} or \cite{Kane2001}.

\begin{definition}
Let \,$\apartment\ssb\approx\ssb\R^n$ be a Euclidean space.
\begin{itemize}
\item
A $($\ssb crystallographic\ssf$)$ root system in \hspace{.4mm}$\apartment$ is
a finite set \ssf$R$ of nonzero vectors
satisfying the following conditions\,$:$
\vspace{.5mm}

\noindent
{\rm(a)}
for every \ssf$\alpha\!\in\!R$\ssf, the reflection \,$r_\alpha(x)\ssb
=\ssb x-2\,\tfrac{\langle\ssf\alpha\ssf,\,x\ssf\rangle}{|\alpha|^2}\,\alpha$
\ssf maps \ssf$R$ \ssf onto itself,

\noindent
{\rm(b)} \,$2\,
\smash{\tfrac{\langle\ssf\alpha\ssf,\,\beta\ssf\rangle}{|\alpha|^2}}
\ssb\in\Z$
\,for all \,$\alpha,\beta\!\in\!R$\ssf.
\item
A root system \ssf$R$ is reducible
if it can be splitted into two orthogonal root systems,
and irreducible otherwise.
\end{itemize}
\end{definition}

\begin{remark}
${}$
\begin{itemize}
\item
Unless specified, we shall assume that \ssf$R$ spans \ssf$\apartment$\ssf.
\item
$\,\alpha^{\ssb\vee}\!=\ssb\frac2{|\alpha|^2}\,\alpha$
\ssf denotes the coroot corresponding to a root \ssf$\alpha$\ssf.
If \ssf$R$ is a root system,
then \ssf$R^\vee$ \ssb is again a root system.
\item
Most root systems are reduced, which means that

\noindent
{\rm(c)} the roots proportional to any root \ssf$\alpha$
are reduced to \ssf$\pm\ssf\alpha$\ssf.

\noindent
Otherwise the only possible alignment of roots is

\centerline{$
-\ssf2\ssf\alpha\ssf,-\ssf\alpha\ssf,+\ssf\alpha\ssf,+\ssf2\ssf\alpha\,.
$}\vspace{-.5mm}

\noindent
A root \ssf$\alpha$ \ssf is called
\begin{itemize}
\item[$\circ$]
indivisible if \,$\frac\alpha2$ is not a root,
\item[$\circ$]
non--multipliable if \,$2\ssf\alpha$ is not a root.
\end{itemize}
\item
We shall also consider
non--crytallographic reduced root systems \ssf$R$\ssf,
which satisfy \ssf{\rm(a)} and \ssf{\rm(c)}, but not necessarily \ssf{\rm(b)}.
\end{itemize}
\end{remark}

\begin{definition}
${}$
\begin{itemize}
\item
The connected components of
\vspace{.5mm}

\centerline{$
\{\,x\!\in\!\apartment\,|\,
\langle\ssf\alpha\ssf,x\ssf\rangle\ssb\ne\ssb0
\hspace{2mm}\forall\;\alpha\!\in\!R\,\}
$}\vspace{.5mm}

\noindent
are called Weyl chambers.
We choose any of them,
which is called positive and denoted by \,$\chamber$.
\ssf$R^{\vsf+}$ denotes the set of roots which are positive on \,$\chamber$.
\item
The Weyl or Coxeter group \,$W$ \ssb associated with \ssf$R$
is the finite subgroup of the or\-thog\-o\-nal group \,$\text{\rm O}(\apartment)$
generated by the root reflections \,$\{\,r_{\vsb\alpha}\,|\,\alpha\hspace{-.5mm}\in\!R\,\}$\ssf.
\end{itemize}
\end{definition}

\begin{remark}
${}$
\begin{itemize}
\item
The group \,$W$ \ssb acts simply transitively on the set of Weyl chambers.
\item
The longest element \,$w_{\ssf0}$ in \,$W$ \ssb
interchanges \,$\chamber$ and \,$-\ssf\chamber$.
\item
Every \,$x\hspace{-.5mm}\in\!\apartment$
\ssf belongs to the \ssf$W$\hspace{-.5mm}--\ssf orbit of a single
\,$x^+\hspace{-1mm}\in\!\overline{\smash{\chamber}\vphantom{X}}$.
\end{itemize}
\end{remark}

There are six classical families of irreducible root systems\,:
\begin{itemize}
\item
$\,\text{\rm A}_{\ssf n}$ ($n\hspace{-.5mm}\ge\!1$)\,:
\hspace{5.6mm}
$\apartment=\{\,x\!\in\!\R^{\ssf n\ssf+1}\,|\,
x_{\vsf0}\hspace{-.5mm}+\ssb x_1\!+\ssf\dots\ssf+\ssb x_{\vsf n}
\hspace{-.5mm}=\ssb0\,\}$

\hspace{25mm}
$R=\{\,e_{\vsf i}\ssb-\ssb e_j\,|\,
0\hspace{-.5mm}\le\ssb i\ssb\ne\ssb j\ssb\le\ssb n\,\}$

\hspace{25mm}
$\chamber\!=\{\,x\!\in\!\apartment\,|\,
x_{\vsf0}\hspace{-.5mm}>\ssb x_1\hspace{-.5mm}>\vsf\dots\vsf>\ssb x_{\vsf n}\,\}$

\hspace{25mm}
$W\hspace{-.5mm}=\ssf\text{S}_{\ssf n+1}$

\item
$\,\text{\rm B}_{\ssf n}$ ($n\hspace{-.5mm}\ge\hspace{-.5mm}2$)\,:
\hspace{5.55mm}
$\apartment=\Rn$

\hspace{25mm}
$R=\{\,\pm\,e_{\vsf i}\,|\,
1\hspace{-.5mm}\le\ssb i\ssb\le\ssb n\,\}
\cup\{\,\pm\,e_{\vsf i}\ssb\pm\ssb e_j\,|\,
1\hspace{-.5mm}\le\ssb i\ssb<\ssb j\ssb\le\ssb n\,\}$

\hspace{25mm}
$\chamber\!=\{\,x\!\in\!\R^n\,|\,
x_1\hspace{-.5mm}>\vsf\dots\vsf>\ssb x_{\vsf n}\hspace{-.5mm}>\ssb0\,\}$

\hspace{25mm}
$W\hspace{-.5mm}=\{\ssf\pm\ssf1\ssf\}^n\hspace{-.5mm}\rtimes\text{S}_{\ssf n}$
\vspace{.5mm}

\item
$\,\text{\rm C}_{\vsf n}$ ($n\hspace{-.5mm}\ge\hspace{-.5mm}2$)\,:
\hspace{5.65mm}
$\apartment=\Rn$

\hspace{25mm}
$R=\{\,\pm\,2\,e_{\vsf i}\,|\,
1\hspace{-.5mm}\le\ssb i\ssb\le\ssb n\,\}
\cup\{\,\pm\,e_{\vsf i}\ssb\pm\ssb e_j\,|\,
1\hspace{-.5mm}\le\ssb i\ssb<\ssb j\ssb\le\ssb n\,\}$

\hspace{25mm}
$\chamber\!=\{\,x\!\in\!\R^n\,|\,
x_1\hspace{-.5mm}>\vsf\dots\vsf>\ssb x_{\vsf n}\hspace{-.5mm}>\ssb0\,\}$

\hspace{25mm}
$W\hspace{-.5mm}=\{\ssf\pm\ssf1\ssf\}^n\hspace{-.5mm}\rtimes\text{S}_{\ssf n}$

\item
$\,\text{\rm BC}_{\vsf n}$ ($n\hspace{-.5mm}\ge\!1$)\,:
\hspace{2.95mm}
$\apartment=\Rn$

\hspace{25mm}
$R=\{\,\pm\,e_{\vsf i}\ssf,\pm\,2\,e_{\vsf i}\,|\,
1\hspace{-.5mm}\le\ssb i\ssb\le\ssb n\,\}
\cup\,\{\,\pm\,e_{\vsf i}\ssb\pm\ssb e_j\,|\,
1\hspace{-.5mm}\le\ssb i\ssb<\ssb j\ssb\le\ssb n\,\}$

\hspace{25mm}
$\chamber\!=\{\,x\!\in\!\R^n\,|\,
x_1\hspace{-.5mm}>\vsf\dots\vsf>\ssb x_{\vsf n}\hspace{-.5mm}>\ssb0\,\}$

\hspace{25mm}
$W\hspace{-.5mm}=\{\ssf\pm\ssf1\ssf\}^n\hspace{-.5mm}\rtimes\text{S}_{\ssf n}$

\item
$\,\text{\rm D}_n$ ($n\hspace{-.5mm}\ge\hspace{-.5mm}3$)\,:
\hspace{5.6mm}
$\apartment=\R^n$

\hspace{25mm}
$R=\{\,\pm\,e_{\vsf i}\ssb\pm\ssb e_j\,|\,
1\hspace{-.5mm}\le\ssb i\ssb<\ssb j\ssb\le\ssb n\,\}$

\hspace{25mm}
$\chamber\!=\{\,x\!\in\!\R^n\,|\,
x_1\hspace{-.5mm}>\vsf\dots\vsf>\ssb|\ssf x_{\vsf n}|\,\}$

\hspace{25mm}
$W\hspace{-.5mm}=\{\,\epsilon\hspace{-.5mm}\in\!\{\ssf\pm\ssf1\ssf\}^n\,
|\,\epsilon_1\ssb\dots\ssf\epsilon_{\vsf n}\!=\!1\,\}\hspace{-.5mm}
\rtimes\text{S}_{\ssf n}$

\item
$\,\text{I}_{\vsf 2}(m)$ ($m\hspace{-.5mm}\ge\hspace{-.5mm}3$)\,:
\hspace{-.55mm}
$\apartment=\C$

\hspace{25mm}
$R=\{\,\smash{e^{\hspace{.5mm}i\hspace{.5mm}\pi\ssf\frac jm}}\,|
\,0\ssb\le\ssb j\ssb<\ssb2\hspace{.5mm}m\,\}$

\hspace{25mm}
$\chamber\!=\{\,z\!\in\!\C^*\,|
\,\bigl(\frac12\!-\!\frac1m\bigr)\hspace{.5mm}\pi\ssb<\ssb\arg z\ssb<\ssb\frac\pi2\,\}$

\hspace{25mm}
$W\hspace{-.5mm}=(\vsf\Z\ssf/m\hspace{.5mm}\Z\vsf)\!\rtimes\!(\vsf\Z\ssf/\ssf2\ssf\Z\vsf)$ \,(dihedral group)
\end{itemize}

The full list of irreducible root systems (crystallographic or reduced)
includes in addition a finite number of exceptional cases\,:

\centerline{$
\,\text{\rm E}_{\ssf6}\ssf,
\,\text{\rm E}_{\ssf7}\ssf,
\,\text{\rm E}_{\ssf8}\ssf,
\,\text{\rm F}_{\hspace{-.4mm}4}\ssf,
\,\text{\rm G}_{\ssf2}\ssf,
\,\text{\rm H}_{\ssf3}\ssf,
\,\text{\rm H}_{\ssf4}\ssf.
$}

\begin{remark}
In the list above,
${}$
\begin{itemize}
\item
the non crystallographic root systems are
\vspace{.5mm}

\centerline{
$\text{\rm H}_{\ssf3}$\ssf,
$\text{\rm H}_{\ssf4}$
and \,$\text{\rm I}_{\ssf2}(m)$
\ssf with \,$\begin{cases}
\,m\ssb=\ssb5\ssf,\\
\,m\ssb\ge\ssb7,\\
\end{cases}$
}\vspace{1mm}

\item
all root systems are reduced, with the exception of \,$\text{\rm BC}_{\ssf n}$\ssf,
\item
there are some redundancies in low dimension\,$:$
\vspace{1mm}

\centerline{$\begin{cases}
\;\text{\rm A}_1\!\times\hspace{-.5mm}\text{\rm A}_1\!
\approx\hspace{-.5mm}\text{\rm D}_{\ssf2}\hspace{-.5mm}
\approx\hspace{-.5mm}\text{\rm I}_{\ssf2}(2)\\
\;\text{\rm B}_{\ssf2}\hspace{-.5mm}
\approx\hspace{-.5mm}\text{\rm C}_{\ssf2}\hspace{-.5mm}
\approx\hspace{-.5mm}\text{\rm I}_{\ssf2}(4)
&\text{$($up to the root length\ssf$)$}\\
\;\text{\rm A}_{\ssf2}\hspace{-.5mm}
\approx\hspace{-.5mm}\text{\rm I}_{\ssf2}(3)\\
\;\text{\rm G}_{\ssf2}\hspace{-.5mm}
\approx\hspace{-.5mm}\text{\rm I}_{\ssf2}(6)
&\text{$($up to the root length\ssf$)$}\\
\end{cases}$}

\end{itemize}
\end{remark}

\begin{figure}[ht]
\vspace{-10mm}
\psfrag{A1A1}[c]{$\text{\rm A}_{\vsf1}\!\ssb\times\!\text{\rm A}_{\vsf1}$}
\psfrag{BC1BC1}[c]{$\text{\rm BC}_{\vsf1}\times\text{\rm BC}_{\vsf1}$}
\centerline{\hspace{5mm}
\includegraphics[height=52mm]{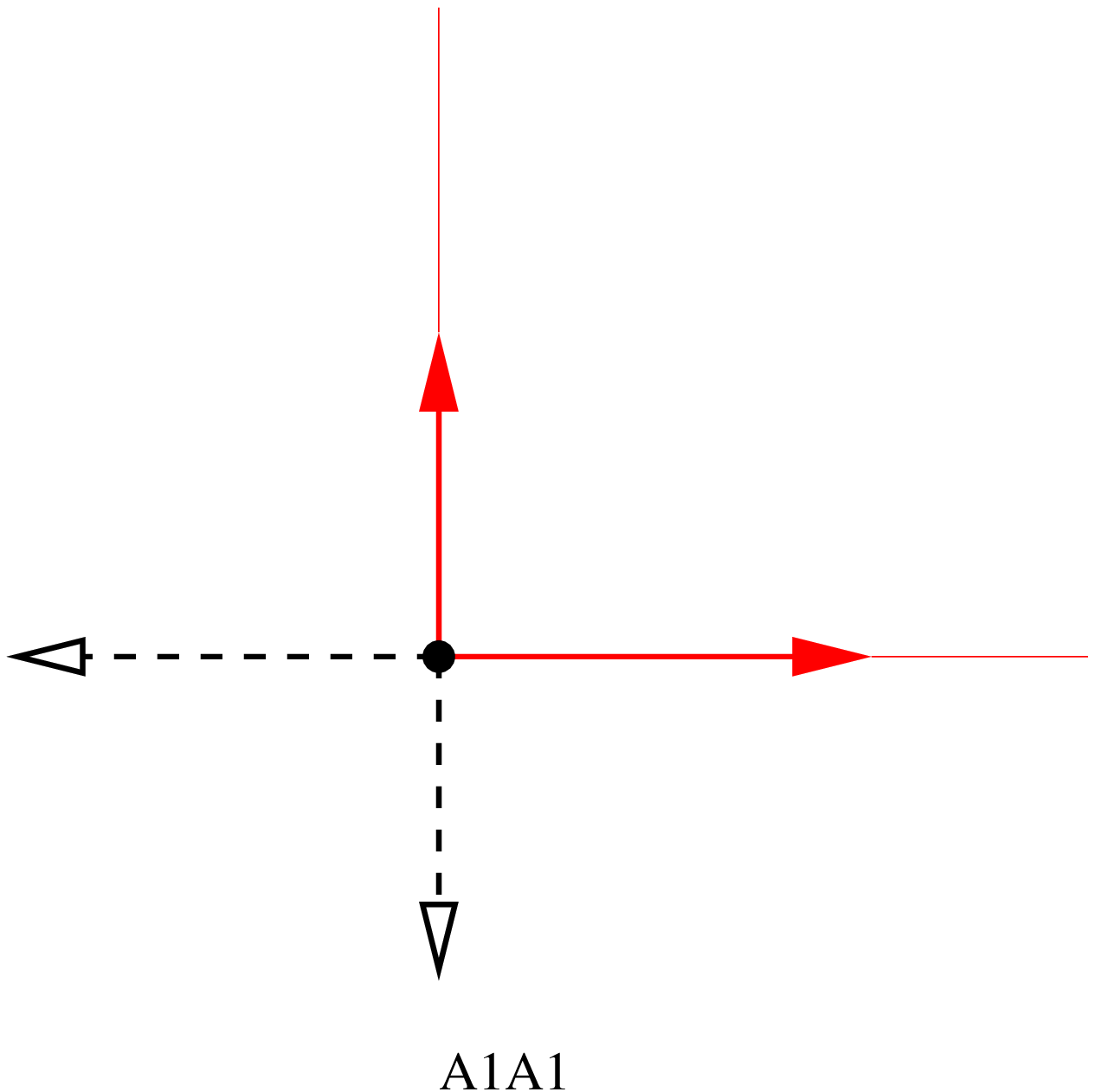}
\hspace{14mm}
\includegraphics[height=45mm]{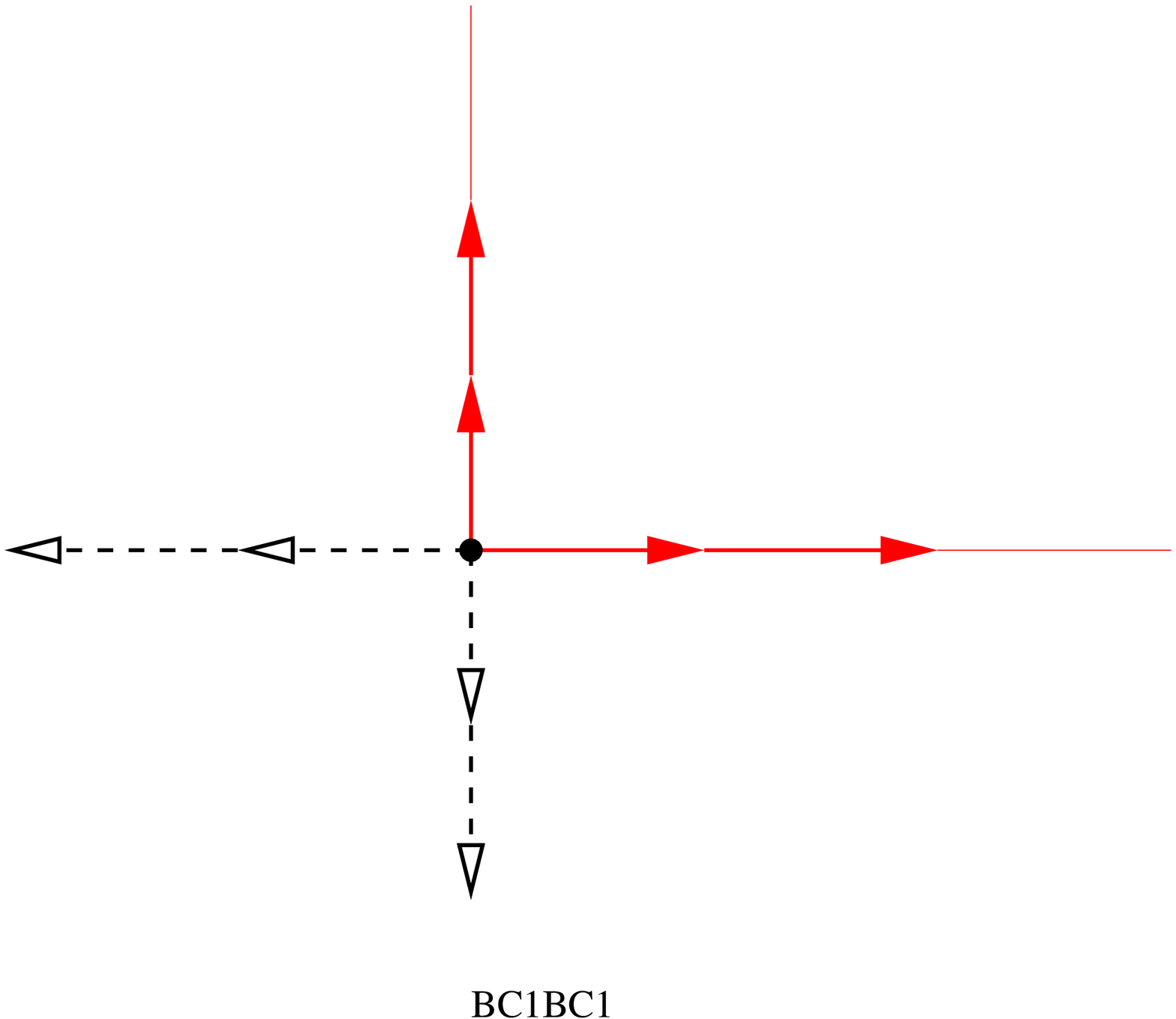}
\hspace{0mm}}
\vspace{5mm}
\psfrag{A2}[c]{$\text{\rm A}_{\ssf2}$}
\psfrag{B2}[c]{$\text{\rm B}_{\ssf2}$}
\centerline{\hspace{4mm}
\includegraphics[height=49mm]{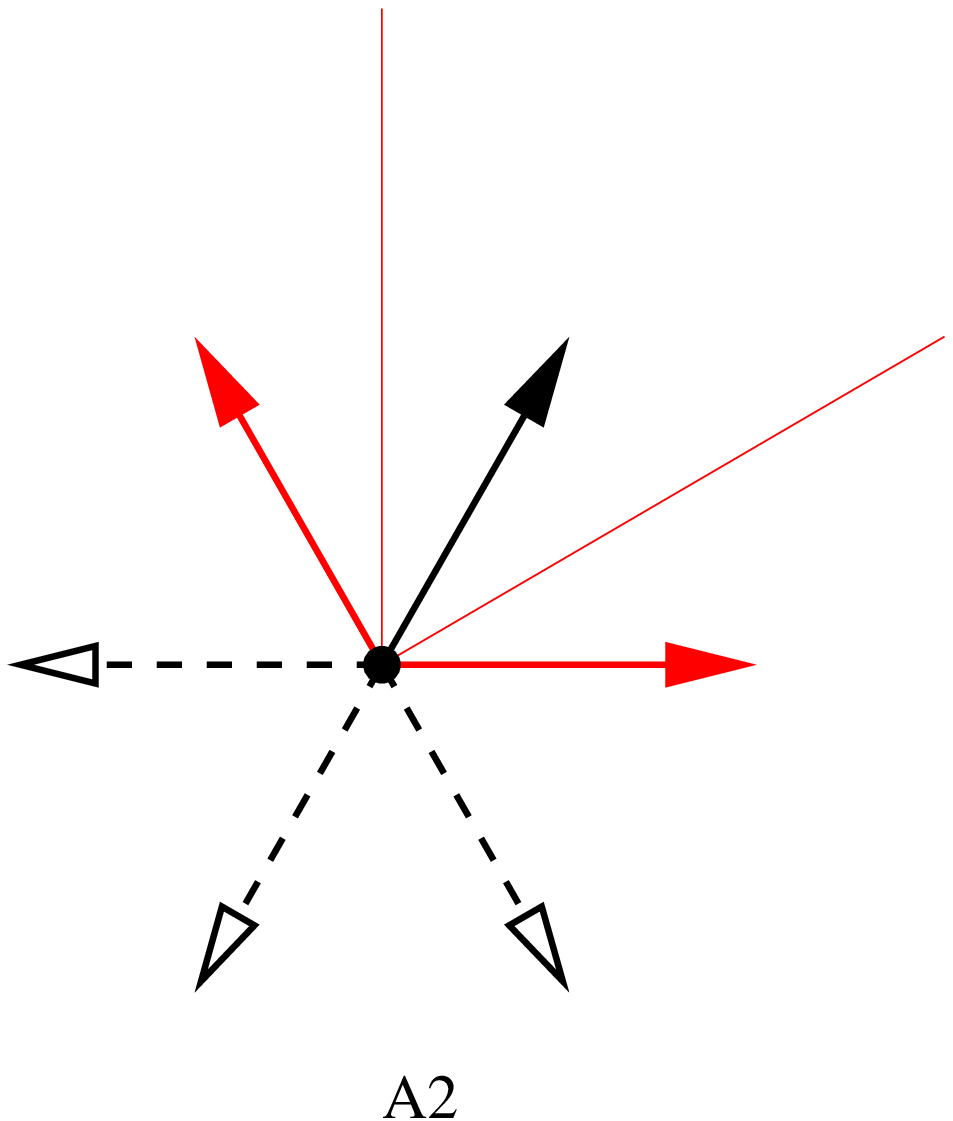}
\hspace{26mm}
\includegraphics[height=50mm]{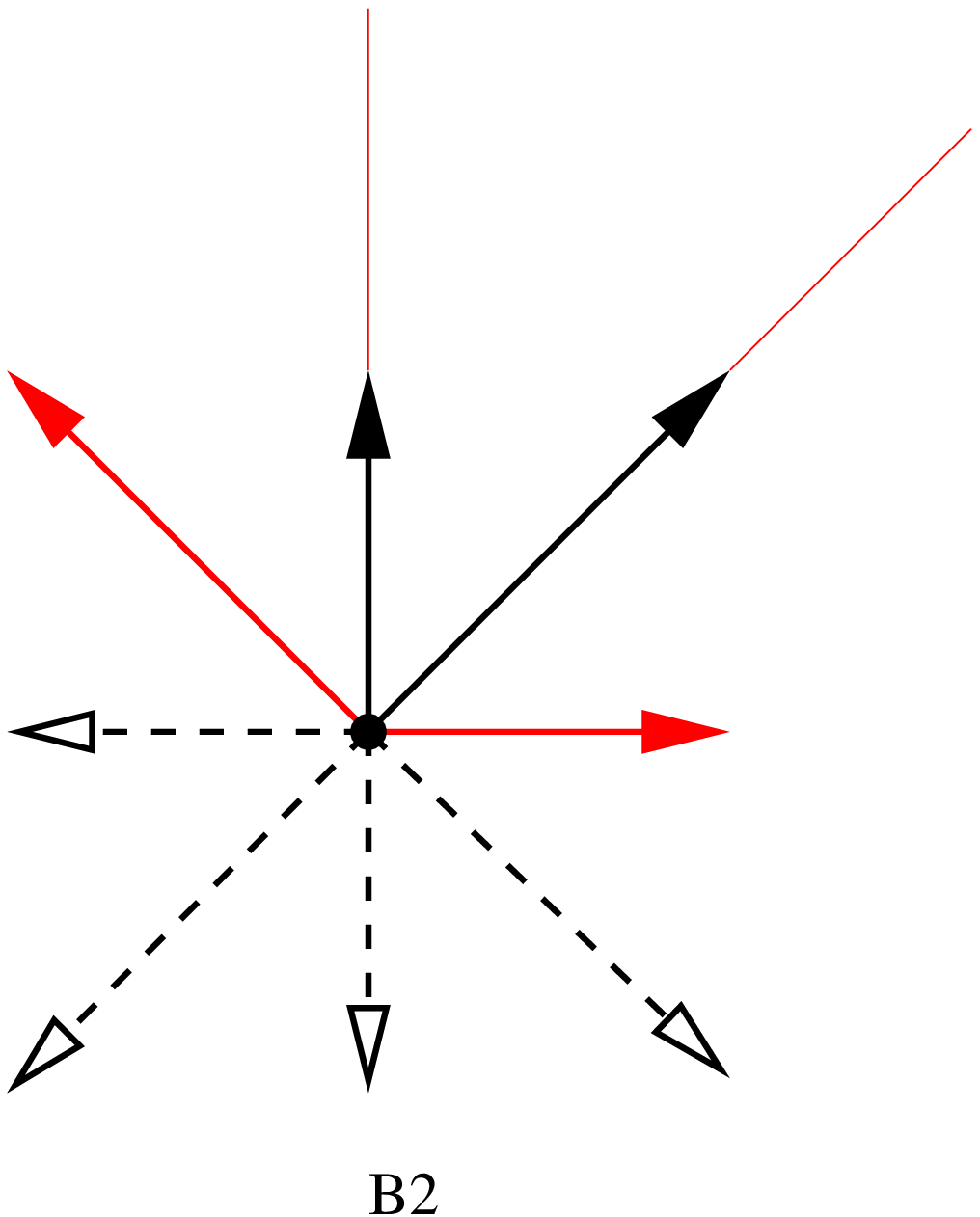}
\hspace{0mm}}\vspace{5mm}
\psfrag{C2}[c]{$\text{\rm C}_{\ssf2}$}
\psfrag{BC2}[c]{$\text{\rm BC}_{\ssf2}$}
\centerline{\hspace{0mm}
\includegraphics[height=56mm]{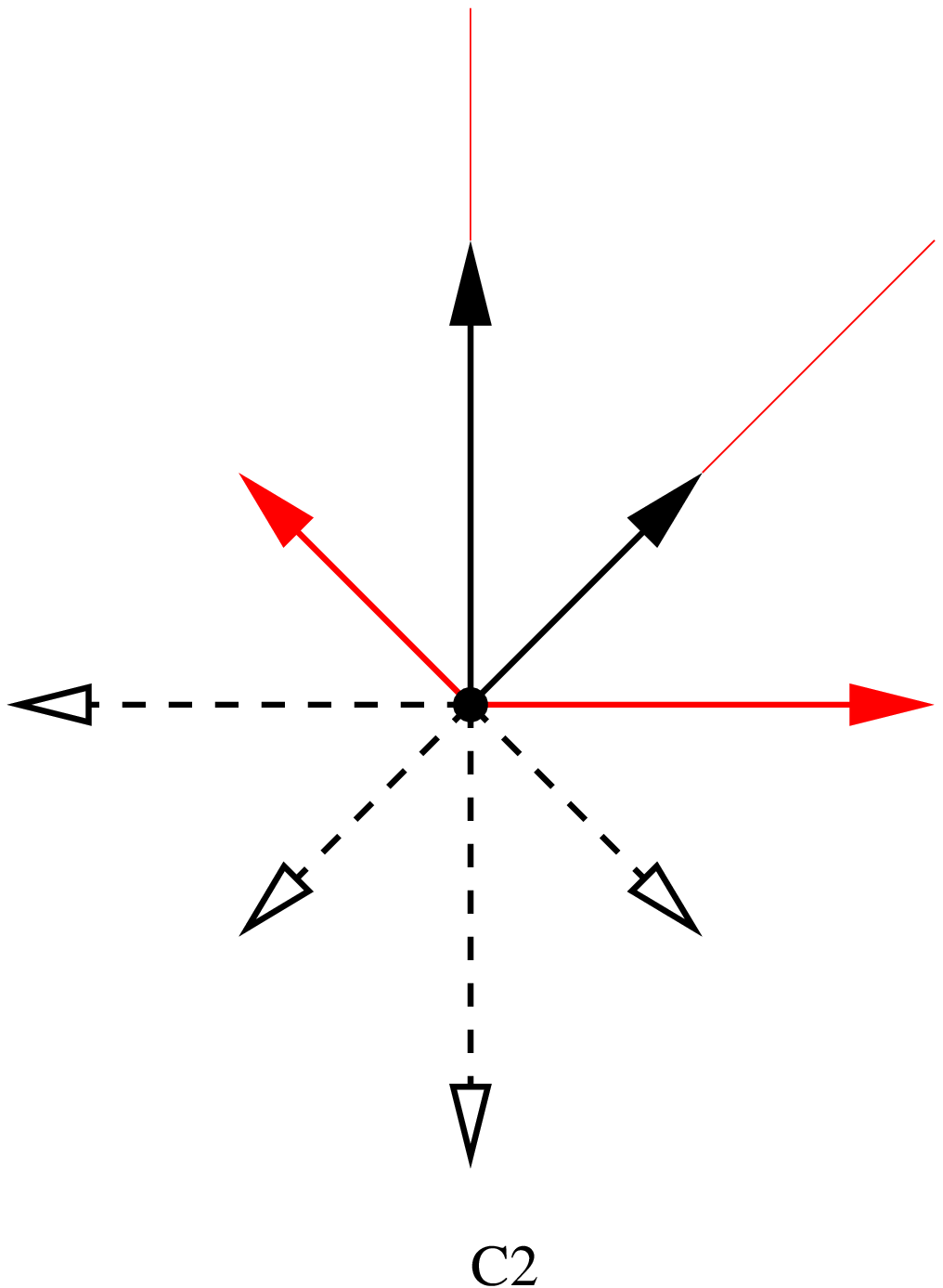}
\hspace{23mm}
\includegraphics[height=55mm]{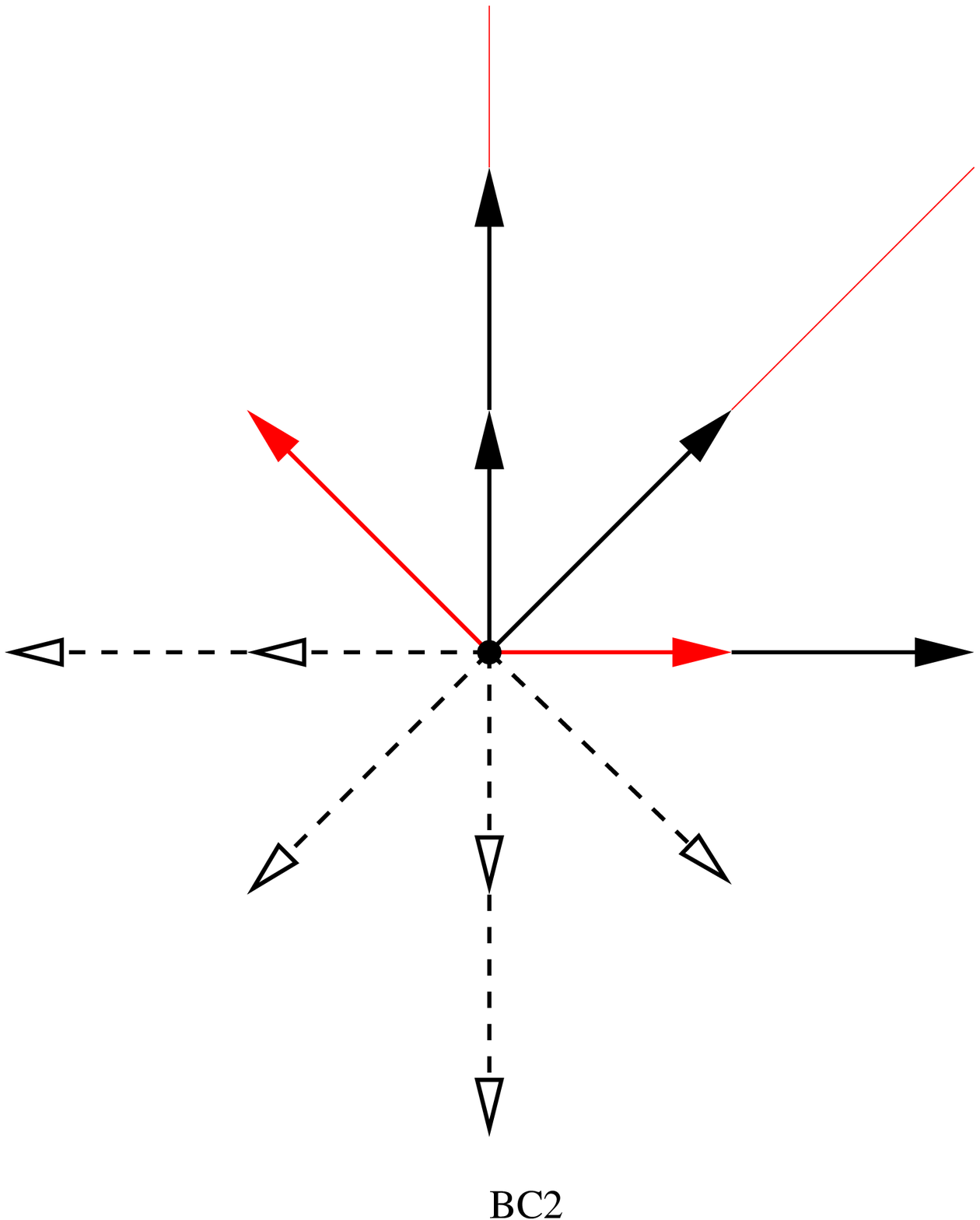}
\hspace{3mm}}\vspace{5mm}
\psfrag{G2}[c]{$\text{\rm G}_{\ssf2}$}
\psfrag{I25}[c]{$\text{\rm I}_{\ssf2}(5)$}
\centerline{\hspace{0mm}
\includegraphics[height=59.5mm]{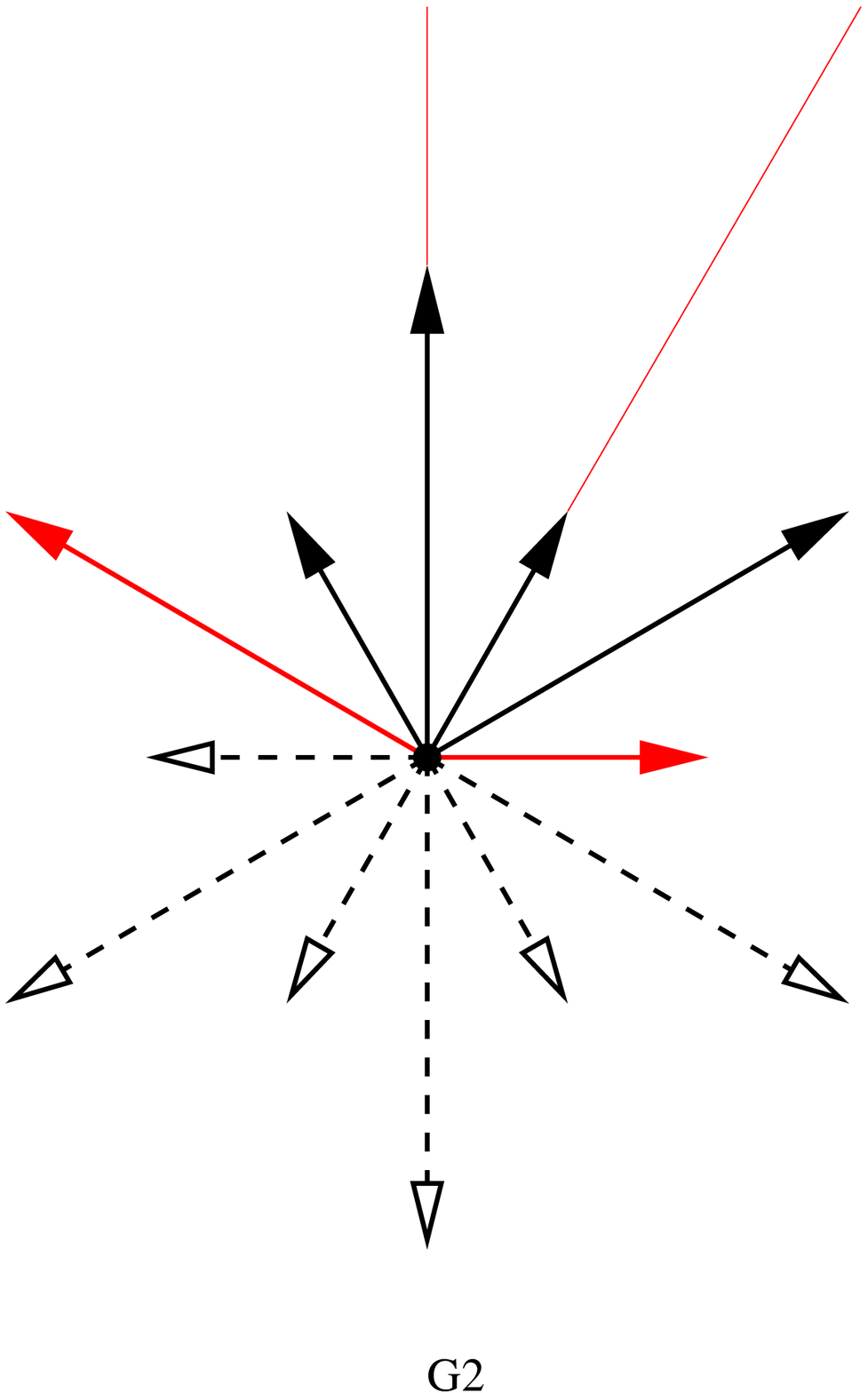}
\hspace{25mm}
\includegraphics[height=60mm]{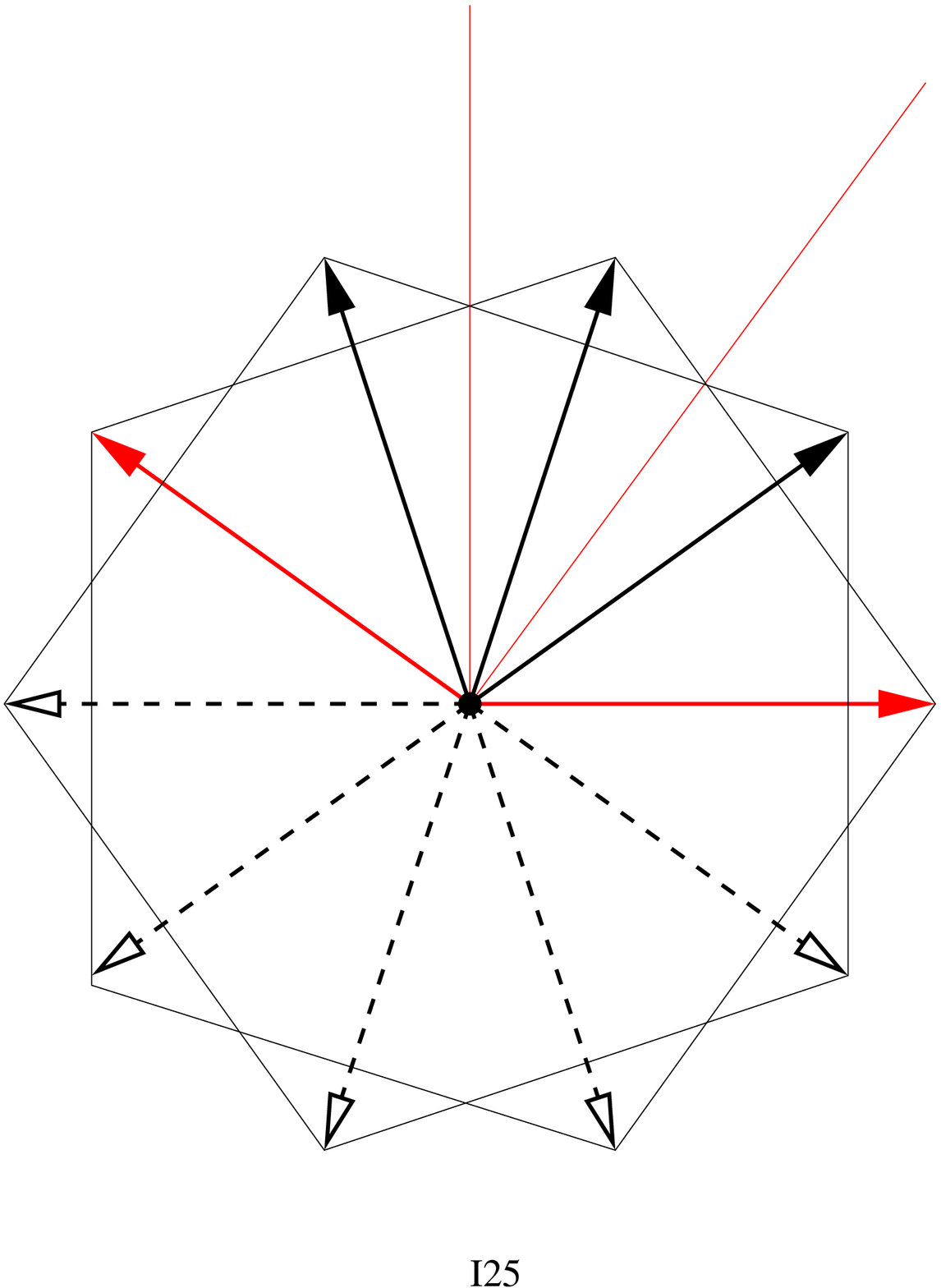}
\hspace{1mm}}
\caption{2\ssf--\ssf dimensional root systems}
\label{2DimensionalRootSystems}
\end{figure}

The 2--dimensional root systems (crystallographic or reduced)
are depicted in Figure \ref{2DimensionalRootSystems}.

\begin{definition}
A multiplicity is a \ssf$W$\hspace{-.5mm}--\ssf invariant function
\,$k\ssb:\ssb R\ssb\longrightarrow\ssb\C$\ssf.
\end{definition}

\begin{remark}
${}$
\begin{itemize}
\item
In Dunkl theory,
one assumes most of the time that \,$k\ssb\ge\ssb0$\ssf.
\item
Assume that \ssf$R$ is crystallographic and irreducible.
Then two roots belong to the same \ssf$W$\hspace{-.5mm}--\ssf orbit
if and only if they have the same length.
Thus \ssf$k$ takes at most three values.
In the non crystallographic case,
there are one or two \ssf$W$\hspace{-.5mm}--\ssf orbits in \ssf$R$\ssf.
Specifically, by resuming the classification of root systems,
\ssf$k$ takes
\begin{itemize}
\item[$\circ$]
$\ssf1$ value in the following cases\,{\rm:}
\vspace{.25mm}

\centerline{$
\text{\rm A}_{\ssf n}\ssf,
\,\text{\rm D}_{\ssf n}\ssf,
\,\text{\rm E}_{\ssf6}\ssf,
\,\text{\rm E}_{\ssf7}\ssf,
\,\text{\rm E}_{\ssf8}\ssf,
\,\text{\rm H}_{\ssf3}\ssf,
\,\text{\rm H}_{\ssf4}\ssf,
\,\text{\rm I}_{\ssf2}(m)
\text{ with \hspace{.5mm}$m$ odd}\ssf,
$}

\item[$\circ$]
$\ssf2$ values in the following cases\,{\rm:}
\vspace{.25mm}

\centerline{$
\text{\rm B}_{\ssf n}\ssf,
\,\text{\rm C}_{\ssf n}\ssf,
\,\text{\rm F}_{\hspace{-.4mm}4}\ssf,
\,\text{\rm G}_{\vsf2}\ssf,
\,\text{\rm I}_{\ssf2}(m)
\text{ with \hspace{.5mm}$m$ even}\ssf,
$}

\item[$\circ$]
$\ssf3$ values in the case of \,$\text{\rm BC}_{\ssf n}$\ssf.
\end{itemize}
\end{itemize}
\end{remark}

\newpage

\end{document}